\documentclass[12pt]{amsart}
\usepackage{amsmath, amssymb}
\newfont {\cyr} {wncyr10}
\renewcommand{\labelenumi}{{(\roman{enumi})}}
\usepackage{rotating}

\usepackage{amsfonts}
\usepackage{amsmath}
\usepackage{amssymb}
\usepackage{setspace}
\usepackage{multicol}
\newtheorem{theorem}{Theorem}[section]
\newtheorem{lemma}[theorem]{Lemma}\newtheorem{notation}[theorem]{Notation}
\newtheorem{proposition}[theorem]{Proposition}
\newtheorem{corollary}[theorem]{Corollary}
\newtheorem{definition}[theorem]{Definition}
\newtheorem{hypothesis}[theorem]{Hypothesis}

\newcounter{claim}[theorem]

\newcounter{cclaim}[theorem]

\def \udot {{}^{\textstyle .}}

\newcommand{\E}{\mathrm{E}}\newcommand{\SU}{\mathrm{SU}}
\newcommand{\F}{\mathrm{F}}\newcommand{\A}{\mathrm{A}}\newcommand{\B}{\mathrm{B}}\newcommand{\M}{\mathcal{M}}
\newcommand{\G}{\mathrm{G}}
\newcommand{\D}{\mathcal{D}}
\newcommand{\C}{\mathcal{C}}
\newcommand{\Q}{\mathrm{Q}}

\newcommand{\Aut}{\mathrm{Aut}}

\newcommand{\Out}{\mathrm{Out}}

\newcommand{\Syl}{\mathrm{Syl}}\newcommand{\syl}{\mathrm{Syl}}

\newcommand{\GF}{\mathrm{GF}}
\newcommand{\GL}{\mathrm{GL}}
\newcommand{\Sp}{\mathrm{Sp}}
\newcommand{\SL}{\mathrm{SL}}

\newcommand{\PGL}{\mathrm{PGL}}
\newcommand{\PSL}{\mathrm{PSL}}\newcommand{\PSp}{\mathrm{PSp}}
\newcommand{\Sym}{\mathrm{Sym}}
\newcommand{\Alt}{\mathrm{Alt}}
\newcommand{\Dih}{\mathrm{Dih}}
\newcommand{\SDih}{\mathrm{SDih}}
\newcommand{\Frob}{\mathrm{Frob}}

\def \OO {\hbox {\rm O}}

\def \Z {\mathbb Z}
\def \GU {\mbox {\rm GU}}

\def \syl {\hbox {\rm Syl}}\def \Syl {\hbox {\rm Syl}}

\def \ov {\overline}

\def \wt {\widetilde}

\def \Aut{ \mathrm {Aut}}

\def \Out{\mbox {\rm Out}}

\def \Fi {\mbox {\rm Fi}}

\def \J{\mbox {\rm J}}

\def \Ly{\mbox {\rm Ly}}

\def \B{\mbox {\rm B}}

\def \M{\mbox {\rm M}}

\def \HN{\mbox {\rm HN}}

\def \HS{\mbox {\rm HS}}

\def \Th{\mbox {\rm Th}}

\def \ON {\mbox {\rm O'N}}

\def \Co {\mbox {\rm Co}}

\def \Ru {\mbox {\rm Ru}}

\def \Suz{\mbox {\rm Suz}}

\def \McL{\mbox {\rm McL}}

\def \He {\mbox {\rm He}}\def \PSU {\mbox {\rm PSU}}

\begin{document}
\renewcommand{\labelenumi}{(\roman{enumi})}

\title  {Strongly $p$-embedded subgroups}
 \author{Chris Parker}
  \author{Gernot Stroth}

\address{Chris Parker\\
School of Mathematics\\
University of Birmingham\\
Edgbaston\\
Birmingham B15 2TT\\
United Kingdom} \email{c.w.parker@bham.ac.uk}

\address{Gernot Stroth\\
Institut f\"ur Mathematik\\ Universit\"at Halle - Wittenberg\\
Theordor Lieser Str. 5\\ 06099 Halle\\ Germany}
\email{gernot.stroth@mathematik.uni-halle.de}

\email {}

\date{\today}

\maketitle \pagestyle{myheadings}

\markright{{\sc }} \markleft{{\sc Chris Parker and Gernot Stroth}}

\section{Introduction}

In this paper we study finite groups which possess a strongly $p$-embedded subgroup for some prime
$p$. Suppose that $p$ is a prime. A subgroup $H$ of the finite group $G$ is said to be
\emph{strongly $p$-embedded in $G$} if the following two conditions hold.
\begin{enumerate}
\item[(i)] $H < G$ and $p$ divides $|H|$; and
\item[(ii)] if $g \in G\setminus H$, then $p$ does not divide $|H \cap H^g|$.
\end{enumerate}
One of the most important properties of strongly $p$-embedded subgroups is that  $N_G(X)\le H$ for
any non-trivial $p$-subgroup $X$ of $H$.

 Groups with a strongly $2$-embedded subgroup
have been classified by Bender \cite{Bender} and Suzuki \cite{Suzuki} and their classification
forms a pedestal upon which the classification of the finite simple groups stands. The
classification of groups with a strongly $2$-embedded subgroup states that if $G$ is a finite group
with a strongly $2$-embedded subgroup, then $O^{2'}(G/O(G))$ is a simple rank $1$ Lie type group
defined in characteristic 2 or $G$ has quaternion or cyclic Sylow 2-subgroups. Of course rank $1$
Lie type groups in characteristic $2$ are the building blocks of the groups of Lie type in
characteristic $2$. When dealing with groups defined in characteristic $p$ for odd $p$, strongly
$p$-embedded subgroups play an equally influential role.

Assume from here on that $p$ is an odd prime and $G$ is a finite group. If $G$ has cyclic Sylow
$p$-subgroup $P$, then $N_G(\Omega_1(P))$ is strongly $p$-embedded in $H$. There is thus no
prospect of listing all such groups. However, if $m_p(G) \ge 2$, then the almost simple groups with
a strongly $p$-embedded subgroup are known as a corollary to the classification of the finite
simple groups. They include the rank $1$ Lie type groups defined in characteristic $p$, and there
are only additional examples if $p \le 11$. See Proposition~\ref{SE-p2} for a complete list. When
we examine the list of simple groups  with a strongly $p$-embedded subgroup $H$, we see that either
$H$ is a $p$-local subgroup or that $F^*(H) \cong \Alt(p) \times \Alt(p)$ with $p \ge 5$ or
$\Omega_8^+(2)$ with $p = 5$. Thus the fact that $H$ is strongly $p$-embedded severely restricts
its structure. One major application of the results of this paper will be to the investigation of
groups of local characteristic $p$ orchestrated by Meierfrankenfeld, Stellmacher and Stroth
\cite{MSS}. In the described application, we have a subgroup $H$ of $G$ generated by all the
normalizers of non-trivial $p$-subgroups contained in a fixed Sylow $p$-subgroup of $G$. Generally
speaking under such circumstances $H$ will be a Lie type group defined in characteristic $p$. We
would like to assert that $H=G$. Assuming that this is not the case, work in progress by Salarian
and Stroth will show  that $H$ is strongly $p$-embedded in $G$. The conclusion of our second
theorem, Theorem~\ref{MainTheorem1}, is that under mild restrictions on the structure of $H$,  $H$
cannot in fact be strongly $p$-embedded. Recall that a \emph{$\mathcal K$-group} is a group in
which every composition factor is from the list of ``known" simple groups. That is, every simple
section is either a cyclic group of prime order, an alternating group, a group of Lie type or one
of the twenty six sporadic simple groups. A group is \emph{$\mathcal K$-proper} if all its proper
subgroups are $\mathcal K$-groups. Obviously $\mathcal K$-proper groups are the focus of attention
in the proof of the classification of the finite simple groups.

\begin{theorem}\label{MainTheorem0} Suppose that $G$ is a finite group, $p$ is an odd prime and that $H$ is a strongly $p$-embedded
subgroup of $G$ such that $H \cap K$ has even order for all non-trivial normal subgroups $K$ of
$G$. Assume that $F^*(H)=O_p(H)$ and $m_p(C_H(t)) \ge 2$ for every involution $t$ in $H$. If
$N_G(T)$ is a $\mathcal K$-group for all non-trivial $2$-subgroups $T$ of $G$, then there  exists
$n \ge 2$ such that either
 $F^*(G)\cong \PSU_3(p^n)$ or, $p=3$ and $F^*(G) \cong {}^2\G_2(3^{2n-1})$.
\end{theorem}

\begin{theorem}\label{MainTheorem1} Suppose that $G$ is a finite
group, $p$ is an odd prime and that $H$ is a strongly $p$-embedded subgroup of $G$. Assume that
$O_{p'}(H)=1$ and that $m_p(C_H(t)) \ge 2$ for every involution $t$ of $H$. If $N_G(T)$ is a
$\mathcal K$-group for all non-trivial $2$-subgroups $T$ of $G$ and $H$ is a $\mathcal K$-group,
then $F^\ast(H)=O_p(H)$.
\end{theorem}

We choose to compartmentalize our  theorems in this way so that we can be precise about the
$\mathcal K$-group hypothesis we are invoking. For example in the expected application of this work
to the project to classify the finite groups of local characteristic $p$ as explained above, $H$
will be a finite quasisimple $\mathcal K$-group and this will be known independently with just a
$\mathcal K$-group assumption on the normalizers of non-trivial $p$-subgroups of $G$. Our two
theorems combine to give the following theorem which is the main result of the paper.

\begin{theorem}\label{MainTheorem} Suppose that $G$ is a finite $\mathcal K$-proper
group, $p$ is an odd prime and that $H$ is a strongly $p$-embedded subgroup of $G$ such that $H
\cap K$ is of even order for all non-trivial normal subgroups $K$ of $G$. Assume that $O_{p'}(H)=1$
and that $m_p(C_H(t)) \ge 2$ for every involution $t$ of $H$.  Then there  exists $n \ge 2$ such
that either
 $F^*(G)\cong \PSU_3(p^n)$ or, $p=3$ and $F^*(G) \cong {}^2\G_2(3^{2n-1})$.
\end{theorem}

The hypothesis in Theorems~\ref{MainTheorem0} and \ref{MainTheorem} that $|H \cap K|$ is even for
all normal subgroups $K$ of $G$ is there to guarantee
 that  we cannot transfer all the elements of order $2$ from $H$ and finish up with a
configuration where the strongly $p$-embedded subgroup has odd order. In particular, we do not want
our group $G$ to have $G > O^2(G)$ and $H \cap O^2(G)$  a strongly $p$-embedded subgroup of
$O^2(G)$ of odd order. For example, taking $G \cong \PGL_2(27)$, we see that  $G$ has a strongly
$3$-embedded subgroup of even order, but the normal subgroup of index $2$ does not.

Recall that our primary application of our theorem is to the project to classify groups of local
characteristic $p$. With this in mind, we have the following corollary to
Theorem~\ref{MainTheorem}.

\begin{corollary}\label{LiepCor}  Suppose that $p$ is an odd prime and $H$ is a strongly $p$-embedded subgroup of $G$. Then $F^*(H)$ is not a Lie type group
defined in characteristic $p$ of rank at least $3$.
\end{corollary}

In fact, if $G$ and $H$ are as in Corollary~\ref{LiepCor}, then Theorem~\ref{MainTheorem} can be
used to show that if $F^*(H)$ is a Lie type group in characteristic $p$ of rank $2$, then the only
possibility is that  $F^*(H) \cong \PSL_3(p)$ or $\PSp_4(p)$. In work in preparation by the
authors, it is shown that the latter case cannot happen. To eliminate the possibility that $F^*(H)
\cong \PSL_3(p)$, it seems that different techniques need to be developed.

We now proceed to describe the contents of the paper and at the same time give an outline of the
proof of our theorems. In Section~\ref{prel}, we present the preliminary results that we shall call
upon throughout the paper. We begin with our main characterization theorems. These are the results
which formally identify the groups in Theorem~\ref{MainTheorem0}. Thus  we state Aschbacher's
Classical Involution Theorem (Theorem~\ref{clasinv}) which will be used to identify the unitary
groups and the theorem which is the culmination of work by Walter, Bombieri, and Thompson
(Theorem~\ref{E8Sylow}) to recognize the Ree groups. We then move on to results about almost simple
$\mathcal K$-groups. The most  prominent and important of these is Proposition~\ref{SE-p2} which
describes the structure of the almost simple $\mathcal K$-groups of $p$-rank at least $2$ which
possess a strongly $p$-embedded subgroup. This is followed by Proposition~\ref{Hstru1} which
catalogues the structure of a strongly $p$-embedded subgroup in each of the groups listed in
Proposition~\ref{SE-p2}. Next  a series of corollaries to the propositions are presented, perhaps
the most useful  being Corollary~\ref{OrredH}. We then collect some results related to the Thompson
Transfer Lemma and results which limit the structure of $2$-groups which admit certain types of
automorphisms, see, for example, Lemmas~\ref{2group} and \ref{AutosE}.

Our investigation of groups with a strongly $p$-embedded subgroups starts in earnest in
Section~\ref{basic}. We assume that $G$ is a group and that $H$ is a strongly $p$-embedded subgroup
of $G$ and prove some basic properties about such configurations. In particular, easy, but key
facts that are used throughout the work are established such as if $K$ is a subgroup of $G$ which
is not contained in $H$ and $H \cap K$ has order divisible by $p$, then $H\cap K$ is strongly
$p$-embedded in $K$ and, of particular importance, if $m_p(H\cap K) \ge 2$, then $O_{p'}(K) \le H$
and $F^*(C_G(t)/O_{p'}(C_G(t)))$ is an almost simple group of $p$-rank at least $2$ containing a
strongly $p$-embedded subgroup. This is a simple consequence of the properties of strongly
$p$-embedded subgroups and coprime action. A further consequence of these elementary lemmas is that
$F^*(G)$ is a non-abelian simple group. Now the hypothesis in the main theorems that $H \cap K$ has
even order immediately implies that $H \cap F^*(G)$ has even order. Thus for our investigations we
may as well suppose that $G=O^2(G)$ and so, in Section~\ref{sec4}, where  we introduce the main
hypothesis Hypothesis~\ref{hypH}, this is included. In Section~\ref{sec4} we start our
investigation of groups $G$ with a strongly $p$-embedded subgroup $H$ such that $O_{p'}(H)=1$,
$m_p(C_H(t)) \ge 2$ for all involutions $t \in H$, and $G=O^2(G)$. We also impose our $\mathcal
K$-group hypothesis, that $N_G(T)$ is a $\mathcal K$-group for all $2$-groups $T$ of $G$. We may
also suppose that $G$ is not a $\mathcal K$-group.  Thus $H$ is not strongly $2$-embedded in $G$
and so, as $H$ has even order, there is an involution $t \in H$ with $C_G(t) \not \le H$. Since
$m_p(C_H(t)) \ge 2$, we have $O_{p'}(C_H(t)) \le H$ and $C_H(t)$ is strongly $p$-embedded in
$C_G(t)$. This allows us to use the $\mathcal K$-group hypothesis to get hold of the structure of $
C_G(t)/O_{p'}(C_G(t))$. One stark consequence of the information that we obtain in this section is
that for any involution $t \in H$, $C_G(t) \not \le H$ and the pertinent structural information
about $C_G(t)$ is listed in Lemma~\ref{Ctnot-in}. We close Section~\ref{sec4} with an immediate
application of Lemma~\ref{Ctnot-in} which states that if $E(H) \neq 1$, then $E(H)$ is quasisimple.
This is Lemma~\ref{Eqs}.

Sections~\ref{sec5}, \ref{sec7} and \ref{sec6},  study groups satisfying Hypothesis~\ref{hypH}
which have  $Q=F^*(H)=O_p(H)$.  Since $H$ is strongly $p$-embedded in $G$, we then have $H=
N_G(Q)$. Of course we may also assume that $G$ does not have a classical involution. Our objective
then in these three sections is to prove that  $G \cong {}^2\G_2(3^{2n-1})$ for some $n \ge 2$. To
do this we seek to establish the hypothesis of Theorem~\ref{E8Sylow}. So we need to show that $G$
has Sylow $2$-subgroups which have order $8$ and an involution $t$ such that $C_G(t)$ contains a
normal subgroup isomorphic to $\PSL_2(p^a)$ for some $a \ge 2$. We start this work in
Section~\ref{sec5}. For $t$ an involution in $H$, $t$ acts on $Q$ and, of course, either
$C_Q(t)\neq 1$ or $C_Q(t) =1$. In the former case Lemma~\ref{thereisacomponent} shows that there is
a unique component $L_t$ contained in $C_G(t)$ which is normal and has order divisible by $p$.  In
the latter case, we can only establish that $C_G(t)$ has a normal $2$-component
(Lemma~\ref{component}) with the same properties. In both cases, as by hypothesis $C_G(t)$ is a
$\mathcal K$-group, we may use Proposition~\ref{SE-p2} to describe the isomorphism types of
$L_t/O(L_t)$. We then make a careful choice of $t$. We assume that $|C_Q(t)|$ is maximal from among
all involutions $t \in H$ with $C_Q(t)\neq 1$ or that $C_Q(t)=1$. This choice leads to
Lemma~\ref{cyclic} which states that a Sylow $2$-subgroup of $O_{p'}(C_G(t))$ is cyclic and that
$O_{p'}(C_G(t))$ has a normal $2$-complement. In Section~\ref{sec7} our objective is to prove that
$L_t/O_{p'}(L_t) \cong \PSL_2(p^a)$ for some $a \ge 2$. The exact conclusion being posted in
Theorem~\ref{L2pthm}. This section deals with all the other possibilities for $L_t/O_{p'}(L_t)$ the
most resilient case arising when $H$ is strongly $3$-embedded and $L_t/O(L_t)$ is a covering group
of $\PSL_3(4)$. Section~\ref{sec6} is the most technical of the paper. We begin with an involution
$t$ in $H$ such that $L_t/O_{p'}(L_t)\cong \PSL_2(p^a)$ and $O_{p'}(C_G(t))$ has cyclic Sylow
$2$-subgroups. Arguments about the fusion of involutions in $G$ centreing around the Thompson
Transfer Lemma and the fact that $G = O^2(G)$ eventually establish the structure of $C_G(t)$
required to invoke Theorem~\ref{E8Sylow}.

In Section~\ref{sec8}, we start to investigate the situation when $O_{p'}(H)=1$ and $E=E(H)\neq 1$.
By Lemma~\ref{Eqs} we already know that $E$ is a quasisimple group. We further assume that $E$ is a
$\mathcal K$-group. In Lemma~\ref{Op=1} we show that $O_p(H)=1$. We remark that our proof of this
result already requires the $\mathcal K$-group hypothesis on $H$ as well as on $N_G(T)$ for
non-trivial $2$-subgroups $T$ of $G$. A very useful though easy consequence of Lemma~\ref{Ctnot-in}
in conjunction with Propositions~\ref{SE-p2} and \ref{Hstru1} is that if $C_H(t)/O_{p'}(C_G(t))$ is
not soluble, then $F^*(C_H(t)/O_{p'}(C_G(t))) \cong \Alt(p) \times \Alt(p)$ with $p\ge 5$ or to
$\Omega_8^+(2)$ and $p=5$ (see Lemma~\ref{QS1}). Now for most candidates for $E$, we can select an
involution $t \in E$ such that $C_E(t)$ is non-soluble and, since we can show that $p$ divides
$|C_E(t)|$ (Lemma~\ref{QS3}), we finish up only having to consider small Lie rank simple groups
defined over small fields in any significant detail. In these cases, when there is more than one
class of involutions in $E$, it is usually possible to choose a further involution $s$ such that
$m_p(C_H(s)) < 2$. Thus the meat in this section is contained in dealing with the possibility that
$E$ might be a rank $1$ Lie type group or that $E \cong \PSL_3(2^a)$ for some $a \ge 1$.  The
arguments dispatching these possibilities are presented in Lemmas~\ref{QS5} and \ref{QS4}. Finally
we present  proofs of Theorems~\ref{MainTheorem0}, \ref{MainTheorem1} and \ref{MainTheorem} as well
as Corollary~\ref{LiepCor} in Section~\ref{lastsec}.

Our group theoretical notation is mostly standard and can be found in \cite{GLS2}. Particularly we
use the following notation and conventions.  Suppose that $X$ is a finite group. We use $O(X)$ to
represent $O_{2'}(X)$, the largest normal subgroup of $X$ of odd order. As usual $Z^*(X)$ is the
subgroup of $X$ which contains $O(X)$ and satisfies $Z^*(X)/O(X)= Z(X/O(X))$. For $Y \le X$ and $x
\in X$, $x^Y$ denotes the set of $Y$-conjugates of $x$. If $a \in x^Y$, we write $a\sim_Y x$ to
indicate that $a$ and $x$ are conjugate by an element of $Y$.  If $p$ is a prime, $P \in \syl_p(G)$
and $x \in Z(P)^\#$, then $x$ is called a $p$-central element. Our notation for the simple groups
is also standard or self-explanatory. The dihedral group of order $n$ is denoted by $\Dih(n)$, The
semi-dihedral group of order $n$ is denoted by $\SDih(n)$ and the generalized quaternion group of
order $2^a$ is denoted by $\Q_{2^a}$ and is usually referred to as a quaternion group of order
$2^a$. The extraspecial $2$-groups $P$ of order $2^{1+2n}$ are denoted $2^{1+2n}_+$ if $P$ has an
elementary abelian subgroup of order $2^{n+1}$ and are otherwise denoted by $2^{1+2n}_-$. If $A$
and $B$ are groups, then $A*B$ represents the central product of $A$ and $B$, $A:B$ denotes a
semidirect product of $A$ and $B$ with undefined action, and $A\udot B$ denotes a non-split
extension of $A$ by $B$.

\bigskip

 \noindent {\bf Acknowledgement.} The first author thanks the
Institut f\"ur Mathematik, Universit\"at Halle--Wittenberg and the second author for their
hospitality during visits to Halle to carry out research on this project.
\section{Preliminaries}\label{prel}

In this section we gather together an eclectic collection of preliminary results which we will
invoke at various places throughout the paper. We begin with the principal recognition theorems
that we shall apply.

\begin{definition}\label{defclassicalinv}
Suppose that $G$ is a group, $t$ is an involution in $G$ and $C=C_G(t)$. Then $t$ is a
\emph{classical involution in $G$} provided there exists a subnormal subgroup $R$ of $C$ which
satisfies the following conditions.
\begin{enumerate}
\item $R$ has non-abelian Sylow $2$-subgroups and $t$ is
the unique involution in $R$;
\item $r^G \cap C \subset N_G(R)$ for all $2$-elements  $r\in R$;
and
\item $[R^g,R]\le O_{2'}(C)$ for all $g \in C\setminus N_G(R)$.
\end{enumerate}
\end{definition}

Notice that if $G$ is a group, $t$ is an involution in $G$ and  $R$ is a normal subgroup of
$C_G(t)$ containing $t$ which has quaternion Sylow $2$-subgroups, then $t$ is a classical
involution. The following result is of fundamental importance in our work.

\begin{theorem}[Aschbacher]\label{clasinv} Suppose that $G$ is a group and $F^*(G)$ is simple. If
$G$ has a classical involution then $G$ is a $\mathcal{K}$-group.
\end{theorem}

\begin{proof} This is the main theorem in \cite{Aschbacher}. \end{proof}

\begin{theorem}[Walter, Bombieri, Thompson]\label{E8Sylow} Suppose that $G$ is a perfect group, $S
\in \Syl_2(G)$ and $t\in G$ is an involution. If $S$ is abelian of order $8$ and $C_G(t)$ contains
a normal subgroup isomorphic to $\PSL_2(p^a)$, $p$ an odd prime, then $G\cong \J_1$ or
${}^2\G_2(3^{a})$ with $a\ge 3$, $a$ odd.
\end{theorem}

\begin{proof} From the main result in \cite{Walter} we
have that $C_{F^\ast(G)}(t) = \langle t \rangle \times \PSL_2(q)$. The assertion then follows from
\cite{Thompson} and \cite{Bombieri}.
\end{proof}

We now move on to properties of $\mathcal K$-groups with strongly $p$-embedded subgroups.

 {\definition Suppose that $p$ is a prime, $X$ is a group and $P \in \syl_p(X)$. Then the
\emph{$1$-generated $p$-core} of $X$ is defined as follows:
$$\Gamma_{P,1}(X)= \langle N_X(R) \mid 1\not= R \le P\rangle.$$}

We say that a proper subgroup $M$ of $X$ is \emph{strongly $p$-embedded} if and only if $p$ divides
$|M|$ and $\Gamma_{P,1}(X) \le M$ for some Sylow $p$-subgroup $P$ of $M$. That this definition is
equivalent to the one given in the introduction is an easy application of Sylow's Theorem (see
\cite[17.10]{GLS2} for example).

\begin{proposition} \label{SE-p2} Suppose that $p$ is a prime, $X$ is a $\mathcal K$-group and
$K=F^*(X)$ is simple. Let  $P \in \syl_p(X)$ and $Q= P\cap K$. If
$m_p(P) \ge 2$ and $\Gamma_{P,1}(X) <X$, then $\Gamma_{Q,1}(K) <K$
and  $p$ and $K$ are as follows.
\begin{enumerate}
\item $p$ is arbitrary, $a \ge 1$ and $K \cong \PSL_2(p^{a+1})$,
or $\PSU_3(p^a)$,  ${}^2\B_2(2^{2a+1})$ $(p=2)$ or ${}^2{\rm G}_2(3^{2a+1})$ $(p=3)$ and $X/K$ is a
$p'$-group.
 \item $p > 3$, $K \cong \Alt(2p)$ and
$|X/K| \le 2$.
\item $p=3$, $K \cong \PSL_2(8)$ and $X \cong \PSL_2(8):3$.
\item $p=3$, $K \cong \PSL_3(4)$ and $X/K$ is a $2$-group.
\item $p=3$ and $X=K\cong \M_{11}$.
\item $p=5$ and $K \cong {}^2\B_2(32)$ and $X \cong {}^2\B_2(32):5$.
\item $p=5$, $K \cong
{}^2\F_4(2)^\prime$ and $|X/K|\le 2$.
 \item $p=5$, $K \cong \McL$ and $|X/K|\le 2$.
\item $p=5$, $K  \cong \Fi_{22}$ and $|X/K|\le 2$. \item
$p=11$ and  $X=K \cong \J_4$.
\end{enumerate}\end{proposition}

\begin{proof} See \cite [7.6.1]{GLS3}  or \cite[24.1]{GorensteinLyons}.\end{proof}

\begin{notation}We let $\mathcal E$ be the set of pairs $(K,p)$ in the conclusion  of Proposition~\ref{SE-p2}
with $p$ odd.
\end{notation}

\begin{proposition}\label{Hstru1} Suppose that $X$,  $K$, $P$ and $Q$ are as in
Proposition~\ref{SE-p2}. Let $H = \Gamma_{Q,1}(K)$. Then the
following hold.
\begin{enumerate}
\item If \ref{SE-p2} (i) holds, then $H= N_K(Q)$ is a Borel subgroup of
$K$. In particular, $H= JQ$ where $J$ is cyclic and operates
irreducibly on each of $Q/\Phi(Q)$ and $Z(Q)$.
\item If \ref{SE-p2}(ii) holds, then $H$ is isomorphic to the subgroup of
$\Sym(p)\wr \Sym(2)$ consisting of even permutations. In particular,
if $p > 3$, then $O_p(H)=1$.
\item If \ref{SE-p2}(iii) holds, then $Q$ is cyclic of order $9$, $H = N_K(Q)
$ is dihedral of order $18$ and $N_X(P) \cong 3^{1+2}_-.2$. In
particular, $N_K(Q)$ does not act irreducibly on $P/\Phi(P)$.
\item If \ref{SE-p2}(iv) holds, then $Q$ is elementary abelian of order $9$, $H= N_K(Q) $, $C_K(Q) = Z(Q)$, $H/Q
\cong \Q_8$ and the complement to $Q$ in $H$ acts irreducibly on $Q$.
\item If \ref{SE-p2}(v) holds, then $Q$ is elementary abelian of order $9$,
$H=N_K(Q)$, $C_H(Q)= Q$, $H/Q \cong \SDih(16)$ and a
complement to $H$ in $K$ acts irreducibly on $Q$.

\item If \ref{SE-p2}(vi) holds, then $Q$ is cyclic of order $25$, $H = N_K(Q)
$ is a Frobenius group of order $100$, and $N_X(P) \cong 5^{1+2}_-.4$. In particular, a complement
to $Q$ in $N_K(Q)$ does not act irreducibly on $P/\Phi(P)$.

\item If \ref{SE-p2}(vii) holds, $Q$ is elementary abelian of order $25$, $H =
N_K(Q)$, $Q= C_H(Q)$, $H/Q$ is isomorphic to a central product of a
cyclic group of order $4$ with $\SL_2(3)$ and a complement to $Q$ in
$H$ acts irreducibly on $Q$.

\item If \ref{SE-p2}(viii) holds, then $Q$ is extraspecial of order $5^3$,
$H=N_K(Q) $, $C_K(Q)= Z(Q)$, $H/Q$ isomorphic to a non-abelian
extension of a cyclic group of order $3$ by a cyclic group of order
$8$ and a complement to $Q$ in $H$ acts irreducibly on $Q/\Phi(Q)$.

\item If \ref{SE-p2}(ix) holds, then $Q$ has order $25$ and $H \cong
\Aut(\Omega_8^+(2))$. In particular $O_5(H)=1$.
\item If \ref{SE-p2}(x) holds, then $Q$ is extraspecial of order $11^3$, $H=
N_K(Q)$, $C_K(Q)= Z(Q)$, $H/Q$ is isomorphic to the direct product of a cyclic group of order $5$
and $\GL_2(3)$. In particular, a complement to $Q$ in $H$ acts irreducibly on $Q/\Phi(Q)$.
\end{enumerate}

\end{proposition}

\begin{proof} This is mostly \cite[Theorem 7.6.2]{GLS3}, the statement regarding the action of complements to $Q$ in $H$ are easily deduced and are
well-known.
\end{proof}

\begin{corollary}\label{OrredH} Assume that $X$, $K$, $P$, $Q$  and $H$  are as in
Proposition~\ref{Hstru1}. If $O_p(H) \neq 1$, then  $N_H(P)$ acts irreducibly on $P/\Phi(P)$,
unless either
\begin{enumerate}
\item $p=3$ and $X \cong\PSL_2(8):3$; or

\item $p=5$ and $X\cong {}^2\B_2(32):5$.
\end{enumerate}
\end{corollary}

\begin{proof} The assertion follows directly from Proposition~\ref{Hstru1}.
\end{proof}

\begin{corollary}\label{OrredH2} Assume that $X$, $K$, $P$, $Q$  and $H$  are as in
Proposition~\ref{Hstru1} with $p$ odd. Then $H$ contains an involution inverting $P/\Phi(P)$ unless
cases (i) or (ii)  of Corollary~\ref{OrredH} holds or $K \cong \PSL_2(p^a)$ with $p^a\equiv 3\pmod
4$.
\end{corollary}

\begin{proof} As long as $N_H(P)/C_H(P)P$ has a central involution the result follows with Corollary~\ref{OrredH}.
By Proposition~\ref{Hstru1}, if $N_H(P)/C_H(P)$ does not have a central involution, we must be in
case (i) of Proposition~\ref{SE-p2}. Furthermore the Borel subgroup must be of odd order. Therefore
$K \cong \PSL_2(p^a)$ with $p^a\equiv 3\pmod 4$. \end{proof}

\begin{lemma}\label{L34} Let $F^\ast(X) = K \cong \PSL_3(4)$ and $t \in X$  be an involution which induces an outer
automorphism on $L$. Then $C_{K}(t) \cong 3^2.\mathrm \Q_8\cong \PSU_3(2)$, $\PSL_3(2)$ or
$\Alt(5)$.\end{lemma}

\begin{proof} By \cite[2.5.12]{GLS3} $K$ just posseses involutory outer automorphisms which are field-, graph,- or graph-field
automorphisms. A field automorphism centralizes $\PSL_3(2)$, a graph automorphism centralizes
$\PSL_2(4) \cong \Alt(5)$ and a graph-field automorphism centralizes $\PSU_3(2)$. \end{proof}

\begin{corollary}\label{OrredH3} Assume that $X$, $K$, $P$, $Q$   are as in
Proposition~\ref{SE-p2}. Set $H = \Gamma_{P,1}(X)$. If $O_{p'}(H) \neq 1 $, then $K \cong
\PSL_3(4)$ or $\Fi_{22}$, $|X/K|\ge 2$ and $O_{p'}(H)=Z(H)$ has order $2$.
\end{corollary}

\begin{proof} This is \cite[24.2]{GorensteinLyons}.

\end{proof}

\begin{lemma}\label{SchurM} Suppose that $p$ is an odd prime and  $(K,p) \in \mathcal
E$. If two divides the order of the Schur multiplier of $K$, then $K$ is isomorphic to
$\PSL_2(p^a)$ for some $a$, $\Alt(2p)$, $\PSL_3(4)$, or $\Fi_{22}$.
\end{lemma}

\begin{proof} For this we just consult \cite[Theorem 6.1.4]{GLS3}.
\end{proof}

\begin{corollary}\label{CExt}  Assume that $p$ is an odd prime,  $(K,p) \in \mathcal E$ and $H= \Gamma_{Q,1}(K)$ where $Q\in \syl_p(K)$.
If $\widehat K$ is a non-trivial
perfect central extension of $K$ with $Z(K)$ a $2$-group and if
$\widehat H$ has quaternion or cyclic  Sylow $2$-subgroups, then
$\widehat K \cong \SL_2(p^a)$ for some $a$.
\end{corollary}

\begin{proof} Suppose that $K \not \cong \PSL_2(p^a)$. Then, by Lemma~\ref{SchurM}, $K \cong \Alt(2p)$, $\PSL_3(4)$ or $\Fi_{22}$.
In these cases, Proposition~\ref{Hstru1} shows that $H$ does not have cyclic or dihedral Sylow
$2$-subgroups and so these groups do not arise. It is well-known that the Sylow $2$-subgroups of
$\SL_2(p^a)$ are quaternion.
 \end{proof}

\begin{lemma}\label{invs}
 Suppose that $p$ is an odd prime and $X$ is  group with $F^*(X)=K$ where  $(K,p)\in \mathcal E$.
Let $P \in \syl_p(X)$ and $Q = P \cap K$. Set $H=\Gamma_{Q,1}(K)$.  If $|K:H|$ is odd, then $K
\cong \M_{11}$.
\end{lemma}

\begin{proof} This follows by inspection of the groups in $\mathcal
E$ and the subgroups corresponding to $H$.
\end{proof}

\begin{lemma}\label{invs2}
 Suppose that $p$ is an odd prime and $X$ is  group with $F^*(X)=K$ where  $(K,p)\in \mathcal E$.
Let $P \in \syl_p(X)$ and $Q = P \cap K$. Set $H=\Gamma_{Q,1}(K)$.  If $|H|$ is odd, then $K \cong
\PSL_2(p^f)$ for some odd $f$ and $p \equiv 3 \pmod 4$.
\end{lemma}

\begin{proof} This follows by inspection of the groups in $\mathcal
E$ and the subgroups corresponding to $H$ described in Proposition~\ref{Hstru1}.
\end{proof}

\begin{lemma}\label{prankauto} Suppose that $K$ is a simple $\mathcal K$-group and
that $p$ is an odd prime. Then $m_p(\Out(K)) \le 2$ and, if $P \le
\Out(K)$ has $p$-rank $2$, then $P \not\le O^p(\Out(K))$.
\end{lemma}

\begin{proof} Since the alternating groups and sporadic simple groups
 have outer automorphism groups which are $2$-groups (perhaps trivial), we may suppose
 that $E$ is a Lie type group. Then the result can be deduced from
 \cite[Theorem 2.5.12]{GLS3}.
 %Check number (and result)!
\end{proof}

\begin{lemma}\label{wcl} Assume that $G$ is a group with $Z^*(G)= O_{p'}(G)$.
If $t \in G$ is an involution, then $t^G\cap C_G(t)\neq \{t\}$.
\end{lemma}

\begin{proof} Let  $R \in\syl_2(C_G(t))$ and assume that $t^G\cap C_G(t)\neq \{t\}$. Then $t^G \cap R=\{ t\}$ and so $t \in Z(N_G(R))$. Therefore, $R
\in \Syl_2(G)$ and  Glauberman's $Z^*$-Theorem \cite{Glau} implies that $t \in Z^*(G)= O_{p'}(G)$
which is impossible. Hence the lemma holds.
\end{proof}

The next result is the famous Thompson Transfer Lemma.

\begin{lemma}\label{Thompsontransfer} Let $G$ be a group,
$S\in {\rm{Syl}}_2(G)$, $T \unlhd S$ with $S = TA$, $A \cap T = 1$, $A$ cyclic and non-trivial. If
$G$ has no subgroup of index two and $u$ is the involution in $A$, then there is some $g \in G$
with $u^g \in T$ and $C_S(u^g)\in {\rm{Syl}}_2(C_G(u^g))$. In particular $|C_S(u)|\leq |C_S(u^g)|$.
\end{lemma}

\begin{proof} This is \cite[(15.16)]{GLS2}.
\end{proof}

\begin{lemma}\label{ThompsonTransferC} Suppose that $G$ is a group and $S\in \Syl_2(G)$. Then
at least one of the following hold.
\begin{enumerate}
\item $G\neq O^2(G)$;
\item $\Omega_1(Z(S)) \le \Phi(S)$; or
\item $N_G(S)$ acts non-trivially on $\Omega_1(Z(S))$.
\end{enumerate}
\end{lemma}

\begin{proof} Suppose that (i) and (ii) do not hold. Then there is a
maximal subgroup $M$ of $S$ and an involution $t \in \Omega_1(Z(S))$ such that $t \not \in M$. By
the Thompson Transfer Lemma \ref{Thompsontransfer}, there is $x \in G$ such $t^x \in M$ and
$C_S(t^x) \in \Syl_2(C_G(t^x))$. Since $t$ is $2$-central, it follows that $t^x \in
\Omega_1(Z(S))$. Burnside's Lemma \cite[16.2]{GLS2} now implies that $t$ and $t^x$ are conjugate by
an element of $N_G(S)$. Thus (iii) holds.
\end{proof}

\begin{lemma}\label{cyc2}Suppose that $G$ is a group with $G= O^2(G)$, $t\in G$ is an involution and  $S \in \Syl_2(C_G(t))$.
Assume that  $S = \langle y\rangle \times S_0$, $t\in \langle y\rangle$ and $Z(S_0)$ is elementary
abelian. Then $\langle t \rangle =\langle y\rangle$.
\end{lemma}

\begin{proof} Assume that $y$ has order greater than $2$. Then $Z(S) = \langle y\rangle \times
Z(S_0)$ and so $t$ is the unique involution in $\Phi(Z(S))$. Therefore $\langle t \rangle $ is a
characteristic subgroup of $S$. Hence $t\in Z(N_G(S))$ and so $S \in \Syl_2(G)$. As $S = \langle
y\rangle \times S_0$ and $G= O^2(G)$, we may now apply the Thompson Transfer
Lemma~\ref{Thompsontransfer} to see that there exists  $g \in G$ such that $t^g \in S_0$ and
$C_S(t^g)\in \Syl_2(C_G(t^g))$. It follows that $t^g \in Z(S)$. Hence $t$ and $t^g$ are conjugate
in $N_G(S)$ by Burnside's Lemma. But we have already noted that $\langle t\rangle$ is a
characteristic subgroup of $S$ and so we have $\langle t\rangle =\langle t^g\rangle$ which is a
contradiction. Thus $\langle t\rangle =\langle y\rangle$ as claimed.
\end{proof}

\begin{lemma}\label{oneclass} Suppose that $G$ is a group, $r$ is an odd prime and $G/O_2(G)\cong \Dih(2r)$.  Let $R \in \Syl_r(G)$
and assume that $C_{O_2(G)}(R) =1$. Then $N_G(R) \cong \Dih(2r)$ and, for $a$ an involution in $
N_G(R)$, all the involutions of $G\setminus O_2(G)$ are conjugate to $a$ and $|C_{O_2(G)}(a)|^2=
|O_2(G)|$.
\end{lemma}

\begin{proof} Set $Q= O_2(G)$. Since $N_Q(R)= C_Q(R)=1$, the Frattini argument shows that
$N_G(R) \cong \Dih(2r)$. Let $a, b \in N_G(R)$ be involutions with $a\neq b$.  Assume that $ c \in
bQ$ is an involution. Then $\langle a, c\rangle$ is a dihedral group of order divisible by $r$.
Thus $\langle a, c\rangle$ contains a conjugate, $R^x$ of $R$. Since $\langle a,c \rangle \cap Q
\le C_Q(R^x) =1$, we have that $\langle a,c\rangle \cong \Dih(2r)$. Thus $a$ and $c$ are
$G$-conjugate. It follows that every involution in $G\setminus Q$ is conjugate to $a$.  Assume that
$a$ normalizes $R^x$ for some $x\in Q$. Then $xax^{-1} \in N_G(R)$ and $xax^{-1}a \in Q \cap
N_G(R)=1$. Thus $x \in C_Q(a)$. It follows that $a$ normalizes exactly $|C_Q(a)|$ conjugates of
$R$. Let $x \in bQ$ be an involution. Then $\langle a,x\rangle $ contains a conjugate of $R$ and if
$y \in bQ$ is an involution with $\langle a,x\rangle = \langle a,y\rangle$, then $x=y$. Thus $a$
normalizes at least $|Q:C_Q(b)|$ conjugates of $R$. It follows that $|Q| \le |C_Q(b)||C_Q(a)|=
|C_Q(b)|^2$. On the other hand, as  $C_Q(b)\cap C_Q(a) \le C_Q(R)=1$, we have $|C_Q(b)|^2 \le
|Q|$ and this completes the proof of the lemma. \end{proof}

We recall that if $A = \Z_{2^k}\times \Z_{2^k}$ and $G= \Aut(A)$. Then $G \cong \GL_2(\Z_{2^k})$
and that
$$O_2(G)=\left\{\left(\begin{array} {cc} a&c\\d&b\end{array}\right)
\mid
 a, b \in \Z_{2^k}^*, c,d \in 2\Z_{2^k}\right\}
 ,$$
where $ \Z_{2^k}^\ast \mbox{ denotes the groups of units of } \Z_{2^k}.$ In particular, we
 have that $|O_2(G)| = 2^{4(k-1)}$ and $G/O_2(G)\cong \SL_2(2)$.

\begin{lemma}\label{autohomocyclic} Suppose that $k\ge 2$,  $A = \Z_{2^k}\times
\Z_{2^k}$ and $G= \Aut(A)$. If $H \le G\cong \GL_2(\Z_{2^k})$ and $H \cong \Alt(4)$, then
$$O_2(H) = \left\langle \left(\begin{array} {cc} 1+2^{k-1}&0\\2^{k-1}&1+2^{k-1}\end{array}\right)
,\left(\begin{array} {cc} 1&2^{k-1}\\2^{k-1}&1\end{array}\right)\right\rangle.$$ In particular,
there is exactly one conjugacy class of subgroups of $G$ isomorphic to $\Alt(4)$.
\end{lemma}
\begin{proof} We first note that $t=\left(\begin{array} {cc}
-1&-1\\1&0\end{array}\right)$ is an element of order $3$ in $G$ and, since $|G|_3=3$, we may
suppose that $t \in H$.

We proceed by induction on $k$. Suppose that $k=2$. Then $O_2(G)$ has order $16$ and is abelian.
Since $O_2(H) \le O_2(G)$  and $Z(G)\neq 1$, we then have $O_2(H)= [O_2(G),t]= \left\langle
\left(\begin{array} {cc} 3&0\\2&3\end{array}\right) ,\left(\begin{array} {cc}
1&2\\2&1\end{array}\right)\right\rangle.$

Now suppose that $k> 2$. Let $B= 2A \cong \Z_{2^{k-1}}\times\Z_{2^{k-1}}$. Then $C_G(B)= \left
\{\left(\begin{array} {cc} 1+2^{k-1}a&2^{k-1}c\\2^{k-1}d&1+2^{k-1}b\end{array}\right)\mid a,b,c, d
\in \{0,1\}\right\}$ which has order $2^4$. Furthermore, $G/C_G(B) \cong \Aut(B)$. If $H \cap
C_G(B) = O_2(H)$, then $O_2(H) = [C_G(B),t] = \left\langle \left(\begin{array} {cc}
1+2^{k-1}&0\\2^{k-1}&1+2^{k-1}\end{array}\right) ,\left(\begin{array} {cc}
1&2^{k-1}\\2^{k-1}&1\end{array}\right)\right\rangle.$ and we are done. Hence $HC_G(B)/C_G(B) \cong
\Alt(4)$. By induction, we have $O_2(H) \le R= \left\langle \left(\begin{array} {cc}
1+2^{k-2}&0\\2^{k-2}&1+2^{k-2}\end{array}\right) ,\left(\begin{array} {cc}
1&2^{k-2}\\2^{k-2}&1\end{array}\right)\right\rangle C_G(B)$. However, an easy calculation then
shows that every element of $R\setminus B$ has order $4$ and this contradicts $O_2(H)$ having
exponent $2$. This proves the lemma.
\end{proof}

%The outcome of the next lemma is similar to \cite[Lemma 8.2]{GLS4}.

\begin{lemma}\label{2group} Let $G$ be a group, $T \le G$ be a $2$-group and $V \le Z(T)$ be a fours group.
Assume that $t, \rho \in N_G(V) \cap N_G(T)$ are elements of order $2$ and $3$ respectively with
$[t,\langle \rho,V\rangle]=1$. If $[V,\rho]=C_T(t)=V$, then $T$ is isomorphic to one of the
following groups.
\begin{enumerate}
\item[(i)] An elementary abelian group of order 16.
\item[(ii)] A homocyclic group of rank $2$, $\Z_{2^n} \times \Z_{2^n}$, $n \geq 1$.
\item[(iii)] A Sylow $2$-subgroup of $\PSL_3(4)$.
\item[(iv)] A Sylow $2$-subgroup of $\PSU_3(4)$.
\end{enumerate}
\end{lemma}
\begin{proof}   Because $C_T(t) = V$ and $[t,\rho]=1$, we have that $C_T(\rho) = 1$. If $T=V$, then (ii) holds with $n=1$. So assume that $T > V$ and
let $W$ be the preimage of $C_{T/V}(t)$. Then, since $C_T(\rho)=1$ and $C_T(t) = V$, we have that $|W : V| = 4$. In
particular, $W/V$ has rank $2$ which means that $|W'|\le 2$. Thus,  as $\rho$ acts on $W^\prime$,
we have that $W$ is abelian and so  $W$ is either elementary abelian or homocyclic $\Z_4\times
\Z_4$

Assume first that $W$ is elementary abelian. Then all involutions in $Wt$ are conjugate in $\langle
W,t\rangle$. As $C_T(t) = V$, this means that $T = W$, which is (i).

Assume next that $W$ is homocyclic. Then, as $W/V = C_{T/V}(t)$ is elementary abelian of order $4$,
$T/V$ satisfies the hypothesis of the lemma. Hence we may assume by induction $T/V$ is isomorphic
to one of the groups in (i) - (iv).

Assume  that $T/V$ is elementary abelian of order 16. If $V = \Omega_1(T)$, then, as $T/V$ can be
considered as a direct some of two irreducible $2$-dimensional  modules for $\langle \rho\rangle$,
we have that $T = WW_1$, where $W_1$ is also homocyclic of order 16 and is normalized by $\rho$.
Furthermore, $C_{W_1}(W) = C_W(T)=V$. Using Lemma~\ref{autohomocyclic}  shows that the action of
$W_1$ on $W$ is uniquely determined. As $W = [t,W_1]$, we see that $T$ is uniquely determined and
so we have that $T$ is isomorphic to  a Sylow 2-subgroup of $\PSU_3(4)$. If $\Omega_1(T) \neq V$,
then we have $T = WW_1$, where $W_1$ is elementary abelian of order 4 and is normalized by $\rho$.
If $C_{W_1}(W) = 1$, then the action of $W_1$ on $W$ is uniquely determined as  above. We therefore
have that $VW_1$ and $(VW_1)^t$ are elementary abelian of order 16 and the action of $W_1$ on
$(VW_1)^t$ is uniquely determined. It follows that $T$ is isomorphic to a Sylow 2-subgroup of
$\PSL_3(4)$. Hence we have $C_{W_1}(W)= W_1$ and consequently $T$ is abelian. Then $\Omega_1(T)$ is
elementary abelian. But this forces $W=\Omega_1(T)$ which is a contradiction as $W \cong \Z_4
\times \Z_4$.

Assume next that $T/V$ is isomorphic to a Sylow 2-subgroup of $\PSU_3(4)$ or $\PSL_3(4)$. Then
there is a homocyclic group $W_1 \cong \Z_8 \times \Z_8$ on which $\langle T/W_1, \rho \rangle$
acts as $\Alt(4)$. Using Lemma~\ref{autohomocyclic} we get  $[W_1,T] =V$ and this contradicts
$Z(T/V)= W/V$. Hence we finally have $T/V$ is homocyclic. As $T'= \langle [x,y]\rangle$ where $T =
\langle x, y\rangle$,  the fact that $C_{T^\prime}(\rho) = 1$ implies that $T$ is abelian. As
$\Omega_1(T) = V$ and $C_T(\rho)=1$, we have that $T$ is homocyclic.
\end{proof}

\begin{corollary}\label{2groupCor} Assume the hypothesis of Lemma~\ref{2group}. If the coset $T t$ contains
more than one $T$-conjugacy class of involutions, then $T$ is homocyclic and $Tt$ contains exactly
four $T$-conjugacy classes of involutions.
\end{corollary}

\begin{proof} We consider each of the possibilities for $T$ given in Lemma~\ref{2group}. If $T$ is
elementary abelian, then  $C_T(t)= V$ has index $4$ in $T$ and so $Tt$ contains exactly $4$
involutions and they are all $T$-conjugate. Assume now that $T$ is non-abelian. Let $W = C_{T/V}(t)$.
We note that $|t^T|= |T/V|=
|T|/4$ and, if $x \in T$ and $(xt)^2=1$, then $t$ inverts $x$. As $T$ is non-abelian, then $T$ has
exponent $4$ and so $x^2 \in V$. We infer that, if $x \in T$ is inverted by $t$, then $x \in W$. It
follows that $Tt$ contains at most $|W|$ involutions. Since $|W| = |T/V|$, we have that $Tt$
contains exactly one $T$-conjugacy class of involutions in this case. Finally, if $T$ is
homocyclic, we have that $t$ inverts every element of $T$. Thus $Tt$ consists of involutions and it
follows that $Tt$ contains exactly four $T$-conjugacy classes of involutions.
\end{proof}

\begin{lemma}\label{smallgroup} Let $t \in G$ be an involution and $T$ be a Sylow $2$-subgroup of $C_G(t)$ such that $T \cong \langle t \rangle \times D$,
where $D$ is dihedral of order $4$ or $8$. Assume that there is a fours group $V \leq D$, such that
$C_G(t)$ contains a 3-element $\rho$ which acts non-trivially on $V$. Then there is a Sylow
$2$-subgroup $R$ of $C_G(V)$ which is normalized by $\langle T, \rho \rangle$ such that $R =
U\langle t \rangle$ where $U$ is isomorphic to one of the groups listed in the conclusion of
Lemma~\ref{2group}.
\end{lemma}

\begin{proof} By considering a minimal counter example to the lemma, we may assume that $G= N_G(V)= C_G(V)T \langle
\rho\rangle$. Set $E = \langle V, t\rangle$. As $G$ is a counter example to the lemma, we have $V
\not \in \syl_2(C_G(V))$. Set $E = \langle V, t\rangle$. Then $E$ is elementary abelian of order
$2^3$. As $T \in \Syl_2(C_G(t))$, we have $E \in \syl_2(C_G(E))$. Thus $C_G(E) = E \times
O_{2'}(C_G(E))$. Assume that $S \in \Syl_2(C_G(V))$ is normalized by $T$. Then $W=N_S(E) > E$ and
since $W=N_S(E)$ centralizes $V$, we deduce that $\langle \rho \rangle WC_G(E)/C_G(E) \cong
\Alt(4)$. In particular, we have $|W|= 2^5$ and $W \in \syl_2(N_G(E))$. If $W^\rho \neq W$, then
there exists $c \in O_{2'}(C_G(E))$ such that $\rho c$ normalizes $W$. Let $\rho^*$ be the $3$-part
of $\rho c$. Then $\rho^*$ acts non-trivially on $V$ and centralizes $t$. Thus we may as well
suppose that $\rho$ normalizes $W$. Set $\ov G = G/V$. Then $\ov {\langle T \rangle W} \in
\syl_2(C_{\ov {N_G(V)}}(\ov t))$ and $\ov \rho $ is a $3$-element which acts non-trivially on $\ov
[W,\rho]$ and centralizes $\ov t$. Therefore, by induction, there is a subgroup $U \ge V$ such that $\ov
U\langle \ov t\rangle \in \syl_2(C_{\ov{G}}(\ov {[W\rho]}))$ where $\ov U$ is one of the groups listed in
the conclusion of Lemma~\ref{2group} and  $\ov {\langle \rho, T\rangle}$ normalizes $\ov U$. By
considering $U\langle T, \rho\rangle$, we have that $U$ is also listed in the conclusion of
Lemma~\ref{2group}. We may assume that $U\langle t \rangle \le S$ and, as $G$ is a counter example
to the lemma,  $S \neq U\langle t \rangle$. Since $t$ is not centralized by an abelian group of
order $16$ and $U \neq V$, we see that $t^G \cap U = \emptyset$.

If $U$ is either elementary abelian or is isomorphic to a Sylow $2$-subgroup of $\PSL_3(4)$ or
$\PSU_3(4)$, then, by Corollary~\ref{2groupCor},  every involution in $U\langle t \rangle \setminus
U$ is conjugate to $t$ by an element of $U$. Setting $X = N_S(U\langle t\rangle)$, we then have
$U\langle t \rangle <X = C_X(t)U\langle t\rangle= U\langle t\rangle$, a contradiction. Thus $U$ is
homocyclic of rank $2$. Since $|U| \ge 16$, we have that $U$ is a characteristic subgroup of
$U\langle t\rangle$. It follows that $X $ normalizes $\Omega_2(U) = [W\rho]$. Since $X$ centralizes $V$, we have
$[[W\rho],X]
\le V$ and thus $\ov X \le C_{\ov G}(\ov {[W,\rho]})$ and this is our final contradiction as $\ov {U\langle
t\rangle} \in \syl_2(C_{\ov G}(\ov {[W,\rho]}))$. This concludes the verification of Lemma~\ref{smallgroup}.
\end{proof}

\begin{lemma}\label{AutosE} Assume that $p$ is an odd prime, $W$ is a $2$-group and that $E \le
\Aut(W)$ is a non-cyclic abelian $p$-subgroup. If $C_W(e)$ contains at most one involution for each
$e \in E^\#$, then $|E|=9$, $W \cong 2^{1+4}_+$, $2^{1+6}_-$ or $2^{1+8}_+$ and $C_W(E)= Z(W)$.
\end{lemma}

\begin{proof} Let $W$ be a minimal counter example to the claim.
Since $E$ is abelian and non-cyclic, $W = \langle C_W(e)\mid e \in E^\#\rangle$. Choose $e$ such
that $C_W(e) \not =C_W(E)$. Then, as $C_W(e)$ contains exactly one involution, $C_W(e)$ is cyclic
or generalized quaternion. Because $E$ normalizes and does not centralize $C_W(e)$,  it follows
that $C_W(e)$ is a quaternion group of order $8$ and $p=3$. In particular, we have that $C_W(E) =
Z(C_W(e))$ has order $2$ and if $C_W(f) > C_W(E)$ for some $f \in E^\#$, then $C_W(f)$ is a
quaternion group of order $8$. We also have that $W$ is non-abelian and $Z(W)= C_W(E)$. We now
choose $e\in E^\#$ such that $C_W(e)$ is contained in the second centre of $W$. Then $C_W(e)$ is
normal in $W$. Since the automorphism group of $C_W(e)$ is isomorphic to $\Sym(4)$, we get that
$W=C_W(e)C_W(C_W(e))$. Set $W_0 = C_W(C_W(e))$. Then $W_0 \not \le C_W(e)$, for otherwise $e$ would
be the trivial automorphism of $W$. Thus $W_0 $ is a non-trivial $2$-group. Suppose that $W_0 \le
C_W(f)$ for some $f \in E^\#$. Then, as $C_W(f) $ is quaternion of order $8$ and $W_0$ is
normalized by $E$, we infer that $C_W(f)= W_0$ and $W= C_W(e)C_W(f) \cong 2^{1+4}_+$. Furthermore,
we have that $E$ is elementary abelian of order $9$. Thus we may suppose that every element of $E$
induces a non-trivial automorphism of $W_0$. As $W_0 < W$, we have that $|E| =9$ and that $W_0
\cong 2^{1+4}_+$, $2^{1+6}_-$ or $2^{1+8}_+$ by induction. The first two cases immediately deliver
$W \cong 2^{1+6}_-$ or $W \cong 2^{1+8}_+$ and the theorem then holds. So suppose that $W_0 \cong
2^{1+8}_+$. Then, for each $f \in E^\#$, $C_{W_0}(f)$ is a quaternion group of order $8$. But $e
\in E^\#$ and $C_W(e) \cap W_0 = Z(W_0)$ and so we have a contradiction. This completes the proof
of the lemma.
\end{proof}

\begin{lemma}\label{AutSL3Sylow} Suppose that $X \cong \Aut(\PSL_3(2^a))$ with $a \ge 2$, $S \in \Syl_2(X)$ and $T
= S \cap F^*(X)$. If $U \le S$ and $U \cong T$, then $U=T$.
\end{lemma}

\begin{proof} Set $K= F^*(X)$. Suppose that  $U \le S$ with $U \cong T$ and $U \neq T$. Let $F_1$ and $F_2$ be the two elementary abelian subgroups
of $T$ of order $2^{2a}$ and note that every involution of $T$ is contained in $F_1 \cup F_2$ and
that $F_1 \cap F_2 = Z(T)$.  Since the group of diagonal outer automorphisms of $K$ has order $3$,
$\Out(K)$ has abelian Sylow $2$-subgroups. Therefore $Z(U)=\Phi(U)=U'\le T$ and $UK/K$ is
elementary abelian of order at most $4$. Let $E_1$ and $E_2$ be the elementary abelian subgroups of
$U$ of order $2^{2a}$. As $U = E_1E_2$, we may suppose that $E_1\not \le T$. As $|UK/K|\le 4$,
$|E_1 \cap T|\ge 2^{2a-2}$. Assume that $(E_1 \cap T) Z(T)\not = Z(T)$. Then we may suppose that $(E_1
\cap T)Z(T) \le F_1$. But then $F_1^u \ge (E_1\cap T)Z(T)
> Z(T)$ for all $u\in U$. Hence $U$ normalizes $F_1$ and therefore $UK/K$ consists of field outer automorphisms.
Hence $|E_1 \cap T | \ge 2^{2a-1}$. Now let $e \in E_1 \setminus T$. Then $e$ centralizes either
$(E_1 \cap T)Z(T)/Z(T) = F_1/Z(T)$ or $Z(T) \le E_1$ and $e$ centralizes $Z(T)$. In either case $e$ cannot be a field
automorphism of $K$ and we have a contradiction. Thus $E_1\cap T = Z(T)$, $|UT/T|= 4$ and $a=2$.
But then $Z(U) = E_1 \cap T=Z(T)$ and we see that $U$ centralizes $Z(T)$. This means that $UK/K$
does not contain non-trivial field outer automorphisms. It follows that $|UK/K|\le 2$ and this is a
contradiction.
\end{proof}

\section{Basic properties of groups with a  strongly $p$-embedded subgroups}\label{basic}

In this section we reveal the basic  structural properties of groups with strongly $p$-embedded
subgroups. The first lemma is one which we already alluded to in Section~\ref{prel} and states the
equivalence between the two definitions of strongly $p$-embedded subgroups which we have given.

\begin{lemma}\label{def} Assume that $G$ is a group, $p$ is a prime, $H\le G$ and $S \in \syl_p(H)$.
Then $\Gamma_{S,1}(G) \le H$  if and only if $p$ divides $|H|$ and   $|H\cap H^g|$ is not divisible by $p$ for
all $g \in G\setminus H$.
\end{lemma}

\begin{proof} See \cite[Lemma~17.11]{GLS2}.
\end{proof}

The fundamental lemmas which gets us started in the proof of Theorem~\ref{MainTheorem} are as
follows.

\begin{lemma}\label{main1} Suppose that $G$ is a group, $p$ is a prime, $H$ is a strongly
$p$-embedded subgroup of $G$ and $K \leq G$ such that $H \cap K$ has order divisible by $p$. Then the following statements hold.

\begin{enumerate}
\item $\syl_p(H) \subseteq \syl_p(G)$;
\item if $K \not \le H$ and , then $H \cap
K$ is a strongly $p$-embedded subgroup of $K$;
\item $\syl_p(H \cap K) \subseteq \syl_p(K)$;
\item $\bigcap H^G \le O_{p'}(G)$; and
\item if $m_p(H) \ge 2$, then $O_{p'}(G) = \bigcap H^G$.
\end{enumerate}
\end{lemma}

\begin{proof} Let $S \in \Syl_p(H)$. Then, as $H$ is strongly $p$-embedded in $G$, $N_G(S) \le H$.
Thus $S \in \syl_p(G)$ and (i) holds.

Suppose that $S$ is chosen so that $S \cap K \in \syl_p(H\cap K)$. Let $T$ be a non-trivial
subgroup of $S\cap K$. Then, as $H$ is strongly $p$-embedded, $N_K(T) \le N_G(T) \le H$, thus
$N_K(T) \le H\cap K$ and so we have that $H\cap K$ is strongly $p$-embedded in $K$. Thus (ii)
holds.

Part (iii) follows from (i) and (ii).

Since, by Lemma~\ref{def}, $p$ does not divide $|H \cap H^g|$, we have that $\bigcap H^G \le
O_{p'}(G)$. So (iv) holds.

Finally, assume that $m_p(H) \ge 2$ and suppose that $A \le S$ is elementary abelian of order
$p^2$. Since $H$ is strongly $p$-embedded, for $a \in A^\#$, we have $C_{O_{p'}(G)}(a) \le C_G(a)
\le H$. Therefore, from coprime action, $$O_{p'}(G) = \langle C_{O_{p'}(G)}(a)\mid a \in
A^\#\rangle\le \langle C_{G}(a)\mid a \in A^\#\rangle\le H.$$ Thus (v) follows from
(iv).\end{proof}

\begin{lemma}\label{main2} Suppose that $G$ is a group, $p$ is a prime and $H$ is a strongly
$p$-embedded subgroup of $G$. Set $\ov G = G/O_{p'}(G)$ and assume  further that $\ov H \neq \ov
G$. Then
\begin{enumerate}
\item $\ov H$ is strongly $p$-embedded in $\ov G$.
\item $F^*(\ov G)$ is a non-abelian simple group; and
\item if $G$ is a $\mathcal K$-group and $m_p(G) \ge 2$, then $(F^*(\ov G),p) \in \mathcal E$.
\end{enumerate}
\end{lemma}

\begin{proof} Let $S \in \syl_p(H)$. Then $\ov S \in \syl_p(\ov H)$. Choose $X$ a non-trivial
subgroup of $\ov S$. Then there exists $T \le S$ such that $X = \ov T$. Therefore, as $H$ is
strongly $p$-embedded in $G$, we have  $$H O_{p'}(G) \ge N_G(T)O_{p'}(G)= N_G(X)$$ by the Frattini
Argument. Hence $\ov H \ge N_{\ov G}(X)$. Since $\ov H \neq \ov G$, we conclude that $\ov H$ is
strongly $p$-embedded in $\ov G$. Hence (i) holds.

Since $\ov H$ is strongly $p$-embedded in $\ov G$ and $O_{p'}(\ov G)=1$, we have $F(\ov G) = 1$.
Assume that $E(\ov G)$ is not simple. Let $K_1$ be an arbitrary component of  $\ov G$ and let $K_2$
be a component with $K_1 \neq K_2$. Then, as $O_{p'}(\ov G) =1$, $p$ divides  $|K_2|$. Let $T = \ov
S \cap K_2$. Then $T \in \syl_p(K_2)$. Thus, as $\ov H$ is strongly $p$-embedded in $\ov G$ and
$[K_1,K_2]=1$, $K_1 \le N_{\ov G}(T) \le \ov H$. Since $K_1$ was chosen arbitrarily, we have that
$E(\ov G) \le \ov H$. Then Lemma~\ref{main1} (iv) implies that $E(\ov G) =1$ which is impossible. Thus
$F^*(\ov G)= F(\ov G) E(\ov G)=E(\ov G)$ is a simple group and (ii) holds

Finally (iii), follows directly from the definition of $\mathcal E$ and Proposition~\ref{SE-p2}.

\end{proof}

\begin{corollary}\label{FstarGsimple} If $m_p(H)\ge 2$ and $O_{p'}(H) = 1$, then $F^*(G)$ is a non-abelian simple group.
\end{corollary}
\begin{proof} We  have $m_p(H)=m_p(G)\ge 2$. Therefore, Lemma~\ref{main1}(iv) implies $O_{p'}(G) = \bigcap H^G \le  O_{p'}(H)=1$. Then
Lemma~\ref{main2}(ii) implies that $F^*(G)$ is a non-abelian simple group as claimed.
\end{proof}

%\begin{lemma}\label{controlp} Suppose that $H$ is a strongly $p$-embedded subgroup
%of $G$. Then $H$ controls $p$-fusion in $G$, and, in particular, $O^p(H)= H$ if and only if
%$O^p(G)= G$.
%\end{lemma}
%
%\begin{proof} We use \cite[Proposition 17.11]{GLS2} to see that $H$ controls $p$-fusion in $G$ and then \cite[Lemma 15.10
%(ii)]{GLS2} to obtain $O^p(H)= H$ if and only if $O^p(G)= G$.
%\end{proof}

\begin{lemma}\label{control1} Suppose that $p$ is an odd prime and  $H$ is strongly $p$-embedded of  $G$. Assume that for all involutions $t \in H$,  $p$ divides
$|C_H(t)|$. Then for all involutions $t \in H$, $t^G \cap H = t^H$.
\end{lemma}

\begin{proof} Suppose that $t \in H$ is an involution. Obviously $t^H \subseteq t^G \cap H$.
Assume that $g \in G$  and $t^g \in H$. Let $X \in \syl_p(C_H(t))$, $Y\in\syl_p(C_H(t^g))$ and note
that by assumption $X$ and $Y$ are non-trivial. As $H$ is strongly $p$-embedded in $G$, $X$ and $Y$
are Sylow $p$-subgroups of $C_G(t)$ and $C_G(t^g)$ respectively. Thus $X^{g}, Y \in
\syl_p(C_G(t^g))$ and so there exists $c \in C_G(t^g)$ such that $Y=X^{gc}$. Therefore $Y \le H
\cap H^{gc}$ and, since $H$ is strongly $p$-embedded in $G$, we get ${gc} \in H$. But then
$t^{g}=t^{gc} \in t^H$ and we are done.
\end{proof}

\begin{lemma}\label{control2} Suppose that $p$ is an odd prime and  $H$ is strongly $p$-embedded of  $G$. Assume that $t\in H$ is an involution, $t^G \cap H = t^H$ and $F$ is a subgroup of $H$ which contains $t$.
If $N_G(F) \not\le H$, then $C_G(t) \not \le H$.
\end{lemma}

\begin{proof} Aiming for a contradiction, we  assume that $C_G(t)= C_H(t)$.
Let $k\in N_G(F)\setminus H$. Then $t^k  \in F \le H$  and, as $t^G
\cap H = t^H$,  there exists $h \in H$ such that $t^{kh}=t$. But
then $kh \in C_G(t)=C_H(t) \le H$. Hence $k \in H$, which is a
contradiction. Thus $C_G(t) \not\le H$ as claimed.
\end{proof}

\section{Involutions in $H$}\label{sec4}

In this section we initiate the investigation of groups satisfying the hypotheses of
Theorems~\ref{MainTheorem0} and \ref{MainTheorem1}. Specifically throughout the remainder of this
article we assume that the following hypothesis holds.

\begin{hypothesis}\label{hypH}$p$ is an odd prime, $H$ is a strongly
$p$-embedded subgroup in $G$ and $H \cap K$ has even order for each non-trivial normal subgroup $K$
of $G$. Furthermore we assume that the following hold.\begin{enumerate}
\item $O_{p'}(H)=1$.
\item $m_p(C_H(t))$ is at least $2$ for each involution $t$ in $H$.
\item $N_G(T)$ is a $\mathcal {K}$-group, for all non-trivial
$2$-subgroups $T$ of $G$.
\item $O^2(G) = G$.
\end{enumerate}
\end{hypothesis}

Because of Hypothesis \ref{hypH} (i) and (ii), we can apply Corollary~\ref{FstarGsimple} to see
that $F^*(G)$ is a non-abelian simple group. The Frattini Argument then shows  $G = HF^*(G)$.
Furthermore, we have that $F^*(G) \cap H < F^*(G)$ and has even order.

\begin{lemma}\label{HStruct1} Suppose that $K$ is a subgroup of $G$, $K \not\le H$,  $m_p(H\cap K) \ge 2$ and $K$ is a $\mathcal
{K}$-group. Then the following hold.
\begin{enumerate}
\item $O_{p'}(K) \le H$.
\item $(F^*(K/O_{p'}(K)),p)\in \mathcal E$.
\item Either $O_{p'}(H\cap K) =O_{p'}(K)$, or $(F^*(K/O_{p'}(K)),p)=( \PSL_3(4),3)$ or $(\Fi_{22},5)$ and $|O_{p'}(H\cap K)/O_{p'}(K)|=2$.
\end{enumerate}
\end{lemma}

\begin{proof} As $H$ is strongly $p$-embedded in $G$,  we can use Lemma~\ref{main1} (ii) and
(v) to see that $O_{p'}(K) \le O_{p'}(H\cap K)\le H$. Then, as $K$ is a $\mathcal K$-group and
$m_p(H\cap K) \ge 2$, Lemma~\ref{main2}(iii) implies $(F^*(K/O_{p'}(K)),p) \in \mathcal E$  . Thus
(ii) holds. Part (iii) follows from (ii) and Corollary~\ref{OrredH3}.
\end{proof}

\begin{lemma}\label{not2-em}There exists an involution $t \in H$ such that $C_G(t) \not \le H$.
\end{lemma}

\begin{proof} Assume the lemma is false. Set $K = F^*(G)$. Then $K$ is a non-abelian simple group and $C_K(t) \le H \cap K$
for all involutions $t$ in $H\cap K$.  It follows from \cite[Lemma 17.13]{GLS2} that $K$ contains a
strongly $2$-embedded subgroup. Therefore, by  \cite{Bender}, we have  $K\cong \SL_2(2^a)$,
$\PSU_3(2^a)$ or ${}^2\B_2(2^a)$ for some $a$. In particular, $K$ is a $\mathcal {K}$-group and as
$m_p(G) \ge 2$, Proposition~\ref{SE-p2} delivers the contradiction. Thus there exists an involution
$t$ in $H$ such that $C_G(t) \not \le H$.\end{proof}

We can now present our first significant result.

\begin{lemma}\label{Ctnot-in} If $t \in H$ is an involution, then
\begin{enumerate}
\item $C_G(t) \not\le H$ and $C_H(t)$ is strongly $p$-embedded in $C_G(t)$;
\item $O_{p'}(C_G(t)) \le H$;
\item $(F^*(C_G(t)/O_{p'}(C_G(t))),p) \in
\mathcal E$; and
\item $O(C_G(t))= O(O_{p'}(C_G(t)))$.
\end{enumerate} \end{lemma}

\begin{proof}  By Lemma~\ref{not2-em},  there exists an
involution $t\in G$ such that $C_G(t)\not \le H$. Choose $t$ with $|C_H(t)|_2$ maximal. Set $K = C_G(t)$ and note that $C_H(t)$ is
strongly $p$-embedded in $C_G(t)$ by Lemma~\ref{main1} (ii). Then, by Hypothesis~\ref{hypH}(ii) and
Lemma~\ref{HStruct1}(i) and (ii), $O_{p'}(K) \le H$ and $(F^*(K/O_{p'}(K)),p)\in \mathcal E$. Let
$T \in \syl_2(H\cap K)$. We will show that $C_G(s) \not \le H$, for all involutions $s \in T$.
Assume first that $T \not \in \syl_2(K)$. Then $N_G(T) \not \le H$. Hence Hypothesis~\ref{hypH}
(ii) and Lemmas~\ref{control1} and \ref{control2} combine to give us that $C_G(s)\not \le H$ for
all involutions $s \in T$ which is our claim. So assume that $T \in \syl_2(K)$. Then, as $(H \cap
K)/O_{p'}(K)$ is strongly $p$-embedded in $K/O_{p'}(K)$, $(F^*(K/O_{p'}(K),p) \in \mathcal E$ and
$|K:H\cap K|$ is odd, Lemma~\ref{invs} implies that $F^*(K/O_{p'}(K))=K/O_{p'}(K) \cong \M_{11}$.
Let $F \leq T$ be such that $F \cap O_{p'}(K)$ is a Sylow 2--subgroup of $O_{p'}(K)$ and $|F : F
\cap O_{p'}(K)| = 2$. Assume that $s \in T$ is an involution with $C_G(s) \le H$. Then, as
$K/O_{p'}(K)\cong \M_{11}$ has exactly one conjugacy class of involutions, all the involutions of
$T$ are $G$-conjugate to involutions in $F$. Thus $\emptyset \neq s^G \cap F \le s^G \cap H = s^H$
by Lemma~\ref{control1} and we may therefore assume that $s \in F$. However, from the structure of
$\M_{11}$, $N_K(FO_{p'}(K))/O_{p'}(K) \not \le (H\cap K)/O_{p'}(K)$ and so we have $N_G(F) \not \le
H$ by the Frattini Argument. Thus Lemma~\ref{control2} implies that $C_G(s) \not \le H$ which is a
contradiction.  Hence $C_G(s) \not \le H$ for all involutions $s \in T$ and our claim is proved.

 Let $z $ be a $2$-central involution of $H$
which centralizes $t$. Then we may suppose that $z \in T$. It
follows that $C_G(z) \not \le H$. By maximality of $|C_H(t)|_2$ we may assume $t = z$ and so $T$ is a Sylow 2-subgroup of $H$ and this proves
(i).

Suppose that $t$ is an involution in $H$. Then by Hypothesis~\ref{hypH}(ii), $m_p(C_H(t)) \ge 2$
and $C_G(t) \not \le H$ by (i). Hence Lemma~\ref{HStruct1} (i) and (ii) gives parts (ii) and (iii).
Finally, as $O(C_G(t))$ is soluble, (iii) implies (iv).
\end{proof}

As a first application of Lemma~\ref{Ctnot-in}, we show that if $E(H)$ is non-trivial then it is
quasisimple.

\begin{lemma}\label{Eqs} If $E(H)\neq 1$, then $E(H)$ is quasisimple.
\end{lemma}

\begin{proof}  Set $E = E(H)$ and suppose that $E$ is the product of at least two components of $H$.
Assume that $L_1$ is a component of $H$ and let $L_2$ be the product of all the components of $H$
which are distinct from $L_1$. Then $E=L_1L_2$ and, since $O_2(H)\le O_{p'}(H)=1$ by
Hypothesis~\ref{hypH} (i), there exists an involution $t\in L_1 \setminus L_2$. We have $C_{G}(t)
\ge L_2$ and, as $p$ divides the order of each component in $L_2$, we have that $L_2/Z(L_2)$ is
isomorphic to a subnormal section of $C_H(t)/O_{p'}(C_H(t))$. Using Lemma~\ref{Ctnot-in} (ii) and
(iii), Proposition~\ref{Hstru1} and  the fact that $O_2(H)=1$, we deduce that  $(L_2,p)$ is
$(\Omega_8^+(2),5)$, $(\Alt(p),p)$ or $(\Alt(p)\times \Alt(p),p)$  where in the latter two cases we
have $p\ge 5$. Furthermore, we have that $O_p(C_H(t)/O_{p'}(C_H(t)))= 1$ and thus, as $O_p(H) \le
C_G(t)$, we have $O_p(H)=1$. By applying the above argument to $L_1$ with an involution taken from
$L_2$, we have that $(L_1,p)$ is $(\Omega_8^+(2),5)$ or $(\Alt(p),p)$. The possibilities for the
isomorphism type of $E$ when $p=5$ are thus $\Omega_8^+(2)\times  \Omega_8^+(2)$, $\Omega_8^+(2)
\times \Alt(5)\times \Alt(5)$, $\Omega_8^+(2) \times \Alt(5)$, $\Alt(5)\times \Alt(5)$ or $\Alt(5)
\times \Alt(5)\times \Alt(5)$ and, for $p>5$, we have $E\cong \Alt(p)\times \Alt(p)\times \Alt(p)$
or $\Alt(p)\times \Alt(p)$. Since $O_{p'}(H)=1=O_p(H)$ we have $E = F^*(H)$. Therefore, as $p \ge
5$, the structure of $E$ and the fact that the given components have no outer automorphisms of
order $p$ implies that $E(H)$ contains a Sylow $p$-subgroup of $H$.
 Select an involution
$d \in L_1L_2$ which projects non-trivially onto each component as a $2$-central involution. Then
$p$ does not divide $C_E(d)$ and we have a contradiction to Hypothesis~\ref{hypH}(ii). Thus $E$ is
 quasisimple as claimed.
\end{proof}

\section{Centralizers of involutions in $H$}\label{sec5}

In this section we work under Hypothesis~\ref{hypH} and aim to uncover the basic structure of the
centralizers of involutions from $H$ under the additional hypothesis that $H$ has no components and
that none of the involutions in $H$ are classical involutions (see
Definition~\ref{defclassicalinv}). We formalize the configuration we shall be investigating in the
following hypothesis.

\begin{hypothesis}\label{EH=1} Hypothesis~\ref{hypH} holds,
\begin{enumerate}\item $E(H) = 1$; and
\item  $G$ does not contain a classical involution.\end{enumerate}
\end{hypothesis}

Assume Hypothesis~\ref{EH=1}. Then, as $O_{p'}(H)=1$ by Hypothesis~\ref{hypH} (i), we have that
$F^*(H)= O_p(H)$. Set $Q= O_p(H)$ and note that, as $H$ is strongly $p$-embedded in $G$,we have $H=
N_G(Q)$ and, as $Q = F^*(H))$, $C_G(Q) \le Q$.

\begin{lemma}\label{thereisacomponent} Suppose that $s \in H$ is an
involution with $C_Q(s)>1$. Then there exists a normal component $L$ in $C_G(s)$ such that
$$LO_{p'}(C_G(s))/O_{p'}(C_G(s))= F^*(C_G(s)/O_{p'}(C_G(s))).$$ Furthermore, if $Z(L)$ has order
divisible by $2$, then $L/Z(L) \cong \PSL_3(4)$.
\end{lemma}

\begin{proof}
Since $C_Q(s)>1$ and $O_{p'}(C_G(s))\le H$ by Lemma~\ref{Ctnot-in}(ii), we have that
$[C_Q(s),O_{p'}(C_G(s))]=1$. Hence $F^*(C_G(s)) \not \le O_{p'}(C_G(s))$. Thus, as
$C_G(s)/O_{p'}(C_G(s))$ is an almost simple group by Lemma~\ref{Ctnot-in} (iii), we have $E(C_G(s))
\not \le O_{p'}(C_G(s))$ and so there is a component $L$ of $C_G(s)$ such that
$LO_{p'}(C_G(s))/O_{p'}(C_G(s))= F^*(C_G(s)/O_{p'}(C_G(s)))$. Since $L$ is the unique component of
$C_G(s)$ which has order divisible by $p$, we deduce that $L$ is normal in $C_G(s)$. If $2$ divides $|Z(L)|$,
then Lemma~\ref{SchurM}, the fact that $C_Q(s)O_{p'}(C_G(s))/O_{p'}(C_G(s))$ is a non-trivial
normal subgroup of $C_H(s)/O_{p'}(C_H(s))$ and Proposition~\ref{Hstru1}, together imply that $L
\cong \SL_2(p^a)$ or $L/Z(L) \cong \PSL_3(4)$. Since $s$ is not a classical involution by
Hypothesis~\ref{EH=1} (ii), we must have $L/Z(L) \cong \PSL_3(4)$ and thus the lemma is
established.
\end{proof}

\begin{lemma}\label{component} Let $u \in H$ be an involution with $C_Q(u) = 1$. Then there is a
$2$-component $L_u$  in $C_G(u)$ such that $L_u$ is not contained in $O_{p'}(C_G(u))$.
\end{lemma}

\begin{proof}  Let $W$ be a Sylow 2-subgroup of $O_{p'}(C_G(u))$. Then $C_G(u) =
N_{C_G(u)}(W)O_{p'}(C_G(u))$. Set $\overline{C_G(u)}= C_G(u)/O_{\{p,2\}'}(C_G(u))$ and let
$\overline{F} = F^*(\overline {C_G(u)})$. Since, by Lemma~\ref{Ctnot-in}(iii),
$C_G(u)/O_{p'}(C_G(u))$ is an almost simple group, we may assume that $\ov F \le
\overline{O_{p'}(C_G(u))}$ for otherwise there would be a component of $\ov{C_G(u)}$ not contained
in $\ov{O_{p'}(C_G(u))}$. Note that, as $C_H(u)$ is strongly $p$-embedded in $C_G(u)$, $O_{p}(\ov
F)=1$. Suppose that $C_{\overline{C_G(u)}}(\overline{W}) \not\leq \overline{O_{p'}(C_G(u))}$. Then
$\ov F$ is not a $2$-group and consequently $E(\overline{C_G(u)}) \neq 1$. Since $\overline W$
intersects each component of $\overline{C_G(u)}$ in a Sylow $2$-subgroup, we see that
$C_{\overline{C_G(u)}}(\overline W)$ normalizes each component of $\overline{C_G(u)}$. Therefore,
as $C_{\overline{C_G(u)}}(\overline{W})\ov F/\ov F$ is not soluble, the Schreier property of simple
groups delivers a contradiction to $\ov F$ being self-centralizing. Hence
$C_{\overline{C_G(u)}}(\overline{W})\leq \overline{O_{p'}(C_G(u))}$.

As $C_Q(u) = 1$, we have that $H = QC_H(u)$ and $Q$ is abelian. Assume that $s \in W$ is an
involution such that $s\neq u$ and $C_Q(s) \not \in \Syl_p(C_G(s))$. Since $s\neq u$, we have
$C_Q(s) >1$. Therefore Lemma~\ref{thereisacomponent} implies that  $C_G(s)$ has a normal component
$L_s$ with $L_s \not\le O_{p'}(C_G(s))$. Let $P \in \Syl_p(H\cap C_G(s))$. Then, as $H\cap C_G(s)$
normalizes $C_Q(s)$ and $uQ \in Z(H/Q)$, $C_Q(s) \le P$ and $[P,u]\le C_Q(s)$. Since $C_Q(s) < P$
and $C_Q(s)$ is normal in $H\cap C_G(s)$, Corollary~\ref{OrredH} implies that $C_Q(s)\le \Phi(P)$
or $(L_s,p) \cong (\PSL_2(8),3)$ or $({}^2\B_2(32),5)$. In the former case, we have that $[P,u] \leq C_Q(s) \leq \Phi(P)$ and so $P \le
C_G(u)$, which contradicts $[C_Q(s),u]\neq 1$. Thus we have the latter situation, and, as
$P\not=C_Q(s)$, this means that $C_Q(s)$ is  of order either $p$ or $p^2$. As $u \in N_{C_G(s)}(u)$, we see that $C_Q(s) \leq L_s$ and so $C_Q(s)$ is cyclic.
As $[P,u] \leq P$ also $[P,us] \leq Q$ and so also $C_Q(s)$ is not a Sylow $p$--subgroup of $C_G(us)$. This implies that also $C_Q(us)$ is cyclic of order $p$ or $p^2$.
 In particular,
since $\langle u, s \rangle$ acts on $Q$, we get, using coprime action, that $|Q/\Phi(Q)| = p^2$.
But then $H/Q$ is isomorphic to a subgroup of $\GL_2(p)$. It follows that $m_p(C_H(u))\le 1$ and
this contradicts Hypothesis~\ref{hypH}(ii) which asserts  that $m_p(C_H(u)) \geq 2$. We have seen
that for any involution $s \not= u$ in $W$ we have that $C_Q(s)$ is a Sylow $p$-subgroup of
$C_G(s)$.

As, by Hypothesis~\ref{hypH} (ii), $m_p(H/Q)= m_p(C_H(u)) \geq 2$, there is a non-cyclic abelian
$p$-subgroup  $E$ contained in $C_H(u)$. Suppose that $e \in E^\#$ and  $w \in C_W(e)$ is an
involution. Then $C_Q(w)\langle e \rangle$ is a $p$-subgroup of $C_G(w)$. If $w \neq u$, then
$C_Q(w) \in \syl_p(C_G(w))$. Since this is not the case, we see that for all $e \in E^\#$, $u$ is
the unique involution in $C_W(e)$. Therefore using Lemma~\ref{AutosE} we have $p=3$ and $W \cong
2^{1+4}_+$, $2^{1+6}_-$ or $2^{1+8}_+$. By Lemma~\ref{Ctnot-in} (iii) $F^*(C_G(u)/O_{p'}(C_G(u)))$
is a non-abelian simple group and so, since $C_G(W) \le O_{p'}(C_G(u))$, we have that
$N_G(W)/C_G(W) $ is non-soluble. In particular, we have that $W$ is not isomorphic to  $ 2^{1+4}_+$
which has soluble automorphism group. Let $s \in W$ be an involution with $s\neq u$. Then $C_W(s)
\cong 2 \times 2^{1+4}_-$ if $W \cong 2^{1+6}_-$ and $C_W(s) \cong 2 \times 2^{1+6}_+$ if $W\cong
2^{1+8}_+$. Since $u$ inverts $Q$, it inverts $C_Q(s)$. Thus $C_Q(s)$ is abelian and $C_Q(s)$
admits a faithful action of $C_W(s)/\langle s\rangle$. By Lemma~\ref{Ctnot-in} (iii),
$(F^*(C_G(s)/O_{p'}(C_G(s))),p)\in \mathcal E$. Since $C_Q(s) \in \Syl_p(C_G(s))$ and
$C_W(s)/\langle s\rangle$ is extraspecial of order $2^5$ or $2^7$, we have a contradiction to
Proposition~\ref{Hstru1} (i), (iii), (iv) and (v). This finally proves that there is a component in
$C_G(u)/O_{\{2,p\}'}(C_G(u))$.
 \end{proof}

We let $\mathcal I$ be the set of involutions $t\in H$ with $C_Q(t)>1$ and such that $|C_Q(t)| \geq
|C_Q(s)|$  for all involutions $s \in H$.  For any involution $x \in H$ with $C_Q(x)>1$, we let
$L_x$ denote the normal component of $C_G(x)$ described in Lemma~\ref{thereisacomponent}.

\begin{lemma}\label{Ltcentralizer}
Suppose that $s $ is an involution in $H$ with $C_Q(s) >1$. Then
$C_{C_G(s)}(L_sC_Q(s))=C_{C_G(s)}(L_s)= O_{p'}(C_G(s))$.
\end{lemma}

\begin{proof} We clearly have $[C_Q(s),O_{p'}(C_G(s))]=1$. Therefore, we also have
$[\langle C_Q(s)^{C_G(s)}\rangle,O_{p'}(C_G(s))]=1$. Hence, as $L_sC_Q(s)\le \langle
C_Q(s)^{C_G(s)}\rangle$,  $O_{p'}(C_G(s)) \le C_{C_G(s)}(L_sC_Q(s)) \le C_{C_G(s)}(L_s)$. On the
other hand, because $F^*(C_G(s)/O_{p'}(C_G(s)))= L_sO_{p'}(C_G(s))/O_{p'}(C_G(s))$,  we also have $
C_{C_G(s)}(L_s) \le O_{p'}(C_G(s))$. This proves the lemma.
\end{proof}

\begin{lemma}\label{alt} Let $s\in H$ be an involution  such that $C_Q(s)>1$. Then
 $C_H(s)/O_{p'}(C_G(s))$ is $p$-closed. In particular,
 $F^*(C_G(s)/O_{p'}(C_G(s)))$ is not
isomorphic to $\Alt(2p)$, $p \ge 5$ or to $\Fi_{22}$.\end{lemma}

\begin{proof}
As  $C_Q(s)O_{p'}(C_G(s))/O_{p'}(C_G(s))$ is a non-trivial normal $p$-subgroup of
$C_H(s)/O_{p'}(C_G(s))$, the assertion follows from Proposition~\ref{Hstru1}. \end{proof}

\begin{lemma}\label{CQt=1} Assume that $s,t\in H$ are involutions, $C_Q(t)=1$ and $t \not= s \in C_H(t)$. Then
$C_G(\langle s,t\rangle)$ is a $p'$-group. In particular, if $P \in \Syl_p(C_G(t))$, then
$m_2(N_{C_G(t)}(P)) \le 2$.
\end{lemma}

\begin{proof} Since $C_Q(t)=1$, we have $H= QC_H(t)$ and $C_H( s ) =
C_Q(s)C_H(\langle s,t\rangle)$. Let $P\in \syl_p(C_H(\langle s,t\rangle))$.  Then $C_Q(s)P\in
\syl_p(C_H(s)) \subseteq \syl_p(C_G(s))$ by Lemma~\ref{main1} (i). Since $t$ normalizes $C_Q(s)P$
and $t$ centralizes $P$ but inverts $C_Q(s)$, we must have that $C_Q(s) \not \le \Phi(C_Q(s)P)$.
Assume that $P \neq 1$. Then $C_Q(s)\Phi(C_Q(s)P) \neq C_Q(s)P$. In particular,
$C_H(s)/O_{p'}(C_H(s))$ does not act irreducibly on $C_Q(s)P/\Phi(C_Q(s)P)$. Since
$C_H(s)/O_{p'}(C_H(s))$ is $p$-closed by Lemma~\ref{alt}, Corollary~\ref{OrredH} implies that
$(L_s,p)$ is either $(\PSL_2(8),3)$ or $({}^2\B_2(32),5)$. As $t$ inverts $C_Q(s)$, we see that $C_Q(s) \leq L_s$
and then $C_Q(s)$ is a cyclic group of
order dividing $p^2$. Since $Q= C_Q(st)C_Q(s)$, we have that $|Q/\Phi(Q)| \le p^2$. This means that
$H/Q$ is isomorphic to a subgroup of $\GL_2(p)$ and implies that $C_H(t)$ has cyclic Sylow
$p$-subgroups. This of course contradicts Hypothesis~\ref{hypH}(ii). Thus $C_G(\langle s,t\rangle)$
is a $p'$-group as claimed.
\end{proof}

\begin{lemma}\label{cyclic} Suppose that $t  \in \mathcal I$ or $C_Q(t)=1$. Then
$O_{p'}(C_G(t))$ has cyclic Sylow $2$-subgroups. In particular, $O_{p'}(C_G(t))$ has a normal
$2$-complement.
\end{lemma}

\begin{proof}
 Assume that $F=\langle
t,s\rangle$ is a fours group in $O_{p'}(C_G(t))$. Then $Q = C_Q(t)C_Q(s)C_Q(ts)$.

Assume further that $t \in \mathcal I$. Then, as $F \le O_{p'}(C_G(t)) \le H$ and $[F,C_Q(t)] \leq
[O_{p'}(C_G(t)),C_Q(t)] = 1$, the maximality of $C_Q(t)$ implies that $C_Q(t) = C_Q(F)= C_Q(s)$.
But then $Q= C_Q(t)$, which contradicts $C_G(Q) \le Q$. Hence no such $F$ exists and so
$O_{p'}(C_G(t))$ has either quaternion or cyclic Sylow $2$-subgroups.

Assume  that $C_Q(t) =1$. Then, by Lemma~\ref{component}, there is a $2$-component $L$ in $C_G(t)$
with $LO_{p'}(C_G(t))/O_{p'}(C_G(t))= F^*(C_G(t)/O_{p^\prime}(C_G(t))))$. Now we have $[L,F]\le
O(C_G(t))$. In particular, $C_H(F)$ contains a Sylow $p$-subgroup $P$ of $L$  and, as $P$ is
non-trivial, this contradicts Lemma~\ref{CQt=1}. Thus
 once again, as $t$ is not a classical involution, we have that $O_{p'}(C_G(t))$ has cyclic Sylow $2$-subgroups
and a normal $2$-complement.

Since, by
Hypothesis~\ref{EH=1} (ii), $t$ is not a classical involution, we deduce that $O_{p'}(C_G(t))$ has
cyclic Sylow $2$-subgroups and a normal $2$-complement.
\end{proof}

\section{Centralizers of involutions with $F^*(C_G(t)/O_{p'}(C_G(t)))\not \cong \PSL_2(p^f)$}\label{sec7}

In this section we will prove the following theorem.

\begin{theorem}\label{L2pthm} Assume that Hypotheses~\ref{hypH} and \ref{EH=1}
hold. Then  either
\begin{enumerate}
\item there exists $t \in \mathcal I$ with  $F^*(C_G(t)/O_{p'}(C_G(t)))\cong \PSL_2(p^f)$
for some $f > 1$; or \item  for all involutions $t \in H$,  $C_Q(t)=1$ and
$F^*(C_G(t)/O_{p'}(C_G(t)))\cong \PSL_2(p^f)$ with $f > 1$ and  $p \equiv 3 \pmod 4$.
\end{enumerate}
\end{theorem}

We prove Theorem~\ref{L2pthm} via a sequence of lemmas the first of which shows that the theorem
holds if $\mathcal I$ is empty. We continue to use the notation of Section~\ref{sec5}. In
particular, $Q = O_p(H)$.

\begin{lemma}\label{E2} If $\mathcal I = \emptyset$, then Theorem~\ref{L2pthm} (ii) holds.
\end{lemma}

\begin{proof} Suppose that  $\mathcal I =\emptyset$. Then $C_Q(t)=1$ for all involutions $t \in H$ and, in particular,
$m_2(H)=1$. Therefore the Sylow $2$-subgroups of $H$  are either cyclic or quaternion. Let $t \in
H$ be an involution. Then $C_Q(t)=1$. By Lemma~\ref{component} there is a normal  $2$-component $L$
of $C_G(t)$ with $L \not \le O_{p'}(C_G(t))$. Note that Lemma~\ref{Ctnot-in} (ii) and (iv) imply
that $O(L) \le H$. Since $H$ has cyclic or quaternion Sylow $2$-subgroups, so does $C_H(t)$. If $t
\in L$, then $tO(L) \in Z(L/O(L))$ is non-trivial. Therefore, as the Sylow $2$-subgroups of $L\cap
H$ are either cyclic or quaternion, Corollary~\ref{CExt} implies that $L/O(L) \cong \SL_2(p^f)$ for
some $f$. But then $G$ contains a classical involution, and this contradicts Hypothesis~\ref{EH=1}
(ii). Hence  $t \not \in L$. Since $t$ is the unique involution in $C_H(t)$, we deduce that $|H\cap
L|$ is odd. Now Lemma~\ref{invs2} shows that $L/O(L) \cong \PSL_2(p^f)$ where $f>1$ is odd and $p
\equiv 3 \pmod 4$.
 \end{proof}
%
%Let $E = \langle F,t\rangle$ where $F$ is a fours subgroup of $R\cap
%L$. Then, as $C_{C/O_{\{p,2\}'}(C)}(E) =
%ET_0O_{\{p,2\}'}(C)/O_{\{p,2\}'}(C)$, we then see that  $C_G(E) =
%C_{C}(E) \le KET_0$. Thus $ET_0 \in \Syl_2(C_G(E))$.

We may now assume that $\mathcal I \not=\emptyset$ and so we begin to restrict the possibilities
for the structure of $C_G(t)$, where $t \in H$ is an involution with $C_Q(t) > 1$. By
Lemma~\ref{Ctnot-in} $C_G(t) \not\le H$ and $O_{p'}(C_G(t)) \le H$ and by
Lemma~\ref{thereisacomponent}, there is a normal component $L_t$ of $C_G(t)$ such that $L_t\not \le
O_{p'}(C_G(t))$. By Lemma~\ref{Ctnot-in} (iii) we then have $(L_t/Z(L_t), p) \in \mathcal E$. Our
aim now is to show that there is an involution $t\in \mathcal I$ with $L_t \cong \PSL_2(p^f)$.

\begin{lemma}\label{E3}  If $t\in H$ is an involution with $C_Q(t)>1$, then
$(L_t/Z(L_t),p)$ is one of $(\PSL_2(p^a),p)$, with $p$ arbitrary and $a\ge 2$,
$({}^2\G_2(3^{2a-1}),3)$ with $a\ge 2$, $(\PSL_2(8),3)$, $(\PSL_3(4),3)$ or $({}^2\B_2(32),5)$.
\end{lemma}

\begin{proof} Assume the claim is false. Then, as   $(L_t/Z(L_t),p)\in \mathcal E$ and $O_p(H \cap L_t) \not= 1$, we have that $(L_t/Z(L_t),p)$  is one
of $(\PSU_3(p^{a}),p)$ with $a\ge 2$,  $(\M_{11},3)$, $(\McL,5)$, $({}^2\F_4(2)',5)$ or
$(\J_4,11)$. Furthermore, by Lemma~\ref{thereisacomponent}, $Z(L_t)$ has odd order. In particular,
in each of the cases we need to consider there is an involution $s\in L_t$ such that $C_{L_t}(s)$
has structure as indicated in the Table~\ref{Table1}.
\begin{table}
\begin{tabular}{|cc|}\hline
$(L_t/Z(L_t),p)$&$C_{L_t/Z(L_t)}(s)$\\\hline
$(\M_{11},3)$&$\GL_2(3)$\\
$ ({}^2\F_4(2)',5)$&$2^{9}{}.\Frob(20)$\\
 $(\McL,5)$&$2\udot
 \Alt(8)$\\$(\J_4,11)$&$2^{1+12}_+.3\udot\Aut(\M_{22})$\\
$(\PSU_3(p^a),p)$&$O^{p'}(C_{L_t/Z(L_t)}(s))\cong
\SL_2(p^a)$\\\hline
\end{tabular}
\caption{Centralizers of $s$ in $L_t/Z(L_t)$}\label{Table1}\end{table} In particular, we note that
in each case $s \in O^{p'}(O^p(C_{L_t}(s)))$. By Lemma~\ref{main1} (i) and (ii), $H \cap L_t$ is
strongly $p$-embedded in $L_t$. Since in each of the groups under consideration as $L_t$, the
centralizer of $s$ has order divisible by $p$ we may assume that $s \in H\cap L_t$. Since $Q \cap
L_t= C_Q(t)$ is normalized by $C_H(t)$, Proposition~\ref{Hstru1} implies that we either have that
$C_Q(t) \in \syl_p(L_t)$ or $p=11$, $L_t \cong J_4$ and $C_Q(t) $ is the centre of the Sylow
$11$-subgroup of $H\cap L_t$ or $L_t \cong \PSU_3(p^a)$. In each case, we have $C_Q(s) \ge
C_{C_Q(t)}(s)>1$. Therefore Lemma~\ref{thereisacomponent} indicates that $C_{G}(s)$ has a normal
component $L_s$ such that $L_s \not \le O_{p'}(C_G(s))$. Since $[O_{p'}(C_G(s)),L_s]=1$, we have
that $L_s= O^{p'}(O_{p'}(C_G(s))L_s)$. Hence, we either have $O^{p'}(O^{p}(C_G(s)))= L_s$, or $p$
divides $|C_G(s)/O_{p'}(C_G(s))L_s|$. In the latter case, Proposition~\ref{SE-p2} shows that
$C_G(s)/O_{p'}(C_G(s)) \cong \PSL_2(8):3$ and $p=3$ or $C_G(s)/O_{p'}(C_G(s)) \cong {}^2\B_2(32):5$
and $p=5$. In this case we have $L_s= O^{p'}(O^{p}(C_G(s)))$ as well. It follows that in all cases
$s \in L_s$, as $s \in O^{p'}(O^p(C_{L_t}(s)))$. Therefore Lemma~\ref{thereisacomponent} implies  that $L_s/Z(L_s) \cong \PSL_3(4)$ and
that $p=3$. Furthermore, the Sylow $3$-subgroups of $C_G(s)$ have order $9$ and so after consulting
Table~\ref{Table1} we get $L_t \cong \M_{11}$, $\PSU_3(3)$ or $\PSU_3(9)$.   If $t \in
O_{p'}(C_G(s))L_s$, then either $C_G(\langle t,s\rangle)$ would have order coprime to $3$ or would
contain a component $K$ such that $K/Z(K)$ is isomorphic to $\PSL_3(4)$, from the choice of $s$ we
see that both of these scenarios are impossible. Hence $t$ induces an outer automorphism of $L_s$
and so with Lemma~\ref{L34} we get that $C_{L_s/Z(L_s)}(t) \cong 3^2.\mathrm \Q_8$, $\PSL_2(7)$ or
$\Alt(5)$. Since $C_{L_s}(t)$ projects
 to a subgroup of $C_G(t)/O_{p'}(C_G(t))$ which is normal in the centralizer of the
image of $s$ in $C_G(t)/O_{p'}(C_G(t))$, this then gives us our final contradiction.
\end{proof}

\begin{lemma}\label{notree} If $t\in \mathcal I$, then
$L_t/Z(L_t)$ is not isomorphic to ${}^2\G_2(3^{2n+1})$, $n \geq 1$.
\end{lemma}

\begin{proof}
Suppose that $L_t/Z(L_t) \cong {}^2\G_2(3^{2n+1})$ with $n \ge 1$. Note that by \cite[Theorem
2.5.12]{GLS3} the outer automorphisms of $L_t$ have odd order and $Z(L_t)=1$. Let $R \in
\syl_2(L_tO_{p'}(C_G(t)))= \syl_2(C_G(t))$ be such that $R \cap H \in \Syl_2(C_H(t))$ and set  $T =
R \cap O_{p'}(C_G(t))$ and $E = T\cap L_t$. Then, by Lemma~\ref{cyclic}, $T$ is a cyclic group and,
as $L_t\cong {}^2\G_2(3^{2n+1})$, $E$ is elementary abelian of order $8$. Since
$[L_t,O_{p'}(C_G(t))]=1$, $R=T\times E$ is abelian and so $R \in \syl_2(C_G(R))$. By
Lemma~\ref{wcl}, $t$ is conjugate in $G$ to an element $t^x \in R \setminus \{t\}$ and so
Lemma~\ref{cyc2} implies that $T$ has order $2$. Therefore $R$ is elementary abelian of order $16$.

From the structure of ${}^2\G_2(3^{2n+1})$ we have  $N_{L_t}(E) /E \cong \Frob(21)\cong
N_{L_t^x}(E^x)/E^x$. Therefore  $N_G(R)/C_G(R)$ contains distinct subgroups isomorphic to
$\Frob(21)$ and  $N_G(R)$ has orbits of lengths either $8$ and $7$ or $15$ on $R^\#$. As $\GL_4(2)$
has no subgroups of order $3^2 \cdot 5 \cdot 7$ it is impossible for $N_G(R)$ to be transitive on
$R^\#$. Hence $|t^{N_G(R)}|= |t^G \cap R|=8$, $N_G(R)/C_G(R) \cong 2^3.\Frob(21)$ or $\PSL_3(2)$
and $E $ is the unique normal subgroup of $N_G(R)$ of order $2^3$. Let $s \in E^\#$ be chosen so
that $s$ normalizes a Sylow $3$-subgroup $D$ of $C_H(t)$, then $s$ and $t$ are not $G$-conjugate
and $s \in H$. Furthermore $C_{L_t}(s) \cong 2 \times \PSL_2(3^{2n+1})$ and
 $C_{N_{G}(R)}(s)$ has
non-abelian Sylow $2$-subgroups. In particular, as $3$ divides $|C_{G}(\langle s,t \rangle)|$,
Lemma~\ref{CQt=1} implies that $C_Q(s)\neq 1$. By Lemma~\ref{thereisacomponent}, there is a normal
component $L_s$ in $C_G(s)$ with $L_s \not \leq O_{p'}(C_G(s))$. We then have that $t$ induces an
automorphism of $L_s$ which centralizes a subgroup of $L_s$ which is isomorphic to
$\PSL_2(3^{2n+1})$. By Lemmas~\ref{L34} and \ref{E3}, we see that $L_s/Z(L_s)\cong
{}^2\G_2(3^{2m+1})$ for some $m \ge 1$ or $L_s/Z(L_s)\cong \PSL_2(3^m)$ for some $m \ge 2n+1$. In
the former case we have $C_{L_s}(t) \cong C_{L_t}(s)$ and so we infer that $L_t \cong L_s$ (but
perhaps $s \not \in \mathcal I$). Let $B \in \syl_2(C_G(s))$ be chosen so that $R \le B$. Then, as
$\Out(L_s)$ has odd order and $Z(L_s)=1$, $B = (B\cap O_{p'}(C_G(s)))\times (B\cap L_s)$ where
$B\cap L_s$ is elementary abelian of order $8$. Notice that $t \in B$ and $C_{B}(t) \ge \langle
s\rangle(B\cap L_s)$. Since $R\in \Syl_2(C_G(t))$ and $|R|=2^4$, we get that $R = \langle s\rangle
(B\cap L_s)$. However, we then have $ N_{L_s}(B\cap L_s) \le N_G(R)$ and so $N_{L_s}(B\cap L_s)$
leaves both $E$ and $B\cap L_s$ invariant. Hence $E=(B\cap L_s)$ and consequently $s \in E \subset
L_s$ which is a contradiction. Assume that $L_s \cong \PSL_2(3^m)$.
 Then, by Proposition~\ref{SE-p2} and Lemma~\ref{Ctnot-in},
$L_s=O^{3'}(C_G(s))$. Since $O^{3'}(C_{N_G(R)}(s))C_G(R)/C_G(R) \cong \Alt(4)$  and
$|O^{3'}(C_{N_G(R)}(s)) \cap R|=2^2$ or $2^3$, we see that $O^{3'}(C_{N_G(R)}(s))$ contains a non cyclic abelian
group of order at least 8. Hence $L_s$ contains a noncyclic abelian group of order 8. So we have a contradiction to the Sylow $2$-subgroup structure
of $\PSL_2(3^m)$. This completes the proof of the lemma.
\end{proof}

\begin{lemma}\label{l28} If $t\in \mathcal I$, then
$(L_t/Z(L_t),p)$ is not  either $(\PSL_2(8),3)$ or $({}^2\B_2(32),5)$.
\end{lemma}

\begin{proof} Assume $({L_t}/Z({L_t}),p)$ is   either
$(\PSL_2(8),3)$ or $({}^2\B_2(32),5)$. Then $C_G(t)/O_{p'}(C_G(t)) \cong \PSL_2(8):3$ or
${}^2\B_2(32):5$ respectively. In both cases, there is an involution $s \in {L_t}$ which normalizes
$C_Q(t)$. Thus $s\in H$. Let $R$ be a Sylow 2-subgroup of $C_G(t)$ which contains $\langle
s\rangle$. Then $ R = T \times (R \cap {L_t})$, where $T$ is a Sylow 2-subgroup of
$O_{p'}(C_G(t))$. By Lemma~\ref{cyclic} we have that $T$ is cyclic. By Lemma~\ref{wcl}, $t$ is
conjugate to some element $t^x$ in $R \setminus \{t\}$. Since $Z(R \cap {L_t})$ is elementary
abelian, Lemma~\ref{cyc2} implies that $T = \langle t \rangle$ has order $2$.

Assume that $({L_t},p) = (\PSL_2(8),3)$. Then $R$ is elementary abelian of order $2^4$ and
$N_{C_G(t)}(R)/C_{C_G(t)}(R) \cong \Frob(21) \cong N_{C_G(t^x)}(R)/C_{C_G(t^x)}(R)$. Just as in the
proof of Lemma~\ref{notree}, we get $N_G(R)/C_G(R) \cong 2^3.\Frob(21)$ or $\PSL_3(2)$ and the
$N_G(R)$-orbits on $R$ have lengths $8$ and $7$ with $t$ in the orbit of length $8$. Since $s$
normalizes a Sylow $3$-subgroup of $C_H(t)$ and since $s$ is contained in a unique Sylow
$2$-subgroup of $L_t$, we can choose  $P \in \Syl_3(N_G(R)\cap H)$  so that $t \in C_R(P)$. Let $F
= \langle t,r\rangle=C_R(P)$. We may choose notation so that $|r^{N_G(R)}|=7$. Since $3$ divides
$C_G(F)$, Lemma~\ref{CQt=1} implies that $C_Q(r)>1$. Furthermore, as $C_{N_G(R)}(r)/R$ contains a
subgroup isomorphic to $\Alt(4)$, we have $|O^{3'}(C_{N_G(R)}(r))| \ge 2^4.3$.

By Lemmas~\ref{thereisacomponent} and \ref{E3}, we have that $L_r/Z(L_r)$ is isomorphic to one of
$\PSL_2(8)$, $\PSL_3(4)$, $\PSL_2(3^m)$, $m \ge 2$, or ${}^2\G_2(3^{2n+1})$, $n \ge 1$. If $L_r=
O^{3'}(C_G(r))$, then $L_r \ge O^{3'}(C_{N_G(R)}(r))$. As ${}^2\G_2(3^{2n+1})$ and $\PSL_2(8)$ have
Sylow $2$-subgroups of order $8$ and  in $\PSL_2(3^m)$ the normalizer of a four group is isomorphic
to either $\Alt(4)$ or $\Sym(4)$, these possibilities cannot arise. If $L_r/Z(L_r) \cong
\PSL_3(4)$, then $|C_Q(r)|= 3^2$. Now in $C_G(t)$, we see that $|C_{C_G(t)}(r)|_3= 3$ and
$C_{C_G(t)}(r) O_{3'}(C_G(t))/ O_{3'}(C_G(t))$  is soluble. Since $t$ induces an automorphism of
$L_r$, we now deduce a contradiction from Lemma~\ref{L34}. Thus $L_r \not= O^{3'}(C_G(r))$. Then
$L_r \cong \PSL_2(8)$. Since $R$ is abelian, $R \le C_G(r)$. Let $D \in \Syl_2(C_G(r))$ with $R \le
D$. Then $D = (D\cap O_{3'}(C_G(r))) \times (D \cap L_r)$. Thus $C_D(t) \ge \langle r\rangle \times
(D \cap L_r) \le Z(D)$. It follows that $D \le C_G(t)$ and so $R=D$. Since $N_G(R)$ normalizes a
unique subgroup of order $2^3$ in $R$, it follows that $r \in R\cap L_t  =R \cap L_r$ which is of
course impossible. This shows that  $L_r \not\cong \PSL_2(8)$ and so $L_t/Z(L_t) \not \cong
\PSL_2(8)$.

Assume that $({L_t},p) = ({}^2\B_2(32),5)$. Let $P \in \syl_5(C_H(t))$ and let $s \in L_t$ be such
that $|C_G(s) \cap P |=5$. Then $s \in H$ and $|C_{L_tP}(s)| = 2^{10}\cdot 5$. Therefore
$O^{5'}(C_G(s)) $ has a subgroup of order  $ 2^9\cdot 5$ and $2$-rank $5$. By
Lemmas~\ref{thereisacomponent}, \ref{alt} and \ref{E3}, we have that $L_s/Z(L_s)$ is isomorphic to
${}^2\B_2(32)$ or  $\PSL_2(5^m)$ some $m\ge 1$. In the latter case, Proposition~\ref{SE-p2} implies
that $L_s= O^{5'}(C_G(s))$. Since $L_s/Z(L_s)$ has $2$-rank $2$, we have a contradiction. Thus $L_s
\cong {}^2\B_2(32)$.  Since $s \in Z(R)$, $R \le C_G(s)$. Let $D \in \Syl_2(C_G(s))$ with $R \le
D$. Then $D = (D\cap O_{5'}(C_G(s))) \times (D \cap L_s)$. As $t \in \Omega_1(D)$ and $\Omega_1(D)
\le  (D\cap O_{5'}(C_G(s)))Z(D \cap L_s)$, we have that $D \cap L_s \le R$. Therefore $R= (D\cap
L_s) \langle s\rangle$. But then $s\in R' = (D\cap L_s)'$ and this is of course impossible.
\end{proof}

\begin{lemma}\label{l2q}  There exists $t\in \mathcal I$ such that
$(L_t/Z(L_t),p)$ is  $(\PSL_2(p^f),p)$.
\end{lemma}

\begin{proof} Suppose the lemma is false and let $t \in \mathcal I$. Then by  Proposition~\ref{SE-p2} and
Lemmas~\ref{thereisacomponent}, \ref{alt}, \ref{E3},\ref{notree} and \ref{l28}, we may assume that
$(L_t/Z(L_t),p)$ is $(\PSL_3(4),3)$.

Let $R$ be a Sylow 2-subgroup of $C_G(t)$ and set $T = R \cap O_{3'}(C_G(t))$. By
Lemma~\ref{cyclic}, $T$ is cyclic. By Proposition~\ref{Hstru1} (iv), $(L_t \cap H)/Z(L_t) \cong
3^2:\Q_8$ and so  $|C_Q(t)O_{3'}(C_H(t))/O_{3'}(C_H(t))|=|C_Q(t)|=9$. We also remark that as $T$ is
cyclic, $Z(L_t)$ is cyclic and if $Z(L_t)$ is non-trivial then $t \in Z(L_t)$.
\\
\\
We remark that the normalizer of an elementary abelian group $E$ of order $2^4$ in $L_t/Z(L_t)$ has
shape $2^4:\Alt(5)\cong 2^4:\SL_2(4)$. Further $N_{L_t/Z(L_t)}(E)$  acts transitively on
$E^\sharp$. Suppose that $|Z(L_t)| = 2$. Then because of the  transitive action of $L_t$ on
$E^\sharp$, we have that $E$ lifts to an elementary abelian group $F$ of order 32. In particular
involutions lift to involutions. Suppose that $N_{L_t}(F)$ acts decomposably on $F$. Then
$N_{L_t}(F)/E$ is a direct product of a group of order two by $\PSL_2(4)$. But then a Sylow
2--subgroup splits over $Z(L_t)$, which is not possible. Hence $N_{L_t}(F)$ acts indecomposably on
$F$. Assume now that $Z(L_t)$ is cyclic of order four. If the preimage $F_1$ of $E$ would be
abelian, then it is of type $(4,2,2,2,2)$. Now $|\Omega_1(F_1)| = 32$ and $|\Phi (F_1)| = 2$. But
then $N_{L_t}(F_1)$ cannot act indecomposably on $F_1/\Phi (F_1)$, which contradicts the fact just
proved for $|Z(L_t)| = 2$. Hence $F_1$ is nonabelian and so $E$ lifts to a special group of
symplectic type.
%From \cite[page
%24]{Atlas} we see that if $|Z(L_t)| \le 2$, then the elementary abelian subgroups of order $2^4$ in
%$L_t/Z(L_t)$ lift to elementary abelian subgroups of order $2^5$ while if $Z(L_t)$ is cyclic of
%order $4$ they left to special groups of symplectic type. Furthermore, we remark that the
%normalizer of an elementary abelian group of order $2^4$ in $L_t/Z(L_t)$ has shape
%$2^4:\Alt(5)\cong 2^4:\SL_2(4)$.
\\
\\
Assume  that $TL_t$ contains an elementary abelian subgroup $F$ of order 32. Then we have that
$|Z(L_t)| \le 2$.  We select $F$ such that $F\cap H$ has order exactly $4$ and set $F \cap H =
\langle s,t\rangle$ where $s \in L_t$ with $s\neq t$. As $s \in L_t$, $s$ and $st$ invert $C_Q(t)$.
As $[t,Q] \not= 1$ there is $x\in \{s,st\}$ such that $C_Q(x) \not= 1$. If, additionally, $x \not\in \mathcal I$,
then, as $|C_Q(t)|=9$, $|C_Q(x)| = 3$ and, by Lemma~\ref{thereisacomponent} and
Corollary~\ref{OrredH}, $L_x \cong \PSL_2(8)$ and $C_G(x)/O_{3'}(C_G(x)) \cong \PSL_2(8):3$. Since
$C_G(x)/O_{3'}(C_G(x)) \cong \PSL_2(8):3$ and $F \le C_G(x)$, we have that $|F \cap
O_{3'}(C_G(x))|\ge 4$. Since $O_{3'}(C_G(x)) \le H$, we get $F \cap O_{3'}(C_G(x))= F\cap H$. But
then $t \in O_{3'}(C_G(x))$ and as $[O_{3'}(C_G(x)),C_Q(x)]=1$, we have $C_Q(x) \le C_Q(t)$ whereas
$x$ inverts $C_Q(t)$. Hence, if $C_Q(x) \neq 1$, then $x \in\mathcal I$. Hence  we also have
$L_x/Z(L_x) \cong \PSL_3(4)$. Let $U$ be a Sylow 2-subgroup of $O_{3'}(C_G(x))$. Then, as $x \in
\mathcal I$, $U$ is cyclic and $U$ commutes with $L_x$. Since $TF \cap O_{3'}(C_G(x))$ is cyclic
and $TF$ is abelian, we have that $|TF O_{3'}(C_G(x))/O_{3'}(C_G(x))|= 2^4$ and hence $F \le UL_x$
by Lemma~\ref{L34}. In particular, $U$ and $F$ commute.  As $TF$ is a Sylow $2$-subgroup of
$C_G(F)$, we have $TF = UF$ and $|U|= |T|$. Assume that $|T|>2$. Then, on the one hand $t$ is the
unique involution in $\Phi(TF)$, while on the other hand, $x$ is the unique involution in $UF$ and
so, as $x\neq t$,  we must have $T = \langle t \rangle$ has order $2$. As $N_{L_t}(F)$ centralizes
only $t$ in $F$ and $N_{L_x}(F)$ centralizes only $x$ in $F$, we have $N_G(F) \not\le C_G(t)$. Of
course we also have that $N_G(F)/C_G(F)$ is isomorphic to a subgroup of $\GL_5(2)$. From the
structure of $L_t$, we know that $N_{C_G(t)}(F)$ induces orbits of length 1,15,15 on $F^\sharp$ and
so $t$ has 31 or 16 $N_G(F)$-conjugates in $F$. If all involutions in $F$ are conjugate, we get
that $N_G(F)/C_G(F)$ has order $2^2 \cdot 3 \cdot 5 \cdot 31$ if $N_{C_G(t)}(F)$ induces $\Alt(5)$
on $F$ or $2^3 \cdot 3 \cdot 5 \cdot 31$ if $N_{C_G(t)}(F)$ induces $\Sym(5)$ on $F$. As in both
cases, by Sylow's Theorem, the Sylow $31$-subgroup would be normal, this contradicts the structure
of the normalizer of a Sylow 31-subgroup of $\GL_5(2)$. So $t$ has $16$ $N_G(F)$-conjugates in $F$.
Since $x\in \mathcal I$, we may argue in a similar way to see that there are $16$
$N_G(F)$-conjugates of $x$ in $F$.  Recall that $x$ was chosen in $\{s,st\}$ and $s \not\cong st$  in $L_t$,  at
least one of $s$ and $st$
is not $N_G(F)$ conjugate to $t$. It follows that either $C_Q(s)=1$ or $C_Q(st)=1$. Assume that $y
\in\{s,st\}$ is such that $C_Q(y)=1$. Then $Q= C_Q(t)C_Q(x)$ and, as $|C_Q(x)|=|C_Q(t)|=9$,  $Q$ is
elementary abelian of order $3^4$ and $H/Q$ is isomorphic to a subgroup of $\GL_4(3)$. By
Lemmas~\ref{component} and \ref{cyclic} there is a $2$-component $L$ in $C_G(y)$ and
$O_{3'}(C_G(y))$ has a normal $2$-complement and a cyclic Sylow $2$-subgroup. Since $O(L) \le H$
and as $|O(L)|$ has order coprime to both $2$ and $3$ and divides $|\GL_4(3)|$, we have that
$O(L)$ is cyclic of order dividing $65$. Since $p=3$ and $C_H(y)$ is strongly $3$-embedded in
$C_G(y)$, we have that $|O_3(C_H(y)/O(C_G(y)))| \ge 9$ from Proposition~\ref{Hstru1}. But then, as
$O(L)$ is cyclic, $O_3(C_H(y)) \not=1$, but this is impossible as $H= C_H(y)Q$, $Q=O_3(H)$ and $Q
\cap C_H(y)=1$. This contradiction shows that $TL_t$ does not contain an elementary abelian
subgroup of order $32$. Therefore  $L_t$ is a central extension of $\PSL_3(4)$ by a cyclic group of
order 4. In particular, $t \in L_t$ and $L_t$ has precisely two conjugacy classes of involutions.

Assume that there is an $s \in t^G \cap C_G(t)\setminus L_tO_{p'}(C_G(t))$. Then, by
Lemma~\ref{L34} we may assume that $C_{C_Q(t)}(s) \not= 1$ which means that $s \in H$. If $C_Q(t) =
C_Q(s)$, then, as $t \in  \mathcal I$, we get $[Q,s] = 1$, a contradiction. So we have that
$|C_Q(s) \cap C_Q(t)| = 3$. As $t \in \mathcal I$, we see by coprime action that $|Q| = 3^4$.
Furthermore, we have that $Q$ is elementary abelian. Now choose $u \in L_t $ an involution such
that $u$ inverts $C_Q(t)$. By Proposition~\ref{Hstru1} (iv), we have that $u$ is a square in $H$.
In particular, $u$ acts on $Q$ as an element of $\PSL_4(3)$. Since $u \sim_{L_t} ut$,
Lemma~\ref{control1} implies that $u \sim_H ut$. Therefore both $C_Q(u)$ and $C_Q(ut)$ are
non-trivial. Because $u$ and $ut$ act on $Q$ as elements of $\SL_4(3)$, we see that $|C_Q(u)| =
|C_Q(ut)|=3^2$. Finally, we have $C_Q(\langle u,t \rangle) = 1$ and so by coprime action $|Q| =
3^6$, a contradiction.  We have proved that $t^G \cap C_G(t) \subseteq L_tO_{p'}(C_G(t))$. In
particular, by Lemma~\ref{wcl} all involutions in $C_G(L_t)L_t$ are conjugate to $t$.

Let $u \in L_t$ be an involution such that $u Z(L_t)/Z(L_t)\in Z(RZ(L_t)/Z(L_t))$.  We have shown that $u \in
t^G$. Set $R_1 = C_R(u)$ and note that as $\langle t,u\rangle$ is a fours group and $u \not \in
Z(R)$, we have $|R:R_1|=2$. Since $u\sim_G t$, there exists $R_0 \in \syl_2(C_G(u))$ such that
$R_0> R_1$. Obviously $R_1$ is normal in $R_0$. If $R_0 \leq C_G(t)$, then $\Omega_1(Z(R_0))=
\langle t, u\rangle$ and this contradicts the fact that $|\Omega_1(Z(R))|=2$. Hence $R_0$ doesn't
centralize $t$. Let $W \cong \Z_4 \ast \Q_8 \ast \Q_8$ be the preimage in $L_t$ of an elementary
abelian subgroup of $L_t/Z(L_t)$ of order $16$ which is contained in $R$. Then $|W:C_W(u)|=2$ and
$C_W(u) = \langle u \rangle \times V$, where $V \cong \Z_4 \ast \Q_8$.  Set
$R_2 = \langle t^G \cap R_1 \rangle$. Then $R_2 \leq R \cap L_t$ and
$R_0$ normalizes $R_2$. As $C_W(u)$ is generated by involutions, we have that $C_W(u) \leq R_2$. As
$C_{R/Z(L_t)}(C_W(u)/Z(L_t)) = W/Z(L_t)$, we see that $Z(R_2) = Z(C_W(u)) = \langle u
,Z(L_t)\rangle$. In particular, $R_0$ normalizes $\langle u,Z(L_t)\rangle$ and hence $R_0$
normalizes $\langle t\rangle$ which is the characteristic subgroup of $\langle u,Z(L_t)\rangle$
generated by squares. Since $R_0 \not \le C_G(t)$ we have our final contradiction. Hence there
exists $t \in \mathcal I$ such that $L_t/Z(L_t) \cong \PSL_2(p^f)$ for some $f \ge 2$.\end{proof}

\begin{proof}[Proof of Theorem~\ref{L2pthm}] This is a consequence of a combination of all the lemmas in this section culminating in
Lemma~\ref{l2q}.
\end{proof}

\section{Centralizers of involutions with $F^*(C_G(t)/O_{p'}(C_G(t))) \cong \PSL_2(p^f)$}\label{sec6}
We continue to assume that Hypotheses~\ref{hypH} and \ref{EH=1} hold. We use the notation
established in the Sections~\ref{sec4}, \ref{sec5} and \ref{sec7}. In particular, $Q= O_p(H)$.
Because of Theorem~\ref{L2pthm}, we may assume that the following hypothesis is satisfied:

\begin{hypothesis}\label{L2p} Hypotheses~\ref{hypH} and \ref{EH=1}
hold and either
\begin{enumerate}
\item there exists $t \in \mathcal I$ and  $F^*(C_G(t)/O_{p'}(C_G(t)))\cong \PSL_2(p^f)$
for some $f > 1$; \item for all involutions $t \in H$,  $C_Q(t)=1$ and
$F^*(C_G(t)/O_{p'}(C_G(t)))\cong \PSL_2(p^f)$ with $f > 1$ and $p \equiv 3 \pmod 4$.
\end{enumerate}
\end{hypothesis}

Our objective is to prove that if Hypothesis~\ref{L2p} holds, then $F^*(G) \cong
{}^2\G_2(3^{2a+1})$ for some $a \ge 1$. We fix the following notation throughout this section. We
let $t \in H$ be an involution such that either $t \in \mathcal I$ or $C_Q(t)=1$ and such that
$E(C_G(t)/O_{p'}(C_G(t)))\cong \PSL_2(p^f)$ with $f > 1$. By Lemmas~\ref{thereisacomponent} and
\ref{component} there is a normal $2$-component of $C_G(t)$, which we denote by $L$,  with
$LO_{p'}(C_G(t))/O_{p'}(C_G(t)) =F^*(C_G(t)/O_{p'}(C_G(t)))$.  Since $t$ is not a classical
involution by Hypothesis~\ref{EH=1} (ii), we have $L/O(L) \cong \PSL_2(p^f)$.

We fix $T \in \syl_2(O_{p'}(C_G(t)))$, $S \in \syl_2(C_G(t))$ such that $T \le S$, $D = S \cap L$
and $U \in \syl_2(G)$ with $S \le U$.

Obviously $t \in T$ and, by Lemma~\ref{cyclic}, we have that $T$ is cyclic and $O_{p'}(C_G(t)) =
TO(C_G(t))$. Finally we note that $$F^*(C_G(t)/O(C_G(t))) \cong T \times L/O(L) \cong \Z_{2^k}
\times \PSL_2(p^f)$$ for some $k \ge 1$ and $f > 1$.

\begin{lemma}\label{nofieldauto} There is no involution in $C_G(t)$  which induces a non-trivial  field automorphism on $L/O(L)$.
\end{lemma}

\begin{proof} Assume that there is an involution $x\in S$ which induces a non-trivial field automorphism
on $L/O(L)$. Then $C_{L/O(L)}(x)\cong  \PGL_2(p^r)$ where $2r = f$. In particular, $p$ divides
$|C_{L}(x)|$ and so, as $H$ is strongly $p$-embedded, there is a $C_G(t)$-conjugate $s$ of $x$
contained in $C_H(t)$. Since $C_G(\langle s,t\rangle)$ is not a $p'$-group, Lemma~\ref{CQt=1}
implies that $C_Q(t) \neq 1 \neq C_Q(s)$. In particular, Lemma~\ref{thereisacomponent} implies that
$O(L) = Z(L)=1$ and that  $C_G(s)$ has a normal component $L_s$ with $L_s \not\le O_{p'}(C_G(s))$.
 Since
$C_{C_{Q}(t)}(s) = C_{C_{Q}(t)}(st)$ and $C_G(Q) \le Q$, we may assume that $C_Q(s) \not\le
C_Q(t)$. Therefore, $|C_Q(s)| > p^r$ and so $L_s \not \cong \PSL_2(p^r)$. Observe that $t$ induces
a non-trivial automorphism on $L_s$ which centralizes a subgroup of $L_s$ which is isomorphic to
$\PSL_2(p^r)$. Suppose that $p^r \neq 3$. Then Lemma~\ref{E3} shows that either
$(L_s/Z(L_s),p)=(\PSL_2(p^f),p)$, $({}^2\G_2(3^{2a-1}),3)$, $(\PSL_2(8),3)$, $(\PSL_3(4),3)$ or $({}^2\B_2(32),5)$.
As $C_{L_s}(t)$, contains $\PSL_2(p^r)$, we see that $(L_s/Z(L_s),p) = (\PSL_2(p^f),p)$ or $({}^2\G_2(3^r),3)$.
But in the latter we have that $C_{C_G(s)/O_{p'}(C_G(s))}(t) = \langle t \rangle \times \PSL_2(3^r)$, while $C_G(t) \cap C_G(s)$ involves $\PGL_2(p^r)$,
a contradiction.   Thus, if $p^r > 3$, we deduce that $L_s \cong \PSL_2(p^f)$ and $t$ induces
a field automorphism on $L_s$. Assume now that $p^r = 3$. Then we have that $L_t \cong \PSL_2(9)$
and $C_{L_t}(s)\cong \Sym(4)$. Furthermore, we have $C_{C_G(t)}(s)$ is soluble. Since $L_s \not \le
L_t$, $t$ does not centralize $L_s$. Thus $t$ induces an automorphism of $L_s$ and, as
$[C_Q(s),O_{p'}(C_G(s))]=1$, $C_Q(s) \cap C_{L_t}(s)
>1$ and $C_{L_t}(s)\cong \Sym(4)$,
$$C_{L_t}(s) \cap O_{p'}(C_G(s))\le C_{C_{L_t}(s)}(C_Q(s) \cap C_{L_t}(s))\le Z(C_{L_t}(s))=1.$$
Thus $t$ centralizes a subgroup of $C_G(s)/O_{p'}(C_G(s))$ which is isomorphic to $\Sym(4)$.

Using Lemma~\ref{E3}, we have the following possibilities for the isomorphism
type of $L_s/Z(L_s)$: $L_s/Z(L_s) \cong \PSL_2(3^m)$, $m \ge 2$, ${}^2\G_2(3^{2a+1})$, $a\ge 1$,
$\PSL_3(4)$ or $\PSL_2(8)$. In the first case, we get $|C_Q(s)| =3^m$ and so, as $t \in \mathcal
I$, we must have $L_s \cong \PSL_2(9)$ as well, with $t$ inducing a non-trivial field automorphism
on $L_s$. If $L_s \cong {}^2\G_2(3^{2a+1})$, then  $t$ acts as an inner automorphism of $L_s$ and
so $t$ centralizes a subgroup  of $L_s$ isomorphic to  $\PSL_2(3^{2a+1})$. This subgroup would of
course have to be contained in $\PSL_2(9)$ which is impossible. Thus $L_s$ is not a Ree group.
Suppose that $L_s/Z(L_s) \cong \PSL_3(4)$. We now use the fact that $C_{C_G(s)/O_{p'}(C_G(s))}(t)$
contains a subgroup isomorphic to $\Sym(4)$ and is soluble. If $t$ acts as an inner automorphism of
$L_s$, we have $C_{L_s}(t)$ is a $2$-group which is not the case. Thus $t$ acts as an outer
automorphism of $L_s$. Hence Lemma~\ref{L34} indicates that $C_{C_G(s)/O_{p'}(C_G(s))}(t)$ is not
soluble which is also impossible. Finally, if $L_s \cong \PSL_2(8)$, we see
$C_{C_G(s)/O_{p'}(C_G(s))}(t) \cong 2\times \Alt(4)$ which does not contain a subgroup isomorphic
to $\Sym(4)$. Hence this case is also impossible. Thus we have shown that, if $L_t \cong
\PSL_2(9)$, then  $L_s \cong \PSL_2(9)$ and $t$ induces a field automorphism on $L_s$.

We have
shown that, if $L_t \cong \PSL_2(p^f)$ and $s$ induces a non-trivial field automorphism on $L_t$,
then $L_s \cong \PSL_2(p^f)$ and $t$ induces a non-trivial field automorphism on $L_s$. It follows
that $L_s \cap L_t \cong \PGL_2(p^r)$ and that $L_s\cap L_t \cap H = N_{L_s\cap L_t}(C_Q(\langle
s,t\rangle))$. Therefore there is an involution $u \in L_s \cap L_t \cap H$, which centralizes
$\langle s,t \rangle$ and inverts $C_Q(x)$ for  $x \in \{s,t\}$. If $C_Q(st) \le C_Q(t)$, then $u$
inverts $C_Q(st)$ while, if $C_Q(st) \not \le C_Q(t)$, we have that  $L_{st} \cong \PSL_2(p^f)$ by
the argument above. Since $L_s \cap L_t= C_{L_t}(s)= C_{L_t}(st)$, we deduce that $L_s \cap L_t =
L_s \cap L_{st}$.  In particular,  $u \in L_{st}$ inverts $C_Q(st)$ in this case as well.
Therefore, in any case we have that $Q = C_Q(t)C_Q(s)C_Q(st)$ is inverted by $u$. Hence $C_Q(u)=1$.
Let $L_u$ be the $2$-component in $C_G(u)$ given by Lemma~\ref{component}.

If $(L_u \cap H)/O(L_u)$ is $p$-closed, then $\langle s,t\rangle $ normalizes a Sylow $p$-subgroup
of $L_u \cap H$. Thus $p$ divides $|C_G(\langle u,x\rangle)|$ for some $x \in \langle
s,t\rangle^\#$. This of course contradicts Lemma~\ref{CQt=1}. Thus $(L_u\cap H)/O(L_u)$ is not
$p$-closed and hence, either  $p
> 3$ and $L_u/O(L_u) \cong \Alt(2p)$ or $2\udot \Alt(2p)$, or $p =5$ and $L_u \cong \Fi_{22}$ or
$2\udot \Fi_{22}$. Notice that $C_{L_t}(u)$ is a Dihedral group of order $2(p^f\pm 1)$. Thus
$C_{C_G(t)}(u)= C_{C_G(u)}(t)$ is soluble.  Hence $t$ induces a non-trivial automorphism on
$L_u/O(L_u)$. If $L_u/O_{p'}(L_u) \cong \Alt(2p)$, then by considering the possibilities for $t$ we
see that $p=5$ and $L_u/O_{p'}(L_u) \cong \Alt(10)$. Furthermore, in this case we have
that the image of $t$ corresponds to a permutation of cycle shape $1^2.2^4$ or $1^4.2^3$. In each
case $C_{L_u}(t)$ contains a section which is isomorphic to $\Sym(4)$. Since $T$ is cyclic, we see
that $C_{C_G(t)}(u)$  contains no such section. So assume that $L_u/O_{p'}(L_u) \cong \Fi_{22}$. In
this case we select an involution $w$ of $L_u$ which centralizes a non-trivial $5$-subgroup of
$L_u$ (see \cite[page 160]{Atlas} to see that this is possible). Then $w$ may be chosen in $H$ and
we have $C_G(\langle u,w \rangle)$ has order divisible by $5$ in contradiction to
Lemma~\ref{CQt=1}. This eliminates all the possibilities for $L_u/O(L_u)$ and so we conclude that
no involution induces a field automorphism on $L_t$ as claimed.
\end{proof}

\begin{lemma}\label{|T|=2} We have $|T|= 2$.
\end{lemma}
\begin{proof} Suppose that $|T|>2$. We have that $D$ is a dihedral group and, setting  $R = T \times D$, we
have $R\in \syl_2(F^*(C_G(t))/O(C_G(t)))$. Let $v\in S$  denote a $\PGL_2$-automorphism of $L$
(which may be trivial). We remark that $T$ is normal in $S$, $R\langle v\rangle/T$ is a dihedral
group and we note that we do not know the action of $v$ on $T$. Finally, let  $z$ be an  involution
in $C_D(\langle v \rangle)$. Then, by Lemma~\ref{nofieldauto}, $\Omega_1(S) \leq (T \times
D)\langle v \rangle$ and $\Omega_1(C_S(\Omega_1(S))) \le (T \times D)\langle v \rangle$.

If $\Omega_2(T) \le C_S(\Omega_1(S))$, we see that
$$\langle t \rangle = \Omega_1(\Phi(
C_S(\Omega_1(S)))).$$ On the other hand, if  $\Omega_2(T)$ does not centralize $\Omega_1(S)$, then
$v$ is an involution and $[\Omega_2(T),v] \not= 1$. In particular $\Omega_2(T) \leq \Omega_1(S)$. Set $$\mathcal W = \{W\mid W \leq \Omega_1(S),
W \cong \Z_4 \times \Z_2 \times \Z_2 \}.$$ Suppose that there is a $W \in \mathcal W$ with $W \not
\le R$. Then $W \cap R$ has order $8$ and is either elementary abelian or is isomorphic to
$\Z_4\times \Z_2$. In the former case, $W\cap D$ is elementary abelian of order $4$ and as
$R\langle v \rangle/T $ is a dihedral group, we have that $C_{\Omega_1(S)}(W \cap D) \le R$, which
is a contradiction as $W \not \le R$. Therefore $W \cap R \cong \Z_4\times \Z_2$ and we may suppose
that $v \in W$. Since $C_{\Omega_1(S)}(v) = \langle t,v,z \rangle$, this case cannot happen either.
Hence every member of $\mathcal W$ is contained in $R$ and consequently $R_1 = \langle \mathcal
W\rangle \mbox{ char } S$. Thus  $\langle t \rangle = \Omega_1(\Phi (Z(R_1)))$ is a characteristic
subgroup of $S$. We have  shown in both cases that
$$\langle t \rangle \mbox{ char } S.$$
Hence   $S$ is a Sylow 2-subgroup of $G$ and  Burnside's Lemma implies
$$t^G \cap Z(S) = \{t\}.$$
Since every involution in $TL/O(L)$ is conjugate to an involution in $Z(R)O(L)/O(L)$, we  have
$$t^G \cap R = \{t\}.$$ On the other hand, by the $Z^\ast$-Theorem, $t^G \cap S \neq \{t\}$ and so   we may assume that $t
\sim_G v$. As $v \not\in \Phi (C_S(t))$, we have that also $t \not\in \Phi (C_S(v))$. This shows
that $[\Omega_2(T), v] \not= 1$, which means that $v \sim_G vt$. Furthermore, conjugating by
elements from $D$, we see that $v \sim_G vz$. Hence also $v \sim_G vtz$. Thus $z$ and $tz$ are the
only involutions in $\Omega_1(Z(C_S(v))) = \langle t,v,z \rangle$ which are not conjugate to $t$ in
$G$. But then $N_G(C_S(v))$ normalizes $\langle z,tz\rangle$ and consequently $N_G(C_S(v)) \le
C_G(t)$. This contradicts $t \sim_G v$. Thus  $T = \langle t \rangle$ has order $2$ as claimed.
\end{proof}

\begin{lemma}\label{reallynofield} No element of $C_G(t)$ induces a non-trivial field automorphism of even order on $L$.
\end{lemma}

\begin{proof} Let $R= TD$. Suppose that some element of $S$ induces a non-trivial field automorphism on $L$.
Then $f$ must be even and so $D \cong \Dih(2^m)$ with $m \ge 3$ and, by Lemmas~\ref{nofieldauto}
and \ref{|T|=2}, $\Omega_1(S) = R\langle v\rangle$ where $v$ is a perhaps trivial $\PGL_2$
automorphism of $L$. Thus we have $\Omega_1(S) \cong  \Z_2 \times \Dih(2^{m+1})$ if $v \neq 1$ and
otherwise  $\Omega_1(S) \cong \Z_2\times \Dih(2^m)$. Let $z \in Z(D)$. Then $\langle z \rangle =
Z(\Omega_1(S)) \cap \Omega_1(S)^\prime$. Hence
$$\langle z \rangle \mbox{ char } S.$$
Set $$ W = \langle x~|~x^2 = z, x \in \Omega_1(S) \rangle.$$ Then $W$ is characteristic in $S$ and
$$W \cong \Z_4 \times \Z_2.$$ Additionally, we have that $C_{\Omega_1(S)}(W) \cong \Z_{2^m} \times \Z_2$ if $v \neq 1$ and $C_{\Omega_1(S)}(W) \cong \Z_{2^{m-1}} \times \Z_2$ if $v= 1$.
As any field automorphism of $L/O(L)$ contained in $S$ centralizes the cyclic group of order 4 in
$D$, we get
$$C_S(W) = Z\langle y \rangle$$
where $Z \leq D\langle v \rangle$ is  cyclic and has index $2$ and $y$ induces a non-trivial
field automorphism on $L/O(L)\cong \PSL_2(p^f)$. By Lemma~\ref{nofieldauto}, we also have $t \in
\langle y \rangle$. We again aim to show that $t^G \cap R = \{t\}$. Let $y_1 \in \langle y \rangle$
with $y_1^2 = t$ and  assume that $Z$ is generated by $c$. Suppose that $[c,y_1] = 1$. Then we have
that $[D,y_1] = 1$ and so $v=1$ which mean that $\Omega_1(S) = R$. This shows that $C_S(R) =
\langle y_1,z \rangle$ and implies
$$\langle t \rangle = \Omega_1(\Phi (C_S(R))) \mbox{ char } S.$$
Assume now that $[c,y_1] \not= 1$. Then, as $C_{L}(y_1) \cong \PGL_2(p^r)$ where $2r=f$, $c$ has
order at least $8$. Furthermore either $[c,y_1] = z$ or $[c,y_1] = tz$. In both cases $[c^2,y_1] =
1$. For $h \in S$, set
$$X_h = \langle s~|~ s \in C_S(W), s^2 = h \rangle.$$ We have that $(cy_1)^2 = cy_1cy_1 =
tc^2[c,y_1]$. As $c^2$ is of order at least four we have that $cy_1$ is of order at least 8. So for
$h = z$ or $h = tz$ we have that $X_h \subseteq \langle c^2, y_1 \rangle$. But then
$$X_t =\langle y_1,c^{2^{m-1}}\rangle \cong \Z_4 \times \Z_2$$
while
$$X_{tz} =\langle c^{2^{m-2}}y_1\rangle\cong \Z_4.$$
Since $C_S(W)$ is characteristic in $S$, we conclude that $t$ and $tz$  are not conjugate by an
automorphism of $S$. Since $\langle z\rangle$ is characteristic in $S$, we deduce that
$$\langle t \rangle \mbox{ char } S.$$
Therefore $S\in \syl_2(G)$ and, as every involution in $LTO(C_G(t))/O(C_G(t))$ is conjugate into
$Z(R)O(C_G(t))/O(C_G(t))$, Burnside's Lemma \cite[6.2]{GLS2} implies $$t^G \cap R = \{t\}.$$ By the
$Z^\ast$-Theorem, we may assume that $t \sim_G v$. We have $\Omega_1(Z(C_S(v))) = \langle t,v,z
\rangle$. Suppose that there is some field automorphism $x$ such that $[x,v] \in R \setminus D$.
Then as there are exactly three $LT$-classes of involutions in $R \setminus D$, we deduce that $v
\sim_G vt$. Then, using elements of $D$, we have $v \sim_G vz$ and $vt\sim_G vtz$. Hence $v \sim_G
vz \sim_G vt \sim_G tvz$. This shows that $\{z,tz\}$ are the only involutions in $\langle t,v,z
\rangle$, which are not conjugate to $t$ in $G$. Hence $N_G(C_S(v)) \leq C_G(t)$, and this
contradicts $t \sim_G v$. Therefore $S = D\langle v \rangle Y$ where $Y=\langle y\rangle$ is cyclic
$\Omega_1(Y) = \langle t \rangle$ and $Y \cap D \langle v \rangle = 1$. Then by the Thompson
Transfer Lemma~\ref{Thompsontransfer}, we have that $t$ is conjugate to a 2-central involution in
$D \langle v \rangle$ and so $t^G \cap R \not= \{t\}$, a contradiction.
\end{proof}

We now summarize what we have discovered about the structure of $C_G(t)$. By Lemma~\ref{|T|=2}, we
have that $T =\langle t\rangle$ and Lemmas~\ref{nofieldauto} and \ref{reallynofield} imply that no
element of $S$ induces a field automorphism on $L/O(L)$. Thus we have $$F^*(C_G(t))/O(C_G(t)) \cong
\langle t\rangle  \PSL_2(p^f)$$ and  $C_G(t)/O_{p'}(C_G(t))$ is isomorphic to either $\PSL_2(p^f)$,
$\PGL_2(p^f)$ or to $\PSL_2(p^f)\udot 2$  where the last extension is by a product of a field and
diagonal automorphism (which is necessarily non-split). In particular, we have that $S/\langle t
\rangle$ is either a dihedral or semi-dihedral group.

\begin{lemma}\label{S=T=>Abelian}  Assume Hypothesis~\ref{L2p}. Suppose that $\syl_2(C_G(t))
\subseteq \syl_2(G)$. Then $G$ has elementary abelian Sylow $2$-subgroups.
\end{lemma}

\begin{proof} Aiming for a contradiction, suppose that $S\in \Syl_2(C_G(t)) \subseteq \Syl_2(G)$ and that $S$
is not elementary abelian.  Then $Z(S)= \langle t,z\rangle$ where $\langle z\rangle = Z(S) \cap D$.
Since $z$ is a commutator in $S$ and $t$ is not, $z$ and $t$ are not $G$-conjugate and so, as $S
\in \Syl_2(G)$, Burnside's Lemma \cite[6.2]{GLS2} implies that $Z(S)$ contains representatives from
three distinct $G$-conjugacy classes. By Hypothesis~\ref{hypH} (iv), $G= O^2(G)$ and so we must
have that $t \in \Phi(S)$ by Lemma~\ref{ThompsonTransferC}. Since $t \in \Phi(S)$, $S/D$ is cyclic
of order $4$. In particular, $\Omega_1(S) = TD$. Since $T  L$ has exactly three conjugacy classes
of involutions with representatives $z$, $t$ and $zt$, we deduce that $t^G \cap TL = \{t\}$.
Therefore $t^G\cap S \subseteq t^G\cap TD =\{t\}$. Finally the Glauberman $Z^*$-Theorem \cite{Glau}
implies that $t \in Z^*(G)=1$ and we have our contradiction. Hence, if $S\in \syl_2(G)$, then $S$
is abelian and consequently is elementary abelian.
\end{proof}

We recall that $U$ is a Sylow $2$-subgroup of $G$ containing $S$.

\begin{lemma}\label{normalV4} $U$ has a normal elementary abelian subgroup of order~4.
\end{lemma}

\begin{proof} As $U$ is not dihedral or semi-dihedral, this follows from \cite[10.11]{GLS2}.
 \end{proof}

\begin{lemma}\label{D8}  $S$ is either elementary abelian or $S/T \cong \Dih(8)$. Furthermore, if $S$ is not abelian, then there is a fours group in $S$ which
is not contained in $Z(S)$ but is normal in $U$.
\end{lemma}

\begin{proof}  We may suppose that $S$ is non-abelian. Let $\langle z \rangle = Z(S)\cap D$.
Since $S$ is non-abelian, Lemma~\ref{S=T=>Abelian} implies that $U \neq S$. By
Lemma~\ref{normalV4}, there exists a normal elementary abelian subgroup $V$ of $U$ of order 4. As
$G=O^2(G)$, the Thompson Transfer Lemma~\ref{Thompsontransfer} implies that $t$ is conjugate to
some involution $s \in C_G(V)$ such that $U$ contains a Sylow 2-subgroup of $C_G(s)$. Hence we may
assume that $V \leq S$.

Suppose  that $V = Z(S) = \langle t,z \rangle$. Then, as $U \neq S$, $|U : S| =2$ and so we can
write $U = S\langle x \rangle$ for some $x \in U$. Since $z$ is a commutator in $S$ and $t$ is not,
$t$ and $z$ are not $U$-conjugate. Therefore, as $U\neq S$, $t^x = tz$. In particular, $C_Q(t) \neq
1$ and $O(L)=1$.

Let $\langle z,s \rangle$ be a fours group in $D$. Note that, as $S/T$ contains a dihedral subgroup
of order at least $8$, $N_{SL}(\langle z,s\rangle)/T \cong \Sym(4)$ and that $N_{SL}(\langle
z,s\rangle)$ normalizes $E= \langle s,t,z \rangle$. Since $S/\langle t \rangle$ is either dihedral
or semi-dihedral,  $E \in \Syl_2(C_G(E))$. By considering the action of $N_L(E)$ on $E$, we see
that $tz\sim_G ts\sim_G tsz$ and $z\sim_G zs\sim_G z$ and by assumption we have $t^x = tz$. So the
involutions in $E$ are partitioned into two sets $t^G \cap E$ of size $4$ and $z^G\cap E$ of size
$3$. Suppose that $|E \cap L^x|=4$, then $N_{(SL)^x}(E\cap L_x)\cong \Sym(4)$. Then we must have $E
\cap L^x= \langle z,s\rangle$. Since $tz$ is centralized by $L^x$, we infer that $N_G(E)/C_G(E)
\cong \Sym(4)$.

Let $R \in \Syl_2(N_G(E))$. Then $R/E \cong \Dih(8)$. We claim $U$ contains no subgroup $R_0$ which
is isomorphic to $R$. Assume this is false and let $F$ be the subgroup of $R_0$ with $R_0/F \cong
\Dih(8)$. Since $S$ has a cyclic subgroup of index $4$, $R_0$ has a cyclic subgroup $C$ of index
$8$ with $C \le S$.  Since $|R_0|= 2^6$, $|C| \ge 8$. We conclude that  $CF/F$ is
 cyclic of order $4$, $Z(CF)$ has order $2$ and $|C|=8$.

 If $F \le S$, then,
 as the $2$-rank of $S$ is $3$, we have that $\langle t,z\rangle=Z(S) \le F$ and as $C
 \le S$, we have a contradiction to $|Z(CF)|=2$. Thus $F \not \le S$ and $(R_0\cap S)/(F\cap S)
 \cong R_0/F\cong
 \Dih(8)$. Since $C$ is inverted in $R_0\cap S$, we have that
 $[R_0\cap S,C]$ has order $4$. Because, $|F\cap S| = 2^2$,  $[R_0\cap S,C]$ centralizes $F \cap
 S$ and so the structure of $S$ indicates that $F\cap S\le Z(S) = \langle t,z\rangle$. But then
 $F\cap S \le Z(R_0 \cap S)$ which has order $2$ and we have a
 contradiction. This contradiction arose from the assumption that $R$ was isomorphic to a subgroup of $U$ and in turn this followed from the supposition that
  $|E \cap L^x|=4$. Hence we must have $|E \cap L^x|\le 2$.

  It follows that $x$ does not normalize $E$ and $EL^x/O(L^x)
 \cong \PGL_2(p^f)$. Therefore $S = \langle t\rangle  \times D_0$ where $D_0$ is a dihedral group  and $x$ can be chosen so that $S= \langle E, E^x\rangle$
and $D_0 = \langle r,s, s^x\rangle$. In particular, $D_0$ is normalized by $x$ and, as all the
involutions in $D_0$ are contained in $L\cup L^x$ (and consequently $G$-conjugate to $z$),  $t$ is
not $G$-conjugate to an element of $D_0$.

Let $u \in  N_{E^xL}(C_Q(t))$, be an involution conjugate into $D_0$. Then $u \in H$ and, as $u$ is
conjugate to $z$, $C_G(u)$ has Sylow $2$-subgroups isomorphic to $U$. By
Lemmas~\ref{thereisacomponent} and \ref{component} there is a $2$-component $M$ in $C_G(u)$. Since
$|Z(U)| = 2$, we infer that $|Z(M/O(M))| \ge 2$. It follows from Lemma~\ref{SchurM} that
$U/O_{p'}(U) \cong \PSL_2(p^a)$, $\Alt(2p)$, $\PSL_3(4)$ with $p=3$  or $\Fi_{22}$ with $p=5$. In
the first case, we have that $u$ is a classical involutions which is impossible. In the remaining
cases, we have incompatible Sylow $2$-structure, as the sectional $2$-rank of $U$ is $3$ whereas
the sectional $2$-rank of the possible components is at least $4$.

Thus we have shown that there is a fours group different from $\langle t,z \rangle$ which is normal
in $S$. This shows that $S/T$ is dihedral of order at most~8. If $S/T$ is abelian, then $S= DT$ is
elementary abelian. This completes the proof of Lemma~\ref{D8}.
\end{proof}
%
%{\bf Comment: The last lemma was altered to deal with the fact that when $T/D_0$ was cyclic,  it
%could be that there was no complement to $D_0$. }

\begin{lemma}\label{tcentral}  $U=S$ is elementary abelian.
\end{lemma}

\begin{proof} Let $V = \langle r,s \rangle$ be a fours group in $D$. Then there is some element of order three in $L$, which acts non-trivially on
$V$ and, by Lemma~\ref{D8}, $V$ is normal in $S$.   By Lemma~\ref{smallgroup}, we have that
$C_G(V)$ has a Sylow 2-subgroup, which is an extension of $B$ by $\langle t \rangle$ where $B \cong
\Z_{2^n} \times \Z_{2^n}$, a Sylow 2-subgroup of $\PSL_3(4)$, a Sylow 2-subgroup of $\SU_3(4)$ or
is elementary abelian of order $16$. We write $S = \langle S \cap L_t, t, y\rangle$ where, if $S$ is
abelian, $y=1$ and, if $S$ is non-abelian, $y^2 \in T$. Set $U_1 = B\langle t,y\rangle$.

Assume first that $U = U_1\in \syl_2(G)$. If $B= V$, then $U=U_1= S$ and we are done by
Lemma~\ref{S=T=>Abelian}.  So we suppose further $B \not= V$ and seek a contradiction. Then $t$ is
not $G$-conjugate to any involution in $B$, as any such involution centralizes an abelian group of
order $16$. In particular, since $G= O^2(G)$, the Thompson Transfer Lemma~\ref{Thompsontransfer}
implies that $U_1 \neq BT$ and $U_1\neq \langle y \rangle B$. Thus $y$ is an involution and, using
the Thompson Transfer Lemma again, we have that $t$ is $G$-conjugate to some element $y_1$ in $By$
and to $y_2 \in Byt$. For $i=1,2$, we have that $y_i \in N_{C_G(t)}(V)$ and $[y_i,V] \not= 1$. As
$N_G(V)/C_G(V) \cong \Sym(3)$, $y_iB$ inverts some element $\rho B$ of order three. As $y_i$ are conjugate to $t$, we see that $C_B(y_i)$ does not
contain an elementary abelian subgroup of order 8. Hence if $B$ is elementary abelian or a Sylow 2-subgroup of $\PSL_3(4)$, we see that all involutions in $By_i$ are conjugate.
In the other cases we have that $V = \Omega_1(B)$ and so $C_B(\rho) = 1$. Therefore we
may apply Lemma~\ref{oneclass} to get that all the involutions in $\langle y_i,\rho\rangle
B\setminus B$ are conjugate to $y_i$ again. In all cases we have that $|C_B(y_i)|^2= |B|$. In particular, we may suppose
that $y_1 = y$ and that $y_2 = yt$. Since $y$ is conjugate to $t$, $|C_G(y)|_2=2^4$. As $\langle
y,t \rangle$ is abelian, it follows that $|C_B(y)| \le 2^2$. Thus $|B|= |C_B(y)|^2 \le 2^4$ and so
$B$ is either elementary abelian of order $16$ or is isomorphic to $\Z_4\times \Z_4$.  We summarize
the conclusions about fusion of involutions in $U_1$ as follows: all the involutions in $By \cup
Byt$ are $G$-conjugate to $t$ and $t^G \cap B$ is empty. The coset $Bt$ may or may not contain
$G$-conjugates of $z$, where $z \in C_V(y)^\#$. Since $z$ and $t$ are not $G$-conjugate, we have
that $D = V$, $LS \cong \Z_2 \times \PGL_2(p^f)$ and $S \cong \Z_2 \times \Dih(8)$. By considering
$N_U(S) \neq S$, we see that $t$ and $tz$ are $G$-conjugate. Therefore $z^G \cap S = V^\#$.

Suppose that $B$ is  elementary abelian. Then $|C_B(y)|= 2^2$ and $C_B(y)$ contains no conjugates
of $t$. Since $y$ is conjugate to $t$, we see that $C_B(y)$ consists of conjugates of $z$. It
follows that $V$ and $C_B(y)$ are conjugate in $G$. Since $|C_B(U)|=2$, $B$ is the Thompson
subgroup of $U$. Therefore $C_B(y)$ and $V$ are conjugate in $N_G(B)$. Hence $yC_G(B)$, $tC_G(B)$
and, by an argument similar to the one above, $tyC_G(B)$ are conjugate in $N_G(B)/C_G(B)$. In
particular, as $tC_G(B)$ commutes with an element $\rho C_G(B)$ of order $3$ which acts
fixed-point-freely on $B$, $yC_G(B)$ centralizes $\rho_yC_G(B)$ and $tyC_G(B)$ centralizes
$\rho_{ty}C_G(B)$ where $\rho_yC_G(B)$ and $\rho_{ty}C_G(B)$  both have order $3$ and both act
fixed-point-freely on $B$.  The isomorphism between $\GL_4(2)$ and $\Alt(8)$ maps such elements of
order $3$ to $3$-cycles. Up to conjugacy, $T\langle \rho\rangle C_G(B)/C_G(B) = \langle
(1,2,3),(1,2)(4,5),(1,2)(6,7)\rangle$ where $t = (4,5)(6,7)$ and, say, $y=(1,2)(4,5)$. The elements
of order $3$ which commute with $y$ and are inverted by $t$ either move $3$ or $8$. In the former
case we see that $N_G(B)/C_G(B)$ contains a subgroup isomorphic to $\Sym(5)$ and this contradicts
the fact that $|N_G(B)/C_G(B)|_2=4$. Thus $\rho$ and $\rho_y$ commute. Since $\rho$ and $\rho_y$
cannot both commute with $\rho_{ty}$, we have a contradiction. Hence $B$ is not elementary abelian.

%
% We have seen that $V$ is the only one in $C_S(t)$ which has the same property. So we get that $V$
%is conjugate to $C_S(y)$ in $N_G(U)$. As $V$ is normal in $S$, we have that the number of
%conjugates of $V$ under $N_G(U)$ is odd. Hence $V$ has at least 5 conjugates under $N_G(U)$. As
%$N_G(U)/C_G(U)$ is a subgroup of $\Alt(8)$, we see that $N_G(U)/C_G(U) \cong \Alt(5)$. But $t$ is
%centralized by some $\Sym(3)$ in $N_G(U)/C_G(U)$, a contradiction.

So we have shown that $B \cong \Z_4 \times \Z_4$. Suppose that $H$ contains a conjugate $z^g$ of
$z$. Then $U^g \le C_G(z^g)$. By Lemmas~\ref{thereisacomponent} and \ref{component}, $C_G(z^g)$ has
a normal $2$-component $L_r$. Since $U^gL_r$ has Sylow $2$-subgroup $U^g$ and $Z(U^g)=\langle
z^g\rangle$, we infer that $Z(L_r/O(L_r))$ has order divisible by 2. Thus Lemma~\ref{SchurM}
implies that $L_r/O_{p'}(L_r) $ is isomorphic to one of $\PSL_2(p^e)$, $\Alt(2p)$, $\PSL_3(4)$ or
$\Fi_{22}$. Since $|U^g/Z(U^g)| = 2^5$, we can only have $L_r/O_{p'}(L_r) \cong \PSL_2(p^e)$. But
then $z^g$ is a classical involution which is impossible by hypothesis. Thus $H$ does not contain
$G$-conjugates of $z$. It follows that $H \cap L$ has odd order and so Corollary~\ref{OrredH2}
implies that $p^f\cong 3 \pmod 4$. On the other hand, $L\langle y\rangle/O(L) \cong \PGL_2(p^f)$
and so $L\langle y\rangle \cap H$ has order divisible by $2$. Hence we may assume that $H \cap
L\langle t ,y\rangle$ contains $E=\langle t,y\rangle$ and $E \in \Syl_2(C_H(t))$.  Since all the
involutions in $E$ are conjugate to $t$, we have  $E \in \syl_2(H)$. Furthermore, using
Lemma~\ref{control1}, we have $N_H(E)/O(N_H(E)) \cong \Alt(4)$. Since all the involutions in $E$
are conjugate, we also have $C_{Q}(t) \not= 1$. Thus Lemma~\ref{thereisacomponent} implies that $L$
is a component and $[L,O_{p'}(C_G(t))]=[L,O(C_G(t))]=1$.

In $C_G(t)$, $E$ normalizes exactly two $p$-subgroups one of which is $C_Q(t)$. Since $C_G(E)$ acts
on the set of $p$-subgroups of $C_G(t)$ which are normalized by $E$, we deduce that
$|C_G(E):C_H(E)|\le 2$. Noting that $C_{C_G(t)}(E)$ has order divisible by $2^3$, we now have that
$|C_G(E):C_H(E)|=2.$ As $y \not\sim_{C_G(t)} yt$, we deduce that $|N_G(E)/C_G(E)|=3$  and
$O^{2}(N_G(E)) \le H$. In particular, $N_G(E)$ has elementary abelian Sylow $2$-subgroups of order
$8$. Let $E_1\in \syl_2(N_G(E))$ be chosen so that $$E_1 = E\langle z\rangle= \langle t,y,z\rangle
\le S$$ where $z \in Z(S)^\#$. By the Frattini Argument, $N_{N_G(E)}(E_1)C_G(E) = N_G(E)$ and so
there exists an element $\rho \in N_G(E)$ of $3$-power order which normalizes but does not
centralize $E_1$ and additionally $\rho^3 \in C_G(E_1)$. Because $t$, $y$, $ty$, $tz$, $zy$ and
$zty$ are pairwise conjugate and $t$ is not conjugate to $z$, we have that $z^G \cap E_1= \{z\}$.
Thus $\rho$ centralizes $z$.

Since $t$ inverts $B$, the preimage of $C_{B/\langle z\rangle}(t)$ in $B$ is the subgroup $X=\{f
\in B\mid f^2 \in \langle z\rangle\}\cong \Z_4\times \Z_2$. Now $y$ acts non-trivially on
$\Omega_1(B)=V$ and, as $y$ has order $2$, a short calculation shows that $y$ normalizes every
subgroup of $X$ of order $4$. It follows that $$C_{B/\langle z\rangle}(t)= C_{B/\langle
z\rangle}(E_1)=X/\langle z \rangle.$$ Hence $N_U(E_1)= E_1X$ and $N_U(E_1)/E_1$ is the four group
of $\Aut(E_1)$ which centralizes $E_1/\langle z\rangle$. This means that $N_U(E_1)$ is extraspecial
of order $32$ of $+$-type. Since $y$ and $ty$ are not conjugate in $C_G(t)$, we have
$N_G(E_1)/C_G(E_1) \cong \Alt(4)$. In particular, $N_U(E_1) \in \syl_2(N_G(E_1))$.
 Let $F\in\syl_2(N_G(E_1))$ be chosen so that it is normalized by $\rho$. Note that $F$ has index
 $2$ in a Sylow $2$-subgroup $U_0$ of $G$ and that the subgroup of $U_0$ corresponding to $B$
 intersects $F$ in a subgroup isomorphic to $\Z_4 \times \Z_2$.
Also  $C_G(F) = \langle z\rangle C_H(F)$
 and $\rho^3 \in C_H(F)$.  Since $F$ is extraspecial of $+$-type and order $32$, $N_G(F)/FC_G(F)$ is isomorphic to a subgroup
 of $\OO^+_4(2)\cong \Sym(3) \wr \Sym(2)$.  As $N_G(F)/FC_G(F)$ has Sylow $2$-subgroups of order two and, as these subgroups consist of a
non-trivial  element which centralizes $\Z_2\times \Z_4$, we see that $N_G(F)/FC_G(F) \cong
\Sym(3)$ or $(\Z_3 \times \Z_3):2$. In both cases $O_3(N_G(F)/FC_G(F))$ is inverted. Suppose the latter. We
now have that the $9$ involutions of $FO(R)/O(R)$ are conjugate in $R/O(R)$. Since $V \le F$, we
get that all the involutions in $F$ are conjugate to $z$, but $t \in F$ and, as $z$ and $t$ are not
conjugate, we have a contradiction. We have
that $U= F\langle u\rangle$ where $u \in B$ has order $4$. And  $N_G(F) = \langle u,\rho\rangle
FC_G(F)$, $C_F(u) = B \cap F \cong \Z_4\times \Z_2$.
As $t$ inverts $B$ all elements in $Bt$ are involutions. Hence there are involutions in $B\langle t \rangle \setminus F$. Choose
such an involution  $w$. Then $w$ centralizes $\Omega_1(B)$. So $Z(B \cap F)\langle w \rangle) = \Omega_1(B)$.
 In
particular, $B\cap F$ is not centralized by $w$, but $B\cap F$ is normalized by $\langle w\rangle$.
It follows that $W=(B\cap F)\langle w\rangle \cong \Dih(8)\times 2$. Now $W=[F,w]\langle w\rangle$
and so $W$ is normal in $U=F\langle w\rangle$. Since $\Aut(W)$ is a $2$-group, $N_G(W) = C_G(W)U$.
Set $\ov {C_G(z)}=C_G(z)/\langle z\rangle$. Then, as $N_G(W) = C_G(W)U$, we have
 $C_{\ov{F}}(\ov w)\langle \ov{w} \rangle$ is a Sylow $2$-subgroup of $C_{\ov{C_G(z)}}(\ov{w})$. In particular, $\ov w$ is not $\ov{C_G(z)}$-conjugate to a subgroup of $\ov F$.
 Therefore, the
Thompson Transfer Lemma~\ref{Thompsontransfer}, implies that $\ov{C_G(z)}$ has a subgroup $R$ of
index two with elementary abelian Sylow $2$-subgroup $\ov{F}$ of order 16. Recall that $\rho$ acts fixed point freely on $V \leq F$ and
so $[\ov{F},\rho] = \ov{F}$. Set $\widetilde R =
R/O(R)\langle z\rangle$.

 Suppose that $E(\wt R) \neq 1$ and recall that $R$ is a $\mathcal
K$-group by hypothesis. Then, as $\wt R$ has elementary abelian Sylow $2$-subgroups of order $16$
which admit a fixed-point-free element of order $3$, we deduce that  $R/O(R)\langle z\rangle$ has
either one or two components isomorphic to $\PSL_2(l)$ for some prime $l \equiv 3,5\pmod 8$ or has
a  single component which is isomorphic  $\PSL_2(16)$. Since $F \in\Syl_2(R)$ and $N_G(F)/FC_G(F) \cong \Sym(3)$,
we see that just one component $\PSL_2(l)$, $l \equiv 3,5 \pmod 8$ is possible. But then we have some element of order three, whose commutator with $F$ is quaternion of order eigth,
contradicting $[F,\rho] = F$. This contradiction  finally shows that $E(\wt R)=1$. In
particular, $C_G(z) = O(C_G(z))N_G(F)$ and, as $\langle z \rangle\in \syl_2(C_{C_G(z)}(\rho))$,
$\langle z \rangle$  is a Sylow 2-subgroup of $C_G(\rho)$. Therefore $C_G(\rho)$ has a normal
$2$-complement. We claim that $\syl_p(C_G(\rho)) \not=\emptyset$. We first investigate  $\rho^3$
which we assume for now is non-trivial. Recall that $\rho \in H$ and $\rho^3 \in C_G(E_1)\le
C_G(t)$. Thus $\rho^3$ normalizes $C_Q(t)$. Since $C_L(z) \cap H =1$, and $\rho^3 \le C_G(z)$, we
have $\rho^3 \not \in L$. We have $C_{C_G(t)}(L) = O_{p'}(C_G(t))$. If $L\langle \rho^3\rangle$ is
isomorphic to a direct product $A \times L$ (necessarily with $A \le H$), then $A \langle \rho
^3\rangle \cap L \le H \cap C_L(z) =1$, which is impossible unless $\rho^3 \in A$. So in this case
we have $C_Q(t) \le C_Q(\rho^3)$. So suppose that $L\langle \rho^3\rangle$ is not isomorphic to  a
direct product. Hence $\rho^3$ induces an outer automorphism of $L$. But then, as $\rho$ is a
$3$-element, $\rho^3$ induces a field automorphism on $L$ and $C_{C_Q(t)}(\rho^3) \neq 1$. Hence in
each case we have $C_Q(\rho^3)\neq 1$. Since $E \in \syl_2(C_H(t))$, $E$ commutes with $\rho^3$ and
$E$ inverts $C_Q(t)$, we have that, setting $J=C_Q(\rho^3)$, $J= C_J(t)\times C_J(ty)\times
C_J(y)$. Hence $\rho$ centralizes a diagonal element in  this decomposition of $J$ and so
$C_J(\rho) \neq 1$. This then proves our claim. Since $\langle z\rangle $ is a Sylow $2$-subgroup
of $C_G(\rho)$ and $C_G(\rho)$ has a normal $2$-complement, we now see that $\langle z\rangle$ is
contained in a conjugate of $H$. But then $\langle z\rangle $ is conjugate to an element of $E$ and
this means that $z$ and $t$ are $G$-conjugate. We have proved that this is not the case and so at
this stage we deduce $U \not= U_1$.

Assume $B \not= V$ and $U \not= U_1$. If $B$ is not elementary abelian of order 16, then we claim
that $V$ is a characteristic subgroup of $U_1$. Set $U_0 = BT$.  If $U_1=U_0$, then $V= Z(U_1)$ and
we are done. Thus $y \neq 1$ and $|U_1:U_0|=2$. Suppose that $\alpha \in \Aut(U_1)$ with $V \neq
\alpha(V)$. Then $U_0 \neq \alpha(U_0)$. Hence $|U_0:U_0\cap \alpha(U_0)|=2$. Assume that $B$ is
homocyclic. Since $B\cap \alpha(B)$ is centralized by $B\alpha(B)$ and $t$ inverts $B$, we either
have that $tB \not \subset B\alpha(B)$ or $B \cap \alpha(B)$ is elementary abelian. In the latter
case we have that $\alpha(V)= \Omega_1(\alpha(B)) =\Omega_1(B\cap \alpha(B)) = \Omega_1(B)= V$,
which  we supposed was not the case. Therefore, $B\alpha(B) \not = U_1$ and so $|B:B\cap
\alpha(B)|=2$ and again we have $\alpha(V)= V$.  So we may assume that $B$ is non-abelian. In
particular, we have that $|B|= 2^6$. Then $B\cap \alpha(U_0)$ has order $2^5$. Now $U_1/\alpha(B)$
is abelian and so $V=B^\prime \le \alpha(U_0)$ and $\alpha(V) = \alpha(B)' \le B$. Thus $Z(B\cap
\alpha(U_0)) \ge V\alpha(V)>V$, this contradicts the structure of $B$ as the centre of every
subgroup of $B$ of index $2$ is $V$. Thus $V$ is a characteristic subgroup of $U_1$ and our claim
is proved.

Set $U_2 = N_{U}(U_1)$. As $B T \in \syl_2(C_G(V))$, we deduce that $|U_2:U|=2$ and that
$U_1 = BT$. In particular, we have $y=1$ and $S$ is elementary abelian of order $2^3$. Since
$C_{U_2}(t) = S \le U_1$, there must be at least two $BT$-conjugacy classes of involutions in $Bt$.
It follows from Corollary~\ref{2groupCor} that $B$ is homocyclic and there are  four $BT$-conjugacy
classes of involutions in $Bt$. Recall that there is an element $\rho \in N_L(V)$ of order $3$ and
that $\rho$ normalizes $B$, centralizes $t$ and $C_B(\rho)=1$. Thus $Bt$ contains two $B\langle
t,\rho \rangle$ classes of involutions. Now $U_2$ does not centralize $t$ and so we deduce $\langle
U_2,\rho\rangle$ induces a transitive action on the four $BT$-conjugacy classes of involutions in
$Bt$. It follows that $|N_G(U_1):U_1|$ is divisible by $4$. But this contradicts $BT \in
\syl_2(C_G(V))$ and $V$ being normalized by $N_G(U_1)$.

So we have shown that $B$ is elementary
abelian of order $16$.  Note that $B$ is characteristic in $U_1$. By Corollary~\ref{2groupCor}, all
the involutions in $Bt$ are conjugate in $BT$. Therefore, as $U_1 < U$, we deduce that $U_1= BS >
BT$ and so, in particular, $|U_1| = 2^6$. By Lemma~\ref{D8}, there is a non-central fours group $X$
of  $S$ which is normal in $U$. Since $U_1 \in \Syl_2(N_G(V))$ and $U>U_1$, $X \neq V$. We have
that $X\le S \le BS$ and so $X$ normalizes $B$ and $B$ normalizes $X$. Therefore $[B,X] \le B\cap X
\le B\cap S = V$. Since $X \neq V$, we infer that $[B,X]=B\cap X = \langle z\rangle = Z(S)$.
However, $|[B,x]|= 2^2$ for all $x \in BS \setminus B$ and so we have a contradiction.

%
%
%
%
%
%As $B$ is characteristic in $U_1$, we get that $U_2=N_U(U_1)$ normalizes $B$. Now $N_G(B)/C_G(B)$
%is isomorphic to a subgroup of $\GL_4(2) \cong \Alt(8)$. Further $C_{N_G(B)}(t)B/B \cong \Z_2
%\times \Sym(3)$ and $|U_2C_G(B)/C_G(B)| \ge 8$. Thus, the subgroup structure of $\GL_4(2)$
%indicates that $N_G(B)/C_G(B) \cong \Sym(3) \wr \Sym(2)$ or $\Sym(5)$. In the first case let $U_3$
%be a $2$-group which  contains $B$ and projects to  a cyclic group of order 4 and in the second
%case be project to a Sylow 2-subgroup of $\Alt(5)$. In both cases we have for some involution $x
%\in U_3 \setminus B$, that all involutions in $Bx$ (if any) are conjugate to $x$ in $B\langle
%x\rangle$. Hence we get that $|C_{U_2}(x)| = 2^5$. In particular, as $t^G \cap B= \emptyset$, $t$
%is not conjugate to any element in $U_3$, which contradicts the Thompson Transfer Lemma
%\ref{Thompsontransfer} as $G = O^2(G)$.

So we are left with the case that $B = V$ and $U > U_1 = S= \langle V,t,y\rangle$. Let $U_2=
N_U(S)> S$. We claim that $VT$ is normal in $U_2$. This is obviously true if $y=1$. So suppose that
$S$ is non-abelian. Then, by Lemma~\ref{D8}, there is a fours group $X \le S$ which is normal in
$U$ and is not contained in $Z(S)$. It follows that $TX$ is an elementary abelian group of order
$8$. Since $S$ contains exactly two such subgroups, we deduce that both $TX$ and $TV$ are normal in
$U_2$. Let $E = VT$. Since $U_2$ does not centralize $t$ and since $E$ is normalized by an element
of order $3$ in $L$, we see that $|t^G\cap E|\ge 4$. Since $Z(U_2)\cap E \neq 1$, we get $|t^G\cap
E|=4$. Therefore, $|N_G(E)/C_G(E)| = 12$ or $24$. As $N_G(E)/C_G(E)$ is a subgroup of $\PSL_3(2)$,
this means $N_G(E)/C_G(E) \cong \Alt(4)$ or $\Sym(4)$. In both cases $C_E(O_2(N_G(E)/C_G(E)))$ is
non-trivial and normal in $N_G(E)$ and so it must be $V$. But $E=BT \in\syl_2(C_G(V))$ and we have
a contradiction. This is our final contradiction and so Lemma~\ref{tcentral} is proved.
\end{proof}

\begin{lemma}\label{acomp}   $L$ is a component of $C_G(t)$.
 \end{lemma}

 \begin{proof} From Lemma~\ref{tcentral}, we have $G$ has elementary abelian Sylow $2$-subgroups of order $8$ and
 $O^{2'}(C_G(t)/O(C_G(t))) \cong 2 \times \PSL_2(p^f)$ where $f\ge 3$. Assume that the lemma is false.
 By Lemma~\ref{thereisacomponent}, we have
 $C_Q(t)=1$.  Let $P \in \Syl_p(H \cap L)$.
 If $[P,O_{p'}(C_G(t))]= 1$, then
$P$ centralizes $F(O_{p'}(C_G(t)))$ and so $F^*(C_G(t) )\not\le O(C_G(t))$.  Therefore, $C_G(t)$
has a component and of course it must be $L$, which is a contradiction.  Thus, as $T \in
\Syl_2(O_{p'}(C_G(t)))$, there exists an odd prime $r$ such that $r \neq p$, $r$ divides
$|F(O(C_G(t)))|$ and $P$ acts faithfully on $R= O_r(O(C_G(t)))$. Let $C$ be a critical subgroup of
$R$. Then, as $r$ is odd, we may assume that $C$ is of exponent $r$ (see \cite[Theorem
5.3.13]{Gorenstein}). Choose $D \le C$ minimal  so that $C_{C_G(t)}(D) \le O_{p'}(C_G(t))$ and $D$
is normal in $C_G(t)$. Then, as $\Phi(D) < D$, $C_{C_G(t)}(\Phi(D))$ covers
 $O^{2'}(C_G(t)/O_{p'}(C_G(t)) )\cong \PSL_2(p^f)$.  In particular, as $P$ acts faithfully on $D$, $P$ acts faithfully on $D/\Phi(D)$.
 Since $P$ is elementary abelian of order $p^f$ and
 $N_{L}(P)$ has
orbits of length  $(p^f-1)/(2(p-1))$ on the maximal subgroups of $P$, we have $|D/\Phi(D)| \ge
r^{(p^f-1)/(p-1)2}$. Because $f\ge 3$,  $(p^f-1)/(p-1)2 \ge 7$ and so $|D/\Phi(D)| \ge r^7$. In
particular, as $D$ is normal in $C_G(t)$ and $F^*(C_G(t)/O(C_G(t)))$ is generated by three
involutions, we have that $$|D/\Phi(D):C_{D/\Phi(D)}(x)|> r^2$$ for all involutions $x\in C_G(t)$.

On the other hand, as the Sylow $2$-subgroups of $C_G(t)/\langle t\rangle$ have order $4$ with all
involutions conjugate, we also have $|C_{D/\Phi(D)}(x)|
> r^2$ for all involutions $x\in C_G(t)$.
Therefore, for involutions  $x\in C_G(t)$,  $C_{D}(x)$ contains an elementary abelian subgroup of
order $r^3$.  Let $W \le D$ be an elementary abelian subgroup of order $r^3$ centralizing some
involution $s \in L\setminus H$. Then, as  $W \le O_{p'}(C_G(t)) \le H$ and $C_G(Q) \le Q$, $W$
acts faithfully on $Q$. Since $t$ inverts $Q$, $Q$ is abelian  and so there is a non-trivial $w\in
W$ such that $|C_{\Omega_1(Q)}(w)| \ge p^2$. Then, as $s \in C_G(w)$, $C_G(w)\not \le H$ and
$m_p(C_H(w)) \ge 2$. Hence Lemma~\ref{main1} (v) implies that $O_{p'}(C_G(w)) \le H$ and
Lemma~\ref{main2} (ii) implies that $F^*(C_G(w)/O_{p'}(C_G(w)))$ is a non-abelian simple group.
Since $C_Q(w)\not=1$ and $|[C_Q(w),O_{p'}(C_G(w))]| = 1$, we have that $F^*(C_G(w)) \not\le
O_{p'}(C_G(w))$ and so $C_G(w)$ has a normal component $L_w$ with $L_w \not\le O_{p'}(C_G(w))$.
Since $L_w$ has abelian Sylow $2$-subgroups of order at most $8$, we may apply \cite{Walter} to
deduce that $L_w$ is a $\mathcal K$-group. Hence  $L_w\cong \PSL_2(p^b)$, $b \ge 2$,
${}^2\G_2(3^{2c+1})$, $c \ge 1$ or $\PSL_2(8)$. Now, as $ D \le H$,
$C_D(w)O_{p'}(C_G(w))/O_{p'}(C_G(w))$ normalizes $C_Q(w)O_{p'}(C_G(w))/O_{p'}(C_G(w))$ and is an
$r$-group of exponent $r$. It follows from Proposition~\ref{Hstru1} that $|C_D(w):C_D(w) \cap
O_{p'}(C_G(w))| \le r^2$. Recall that $t$ centralizes $w$ and so there is an involution $x\in L_w$
which centralizes $t$. We have $O_{p'}(C_G(w))\le C_G(L_w) \le  C_G(x)$ and so $|C_D(w):C_D(w)\cap
C_G(w)| \le r^2$. Since $|D/\Phi(D):C_{D/\Phi(D)}(w)|> r^2$, we infer that $w \not \in Z(W) \ge
\Phi(D)$. Because $\Phi(D) \cap O_{p'}(C_G(t))$ centralizes $C_Q(w)$, we have $\Phi(D) \cap
O_{p'}(C_G(t)) =1$  for otherwise $w$ could have been chosen in $\Phi(D)$ which we know is not
possible.  On the other hand, as $\langle w\rangle\Phi(D)$ is elementary abelian and
$[D,\Phi(D)]=1$, we know that $|D:C_D(w)|\le |\Phi(D)|$. As $\Phi(D) \le C_D(w)$ and $\Phi(D) \cap
O_{p'}(C_G(w)) =1$, $|C_D(w):\Phi(D)(C_{D}(w)\cap O_{p'}(C_G(w)))| \le r^2/|\Phi(D)|$. Finally we
have $|D/C_D(w)|\le |\Phi(D)|$  and $|C_D(w):C_{C_{D}(w)}(x)| \le r^2/|\Phi(D)|$ which once again
gives $|D/\Phi(D):C_{D/\Phi(D)}(x)| \le r^2$  our final contradiction. Thus $L$ is a component of
$C_G(w)$.
 \end{proof}

\begin{theorem}\label{Ree} Assume Hypothesis~\ref{L2p}. Then $F^\ast(G) \cong$ $^2\G_2(3^{2a+1})$.
\end{theorem}

\begin{proof} By Lemma~\ref{tcentral} the Sylow $2$-subgroups of $G$ are elementary abelian of order
$8$ and, by Lemma~\ref{acomp}, $L\cong \PSL_2(p^f)$ is normal component of $C_G(t)$. Therefore
Theorem~\ref{E8Sylow} gives $F^*(G)\cong {}^2\G_2(3^{2a+1})$ for some $a \ge 1$.
\end{proof}

%
%
%\begin{lemma}\label{HcapL} If $R \cap H \cap L =1$ then $F^\ast(G) \cong {}^2\G_2(3^{2n+1})$ for some $n \geq 1$.
%\end{lemma}
%
%\begin{proof} As $R \cap H \cap L = 1$ in particular $f$ is odd. Further there is no element in $C_G(t)$, which induces a field automorphism of even order.
% Assume $|T| > 2$. As there is no involutory field automorphism  we have that any group $U \leq R$ isomorphic to $\Z_4 \times \Z_2 \times \Z_2$ has to be contained in
%$C_G(\Omega_2(T))$.  Set $R_1 = \langle U~|~ U \leq R, U \cong \Z_4 \times \Z_2 \times \Z_2
%\rangle$. Then $\Phi (Z(R_1)) = \langle t \rangle$. This first of all shows that $R$ is a Sylow
%2-subgroup of $G$. Further it shows that $t^G \cap C_R(\Omega_2(T)) = \{t\}$. By Lemma~\ref{wcl}
%there is some $x \sim t$, $x \in R \setminus \langle t \rangle$. Hence $x$ induces an outer
%automorphism on $L$ and inverts $\Omega_2(T)$. In particular $x \sim xt$. Now we see that $C_R(x) =
%\langle t,x,u \rangle$, where $u \in R \cap L$. We have that $x \sim xu$ and so we get
%$G$-conjugates $t,x,xt,xu,xut$. Hence $u$ and $tu$ are the only involutions in $C_R(x)$ which are
%not conjugate to $t$ in $G$. In particular $N_G(C_R(x))$ acts on $\langle u,ut \rangle$ and so it
%centralizes $t$. This would show that $C_R(x)$ is a Sylow 2-subgroup of $C_G(x)$, which contradicts
%$t \sim x$.
%
%This contradiction now shows that $T = \langle t \rangle$. So Hypothesis~\ref{L2p} is satisfied.
%Then Proposition~\ref{Ree} gives the assertion.
%\end{proof}
%
%
%
%
%
%
%
%
%
%
%
%
%
%
%

\section{Components in $H$}\label{sec8}

From Lemma~\ref{Eqs}, we know that if $E(H) \neq   1$, then $E(H)$ is quasisimple. The objective of
this section is to prove that $E(H)=1$. Thus throughout this section we assume the following
hypothesis.

\begin{hypothesis}\label{EHnot1} Hypothesis~\ref{hypH} holds as well as \begin{enumerate}
\item $E(H)$ is a quasisimple  $\mathcal K$-group; and
\item $G$ does not contain a classical involution.
\end{enumerate}
\end{hypothesis}

Set $E = E(H)$. So $E$ is a quasisimple group. Since $O_{p'}(G)=1$, we have that $F(H)$ is a
$p$-group. In particular, $E$ contains involutions. If $t \in H$ is an involution, we know by
Lemma~\ref{Ctnot-in} that $C_G(t) \not \le H$, $O_{p'}(C_G(t)) \le H$ and $(C_G(t)/O_{p'}(C_G(t)),
p) \in \mathcal E$. Whenever $t$ is an involution from $H$, we use the following bar notation $\ov
{C_G(t)}=C_G(t)/O_{p'}(C_G(t))$.

\begin{lemma}\label{Op=1}  $E=F^*(H)$.
\end{lemma}

\begin{proof} Set $Q=O_p(H)$ and assume that $Q \not\le E(H)$. Of course $[E,Q]=1$.  Suppose that $t$ is an involution in $E$. Then,
by Lemma~\ref{Ctnot-in}, $C_G(t) \not\le H$, $O_{p'}(C_G(t)) \le H$ and $(F^*(\ov{C_G(t)}),p) \in
\mathcal E$. Since $Q$ and $O_{p'}(C_H(t))$ commute, we have that $F^*(C_G(t)) \not \le
O_{p'}(C_G(t))$. As $F^*(\ov{C_G(t)})$ is almost simple, $C_G(t)$ has a unique component $L_t$
which has order divisible by $p$ and we have $\ov{L_t} =\ov{F^*(C_G(t))}$. Note that, as $\ov {Q}
\neq 1$, Proposition~\ref{Hstru1} implies that $\ov{C_H(t)}$ is $p$-closed and soluble.

Assume that $p$ divides $|C_E(t)/Z(E)|$. Then $O_{p}(\ov{C_H(t)})$ contains non-trivial normal
subgroups $\ov Q$ and $ O_{p}(\ov {C_E(t)})$ and they commute. This contradicts the structure of
$\ov{C_H(t)}$ given in Proposition~\ref{Hstru1}. Therefore, $C_{E}(t) \le O_{p'}(C_H(t))$ for all
involutions $t \in E$. Suppose that $C_E(t) \not \le O_{p'}(C_G(t))$. Then using
Corollary~\ref{OrredH2} and noting that $\ov{C_H(t)}$ is soluble, we get  $|\ov{C_E(t)}|=2$ and
$(\ov L_t,p) =(\PSL_3(4),3)$. Then $\ov Q$ is a normal $3$-subgroup of $\ov{C_H(t)}$ and so
Proposition~\ref{Hstru1}(iv) implies that $N_{L_t}(Q)/Q$ is a non-abelian $2$-group. Suppose that
$s \in C_H(E)$ is an involution. Then $E O_{3'}(C_G(s))/O_{3'}(C_G(s))$ is a quasisimple normal
subgroup of $C_H(s)/O_{3'}(C_G(s))$ which by Lemma~\ref{Ctnot-in}(iii) is strongly $3$-embedded in
$C_G(s)/O_{3'}(C_G(s))$. However, as $(C_G(s)/O_{3'}(C_G(s)),3) \in \mathcal E$, we have a
contradiction to Proposition~\ref{Hstru1} as, when $p=3$, $\ov{C_H(t)}$ must be soluble. Thus
$C_H(E)$ has odd order and, in particular, $N_{L_t}(Q)/Q$ maps isomorphically into $H/C_H(E)$ and
so we deduce that $\Out(E)$ has non-abelian Sylow $2$-subgroups. Now $E$ is a $\mathcal K$-group
and since $\Out(E)$ has non-abelian Sylow $2$-subgroups, we infer that $E$ is a Lie type group of
Lie rank at least $4$ defined over a field of characteristic $r$. But these groups all possess an
involution $s$ with $C_E(s)$ involving $\PSL_2(r)$. Since $3$ divides $|\PSL_2(r)|$, this
contradicts $C_E(s) \le O_{3'}(C_H(s))$. Thus $(\ov {L_t},p) \not =(\PSL_3(4),3)$ and we conclude
that $C_E(t) \le O_{p'}(C_G(t))$ for all involutions $t \in E$.

%
%\begin{claim}\label{Q1} If $x \in E$ is an involution, then $F^*(\ov{C_G(x)})\not \cong \PSL_3(4)$.
%\end{claim}
%
%Assume false. Set $L = E(C_G(x))$. As the Schur
%multiplier of $\PSL_3(4)$ is a $\{2,3\}$-group, $L\cong \widehat
%{\PSL_3(4)}$ is a central extension of $\PSL_3(4)$ by a centre of
%$2$-power order. We already know
%$C_{E}(t) \le O_{p'}(C_H(t))$ for all involutions $t \in H$. As now $p = 3$, we see that $E$ is not sporadic and $E$ is not alternating
%$\Alt(n)$, $ \geq 4$. If $E$ is of Lie type of rank at least three and in case of odd characteristic of rank at least two, there
%is some involution $t \in E$ whose centralizer involves $\PSL_2(q)$. So we have that $E/Z(E) \cong \PSL_2(q)$, $\PSL_3(2^n)$,
%$\PSU_3(2^n)$.
%Using Corollary~\ref{OrredH2} we see that $O_{3'}(C_G(x))$ is of index at most two in
%$O_{3'}(C_H(x))$. Furthermore, $L \cap H= N_L(Q)$ and this group
%contains an involution $s$ which inverts~$Q$. Now the only group $E$ which possesses an outer automorphism which centralizes in a Sylow 2-subgroup
%a subgroup of index two is $\PSL_2(q)$, $q$ odd, where $s$ then is a field automorphism. As $m_3(E) > 1$, we get that $q = 3^2f$.
%Then $|C_E(x)| = 3^{2f}-1$. But $C_E(s) \cong \PGL_2(3^f)$ and so $C_{C_E(s)'}(x)| = 3^f-1$ or $3^f+1$.
%This would imply that $3^f + 1 = 2$ or $3^f - 1 = 2$, a contradiction.
%So we have that in any case $[E,s] = 1$. Now $C_H(s)$ contains $E$ as a normal subgroup and, as $p=3$,
%$C_H(s)/O_{3'}(C_H(s))$  is soluble. Thus $E \le
%O_{3'}(C_H(s))$ but $E$ is normal in $H$ and so we have $E \le O_{3'}(H) = 1$, a contradiction.

Suppose now that $T \in \syl_2(E)$ and $t\in Z(T)^\#$ is an involution. Then  $$T\le C_E(t) \le
O_{p'}(C_G(t))\le C_G(L_t).$$ Assume that $s \in T$ is an involution. Then $C_E(s) \le
O_{p'}(C_G(s))$ and thus
$$[C_E(s), \langle Q^{C_G(s)}\rangle] \le [O_{p'}(C_G(s)), \langle
Q^{C_G(s)}\rangle]= \langle [O_{p'}(C_G(s)), Q]^{C_G(s)}\rangle =1.$$ As ${L_t} \le C_G(s)$, we see
$\langle Q^{C_G(s)} \rangle \ge \langle Q^{L_t}\rangle = {L_t}Q$. Thus ${L_t}$ centralizes
$C_E(s)$. It follows that ${L_t}$ centralizes $K=\langle C_E(x)\mid x\in T, x^2=1\rangle$.  If $p$
divides $|K|$, then we have the contradiction  $L_t \le H$ as $H$ is strongly $p$-embedded. Hence
$p$ does not divide $|K|$ and so  $K \neq E$. It follows from \cite[17.13]{GLS2} that $N_E(K)$ is a
strongly $2$-embedded subgroup of $E$. Hence $E \cong \SL_2(2^a)$, $\PSU_3(2^a)$ or
${}^2\B_2(2^{2a-1})$ for some $a \ge 2$ and $K \le N_G(T)$ by \cite{Bender}. Since ${L_t}$
centralizes $T$ and ${L_t}$ is normal in $C_G(t)$, we see that $L_t$ is a characteristic subgroup
of $C_G(T)$. Hence $N_G(T)$ normalizes ${L_t}$. In particular, $N_G(T)= N_H(T)L_t$ by the Frattini
Argument. By Lemmas~\ref{main1}(v) and \ref{main2} (ii) and (iii), we have
$F^*(N_G(T)/O_{p'}(N_G(T)))\cong L_t/Z(L_t)$. Thus using Corollary~\ref{OrredH2} and the fact that
$L_t/Z(L_t)\not\cong \PSL_3(4)$, we have $O_{p'}(N_G(T)) = O_{p'}(N_H(T))$. If $p$ divides
$|N_E(T)|$, then $N_E(T)O_{p'}(N_G(T))/O_{p'}(N_G(T))$ commutes with
$QO_{p'}(N_G(T))/O_{p'}(N_G(T))$ and as before we have a contradiction to the structure of $L_t$
via Proposition~\ref{Hstru1}. Hence $N_E(T) \le O_{p'}(N_H(T))=O_{p'}(N_G(T))$. In particular,
$L_t$ centralizes $N_E(T)$.
 Let $J$ be a complement to $T$ in $N_E(T)$ and set $M= N_E(J)$. Then $E = \langle
N_E(T),M\rangle$. Since $L_t$ centralizes $N_E(T)$, ${L_t} \le C_G(J)$. Therefore,
Lemma~\ref{main2} (ii) and (iii) imply that ${L_t}$ is the unique component in $N_G(J)$. In
particular,  $M$ normalizes ${L_t}$ and so $E = \langle N_E(T),M\rangle$ normalizes $L_t$ as well.
But again by Lemma \ref{main2} (ii) and (iii), we have that $F^*({L_t}E/O_{p'}(L_tE)) =
L_tO_{p'}(L_tE)/O_{p'}(L_tE)$ which means that $E \le  O_{p'}({L_t}E)$ and contradicts
$O_{p'}(H)=1$. This contradiction shows that $Q \le E$ and so $F^*(H)= E$ is quasisimple.
\end{proof}

Because of Lemma~\ref{Op=1}, we may now assume that $F^*(H)/Z(F^*(H))$ is a simple group. We will
now proceed to investigate the possibilities for $H$ using the fact that $H$ is a $\mathcal
K$-group. We begin with two lemmas which will be used frequently.

\begin{lemma}\label{QS1} If $t \in H$ is an involution and
 $\ov{C_H(t)} $ is not soluble, then either
 \begin{enumerate}
\item $(F^*(\ov{C_G(t)}),p) = (\Fi_{22},5)$ and  $\ov{C_H(t)}\cong \Aut(\Omega_8^+(2))$; or
\item $(F^*(\ov{C_G(t)}),p) = (\Alt(2p),p)$ and $F^*(\ov{C_H(t)})\cong\Alt(p)\times \Alt(p)$ with $p \ge 5$.
\end{enumerate}
 Moreover, in
case (ii), the components of $\ov{C_H(t)}$ are not normal in $\ov{C_H(t)}$.
\end{lemma}

\begin{proof}
This follows from Lemma~\ref{Ctnot-in}(iii) and  Proposition~\ref{Hstru1}.
\end{proof}

\begin{lemma}\label{ppart} Assume that $t \in H$ is an involution and  $|C_H(t)|_p\ne |C_E(t)|_p$. Then $(\ov{C_G(t)} ,p) = (\PSL_2(8):3,3)$ or $({}^2\B_2(32):5,5)$.
In particular, $C_G(t)$ has extraspecial Sylow $p$-subgroups of order $p^3$.
\end{lemma}

\begin{proof} Suppose that $|C_E(t)|_p \not= |C_H(t)|_p$. Then $\ov{C_H(t)}>
\ov{C_E(t)}$ and, as $C_H(t)/C_E(t)$ is soluble and has order divisible by $p$,
Proposition~\ref{Hstru1} implies that $\ov{C_H(t)}$ is $p$-closed. Let $T_1 \in \syl_p(C_H(t))$ and
let $T_2 = T_1 \cap E$. Then $\ov{T_1} > \ov{T_2}$ and $\ov{T_1}$ is normal in $\ov{C_H(t)}$. It
follows from Corollary~\ref{OrredH} that either $T_2 \le \Phi(T_1)$ or  $(\ov{C_G(t)},p)
=(\PSL_2(8):3,3)$ or $({}^2\B_2(32):5,5)$. Thus we may assume that $2 \le m_p(T_2/T_1) \le
m_p(H/E)$.  Since the $p$-rank of $\Out(E)$ is at most $2$, we have that $T_1/T_2$ has $p$-rank $2$
by Lemma~\ref{prankauto}.  But then the structure of $\Out(E)$ in Lemma~\ref{prankauto} shows that
$C_H(t)$ cannot act irreducibly on $\ov{T_1}/\Phi(\ov{T_1})$ and Corollary~\ref{OrredH} once again
gives $p\in\{3,5\}$ and the structure of $\ov{C_G(t)}$.
\end{proof}

\begin{lemma}\label{QS2} $E/Z(E)$ is not an alternating group of degree $n \ge 5$.
\end{lemma}

\begin{proof}
Suppose that $E/Z(E)\cong \Alt(n)$ for some $n \ge 5$. Assume that
$Z(E)=1$. Let $n \ge 9$. Let $t$ correspond to a
product of two transpositions in $E$. Then $C_H(t) $ contains a
normal subgroup isomorphic to $\Alt(n-4)$ with $n-4\ge 5$. Thus
$\ov{C_H(t)}$ contains such a normal subgroup and this contradicts
Lemma~\ref{QS1}. So $n < 9$. But then $C_E(t)$ is a $\{2,3\}$--group with cyclic Sylow 3-group and so $m_r(C_H(t)) < 2$
for all odd primes $r$. This contradicts
Hypothesis~\ref{hypH}(ii).

This contradiction shows that $Z(E) \ne 1$. Since
$Z(E)$ is a $p$-group, this means  $E \cong 3\udot \Alt(6)$ or
$3\udot \Alt(7)$ with $p=3$. The first possibility fails as the
centralizers of involutions in $\Alt(6)$ have order $8$. Thus $E
\cong 3\udot \Alt(7)$. In this case, $\ov{C_H(t)}$ has a normal
Sylow $3$-subgroup of order $9$ and $C_H(t)$ does not act
irreducibly on $O_3(\ov{C_H(t)})$, this contradicts
Corollary~\ref{OrredH}. Hence we have shown $E/Z(E)$ is not an
alternating group.\end{proof}

Next we show that $E/Z(E)$ cannot be a sporadic simple group.

\begin{lemma}\label{QSsporadic} $E/Z(E)$ is not a sporadic
simple group.
\end{lemma}

\begin{proof} We use Lemmas~\ref{Ctnot-in} and \ref{QS1} for all the sporadic groups. We first observe that the outer automorphism group of a
sporadic simple group has order dividing $2$ \cite[Table 5.3]{GLS3}. Hence $E$ has index at most
$2$ in $H$. In Table~\ref{TabSpor}, for each possibility for $E/Z(E)$, we give the structure of the
centralizer of some involution $t$ in $E$ (recall $O_2(E)=1$). For each case, except $E \cong
3\udot \M_{22}$ with $p=3$, we see that either $p^2$ does not divide $|C_E(t)|$ or that
$C_H(t)/O_{p'}(C_H(t))$ has a simple section which is not isomorphic to $\Alt(p)$ or
$\Omega_8^+(2)$. Assume that $ E \cong 3\udot \M_{22}$ and $p=3$. Then $C_H(t)$, with $t$ as in
Table~\ref{TabSpor}, is soluble and has Sylow $3$-subgroups of order $9$. However, $\ov{C_H(t)}$
does not act irreducibly on $O_3(\ov {C_H(t)})$ and this then contradicts Corollary~\ref{OrredH}.
Thus this possibility cannot occur either. Therefore $E/Z(E)$ is not a sporadic simple group.
\end{proof}

\begin{table}
\begin{tabular}{ccc||ccc}
$E/Z(E)$&Involution& $C_{E/Z(E)}(t)$&$E/Z(E)$&Involution& $C_{E/Z(E)}(t)$\\
\hline
$\M_{11}$&$2A$&$\GL_2(3)$&$\M_{12}$&$2A$&$2\times \Sym(5)$\\
$\M_{22}$&$2A$&$2^4.\Sym(4)$&$\M_{23}$&$2A$&$2^4:\PSL_3(2)$\\
$\M_{24}$&$2A$&$2^{1+6}.\PSL_3(2)$&$\J_1$&$2A$&$2\times \Alt(5)$\\
$\J_2$&$2A$&$2^{1+4}_-.\Alt(5)$&$\J_3$&$2A$&$2^{1+4}_-.\Alt(5)$\\
$\J_4$&$2A$&$2^{1+12}.3\udot\M_{22}.2$&$\Co_3$&$2A$&$2\udot\PSp_6(2)$\\
$\Co_2$&$2A$&$2^{1+8}_+.\PSp_6(2)$&$\Co_1$&$2C$&$2^{11}.\M_{12}$\\
$\HS$&$2A$&$4*2^{1+4}_+.\Sym(5)$&$\McL$&$2A$&$2\udot\Alt(8)$\\
$\Suz$&$2A$&$2^{1+6}_-.\Omega_{6}^-(2)$&$\He$&$2A$&$2^2.\PSL_3(4).2$\\
$\Ly$&$2A$&$2\udot \Alt(11)$&$\Ru$&$2A$&$2^{11}.\Sym(5)$\\
$\ON$&$2A$&$4.\PSL_3(4).2$&$\Fi_{22}$&$2A$&$2\udot\PSU_6(2)$\\
$\Fi_{23}$&$2A$&$2\udot\Fi_{22}$&$\Fi_{24}'$&$2A$&$2\udot\Fi_{22}.2$\\
$\HN$&$2A$&$2\udot\HS.2$&$\Th$&$2A$&$2^{1+8}_+.\Alt(9)$\\
$\B$&$2A$&$2.{}^2\mathrm E_6(2).2$&$\M$&$2A$&$2\udot\B$\\
\hline

\end{tabular}
\caption{Centralizers of certain involutions in sporadic simple groups.}\label{TabSpor}
\end{table}

We now begin our investigation of the case when $E/Z(E)$ is a group of Lie type.

\begin{lemma}\label{QS3} Suppose that $E/Z(E)$ is a simple
group of Lie type. Then $p$ divides $|C_E(t)|$ for all involutions $t \in H$.
\end{lemma}

\begin{proof} Assume that $E/Z(E)$ is a Lie type group defined in characteristic $r$ and that
$C_E(t)$ is  a $p'$-group. Then, as $m_p(C_H(t)) \ge 2$, $|C_E(t)|_p \ne |C_H(t)|_p$ and hence
Lemma~\ref{ppart} implies that $p\in\{3,5\}$ and $C_H(t)$ has extraspecial Sylow $p$-subgroups. It
follows that $H/E$ and hence $\Out(E)$ has non-abelian Sylow $p$-subgroups. In particular, we see
that $E$ must admit diagonal automorphisms of order $p$. This shows that $E/Z(E)$ is $\PSL_n(q)$ or
$\PSU_n(q)$ and that $p$ divides $n$. Since $H/E$ has non-abelian Sylow $p$-subgroups, we further
infer that $p^2$ divides $n$. Thus $n \ge 9$. But then the canonical form of $t$ shows that the
centralizer of $t$ has a section isomorphic to  $\PSL_4(r)$ or $\PSU_4(r)$. Both these groups have
order divisible by $3$ and $5$, a contradiction. Thus $p$ divides $|C_E(t)|$.
\end{proof}

\begin{lemma}\label{QS5} $E/Z(E)$ is not a rank $1$ Lie type group.
\end{lemma}

\begin{proof}
Suppose that $E/Z(E)$ is defined in characteristic $2$. So $E/Z(E) \cong \PSL_2(2^a)$,
$\PSU_3(2^a)$ or ${}^2\B_2(2^{2a-1})$ for some $a \ge 2$. In the first and third case, the
centralizer of an involution in $C_E(t)$ is a $2$-group  and so  $m_p(\Out(E)) \ge 2$. Using
\cite[Theorem 2.5.12]{GLS3} we see that $\Out(E)$ is cyclic, which is a contradiction. Therefore
$E/Z(E) \cong \PSU_3(2^a)$ for some $a\ge 2$. Let $t \in E$ be an involution and $S \in
\Syl_p(C_H(t))$. Then $C_E(t)$ is $2$-closed, contains a Sylow $2$-subgroup $T$ of $E$ and $T$ has
a cyclic complement in $C_E(t)$ of order dividing $q+1$. In particular, $|C_E(t)|_p \not =
|C_H(t)|_p$ and so $S$ is extraspecial of order $p^3$  and $(\ov{C_G(t)},p) =
(\PSL_2(8)\colon{3},3)$ or $({}^2\B_2(32)\colon 5,5)$ by Lemma~\ref{ppart}.
 Since
$\Out(E)$ has abelian Sylow $p$-subgroups, $p$ must divide $q+1$. Furthermore, the Sylow
$p$-subgroups of  $\mathrm{GU}_3(2^a)$ are abelian and so  $SE$ must involve field automorphisms of
$E$. In particular, we have $p$ divides $a$ and $Z(T)$ is not centralized by $S$. From the
structure of $\ov{C_G(t)}$, $T \le O_2(C_H(t)) \le O_{p'}(C_H(t)) = O_{p'}(C_G(t)) $. Suppose that
$T^g \not=T$ for some $g \in C_G(t)$. Then $T^g \le O_2(C_H(t))$  and $T^gE/E$ is a non-trivial
$2$-group of outer automorphisms of $E$. It follows that $T^gE/E$ is cyclic and acts non-trivially
on the cyclic group $C_E(t)/T$. However $T^g \le O_2(C_H(t))$ and so this is impossible. Thus $T$
is normalized by $C_G(t)$. Since $Z(T)$ is centralized by $C_E(t)$, $L=\langle
C_E(t)^{C_G(t)}\rangle $ centralizes $Z(T)$. Since $p$ divides $|C_E(t)|$, the structure of
$\ov{C_G(t)}$ shows that $\ov L \ge O^{p}(\ov{C_G(t)})$.  Suppose that $2^a \not= 8$. Then, as
$|C_E(t)/T|=(2^{a}+1)/(2^a+1,3)$, there is some element $\omega \in C_E(t)$, $o(\omega) = r$, $r$ a
Zsygmondi prime dividing $2^{a}+1$. Note that $r$ does not divide $|\ov{C_G(t)}|$ and that $\omega$
acts irreducibly on $T/Z(T)$.  Thus $C_G(t)$ acts on $T/Z(T)$ irreducibly. Since $|T|= 2^{3a}$ with
$|Z(T)|= 2^a$, using \cite[Lemma 2.7.3]{PR} we have that $C_G(t)/C_G(T)T$ embeds into
$\SL_2(2^a)\colon a$. Since $C_G(T)\le O_{p'}(C_G(t)) $ we have that $C_G(t)/C_G(T)T$ involves
 $\SL_2(8)\colon 3$ when $p=3$ and involves
${}^2\B_2(32)\colon 5$ when $p=5$  and these groups  must be sections of $\SL_2(2^a):a$.
Furthermore, $r$ divides $O_{p'}(C_G(t)/C_G(T)T)$ which is thus non-trivial.  The structure of
$\SL_2(2^a)$ now implies that $C_G(t)/C_G(T)T$ is soluble, a contradiction. So we are left with
 $2^a = 8$. Let $U$ be a hyperplane of $Z(T)$. Then $T/U$ is extraspecial of order $2^7$. Since
 $L$ centralizes $Z(T)$, $L$ acts on $T/U$. Furthermore, $L/C_L(T)$ has a section isomorphic to
 $\SL_2(8)$ when $p=3$ or to ${}^2\B_2(32)$ when $p=5$. Since $\Out(T/U)$ is isomorphic to  $\mathrm O_6^\pm(2)$, this is also impossible.
 Hence we have shown that $E$ is not a rank 1 group defined in characteristic 2.

Assume that $E/Z(E)$ is $\PSL_2(r^a)$ with $r$ odd. Since $\PSL_2(9) \cong \Alt(6)$,
Lemma~\ref{QS2} implies that $r^a\ne 9$.  Thus the Schur multiplier of $E$ has order $2$ and so  $E
\cong \PSL_2(r^a)$. Let $t$ be an involution in $E$. Then $C_E(t)$ is a dihedral group and hence we
have $|C_H(t)|_p>|C_E(t)|_p$ and Lemma~\ref{ppart} implies that either $p=3$ and $\ov{C_G(t)}\cong
\PSL_2(8):3$ or $p=5$ and $\ov{C_G(t)} \cong {}^2\B_2(32)$. Let $S \in \Syl_p(C_H(t))$. Then $S$ is
extraspecial of order $p^3$. Since the Sylow $p$-subgroups of $\Out(E)$ are cyclic, we must have
$S\cap E$ is cyclic of order $p^2$. In particular, as $C_E(t)$ is a dihedral group, we have
$Q=O_p(C_E(t))>1$. Thus, because $O_{p'}(C_G(t))=O_{p'}(C_H(t))$,
$[Q,O_{p'}(C_G(t))]=[Q,O_{p'}(C_H(t))]=1$. Therefore $L=\langle Q^{C_G(t)}\rangle $ centralizes
$O_{p'}(C_G(t))$ and so, as the Schur multipliers of $\PSL_2(8)$ and ${}^2\B_2(32)$ are both
trivial and $Q$ is inverted in $N_{C_G(t)}(Q)$, we have that $L$ is a normal component of $C_G(t)$.
Thus $C_G(t) = (L \times O_{p'}(C_H(t)))S$. Let $T \in \Syl_2(C_H(t))$. Then as $E$ has one
conjugacy class of involutions, $T \in \syl_2(H)$  and $T$ normalizes $Q$. Thus we can write $T =
(T \cap L) \times (T \cap O_{p'}(C_G(t)))$. Let $T_L = T\cap L$. Then $T_L$ is cyclic of order
$p-1$ ($2$  or $4$) and $T_L\le Z(T)$. In particular, $Z(T) \ge \langle t \rangle T_L$, $t \not \in
T_L$ and $|Z(T)|\ge 4$. Assume that $E$ has non-abelian Sylow $2$-subgroups. Then $T \not= \langle
t \rangle T_L$ and  $T_L \cap E=1$ for otherwise $|Z(T) \cap E| \ge 4$. Let $f \in
\Omega_1(T_L)^\sharp$. Then $f$ centralizes a Sylow $2$-subgroup of $E$ and so we have that $f$
induces a non-trivial a field automorphism of $E$.  Thus, since $r^a\not =3^2$, $C_H(f)$ contains a
component $F$ isomorphic to $\PSL_2(r^{a/2})$ and $a$ is even. As $a$ is even, we have that
$C_E(t)$ is of order $r^{a} - 1$ \cite[II 8.27]{Huppert}. Hence we have that $p$ divides $r^{a}-1$.
In particular, $p$ divides $|\PSL_2(r^{a/2})|= r^{a/2}(r^a-1)/2$. Hence $\PSL_2(r^{a/2})$ is
isomorphic to a subgroup of $C_H(f)/O_{p'}(C_G(f))$ which is consequently not soluble.  Now
Lemma~\ref{QS1} delivers a contradiction. It follows that $E$ has abelian Sylow $2$-subgroups and,
as $|Z(T)|\ge 4$,  $H=E \cong \PSL_2(r)$ where $r\equiv 3,5 \pmod 8$. In particular, $T= \langle t
\rangle T_{L}$. Therefore $T \le E$ and $T_L$ has order $2=p-1$. Therefore  $p=3$ and   $L \cong
\PSL_2(8)$. Let $R \in \syl_2(C_G(t))$ such that $R \ge T$. Then $R$ is elementary abelian of order
$16$. Since $R \in \syl_2(C_G(t))$, $R \in \syl_2(C_G(R))$. Note that $N_L(R)/R$ is cyclic of order
$7$. Thus $N_{L}(R)$ induces orbits of length $1$, $7$ and $7$ on the non-trivial elements of $R$.
Since the non-trivial elements of $T$ are all conjugate, $T_L \le L$, and $te\not\in L$ where $e
\in T_L$, we infer that all the involutions in $R$ are conjugate to $t$. In particular, $R \in
\Syl_2(C_G(r))$ for all $r \in R^\sharp$. Therefore $N_G(R)$ acts transitively on $R^\sharp$ and
$N_G(R)/C_G(R)$ is a subgroup of $\GL_4(2)$ divisible by $3.5.7$ and of odd order. There are no
such subgroups in $\GL_4(2)$ and so we have our final contradiction to this configuration. Hence
$E/Z(E) \not\cong \PSL_2(r^a)$ for odd $r$.

So we are left with $E/Z(E) \cong \PSU_3(r^{a})$ or ${}^2\G_2(r^{a})$, where $r =3$, $a > 1$ in the
latter. As in $\Aut(\PSU_3(3))$, no involution centralizes a group of order $p^2$, we get that
$r^{a} \not= 3$. Therefore in all the cases, there is an involution $t$ such that $E(C_E(t)/\langle
t \rangle) \cong \PSL_2(r^{a})$ is a non-abelian simple group. If $p$ divides $|\PSL_2(r^{a})|$,
then  Lemmas~\ref{QS1} gives a contradiction. So we have that $p$ does not divide
$|\PSL_2(r^{a})|$. As $C_E(t) = \GU_2(r^{a})$ or $\langle t \rangle \times \PSL_2(r^{a})$, we get
that $C_E(t)$ is a $p'$-group and this contradicts Lemma~\ref{QS3}.
\end{proof}

\begin{lemma}\label{QS4}  $E/Z(E) \not \cong \PSL_3(2^a)$ for $a\ge 1$.
\end{lemma}

\begin{proof}
Suppose that $E/Z(E) \cong \PSL_3(2^a)$, $a \ge 1$. Then $E$ has exactly one conjugacy class of
involutions. If $a=1$, then $C_E(t)$ is dihedral of order $8$ and $\Out(E)$ has order $2$, so $H$
does not satisfy our Hypothesis~\ref{hypH}(ii). So we may assume that $a\ge 2$.

Let $t \in E$ be an involution. Then $T = O_2(C_E(t)) \in \syl_2(C_E(t))$. Furthermore, a
complement to $T$ in $C_E(t)$ is cyclic of order dividing $2^a-1$. Hence, again by
Hypothesis~\ref{hypH} (ii), $|C_H(t)|_p
> |C_E(t)|_p$. Thus Lemma~\ref{ppart} implies that $p \in \{3,5\}$ and
$\ov{C_G(t)}$ has extraspecial Sylow $p$-subgroups. Since $\Out(E)$ has abelian Sylow
$p$-subgroups, we have $|C_E(t)|_p>1$ and, in particular, $p$ divides $2^a-1$.

We have that $O_2(C_H(t)/O_{p'}(C_G(t)) \cap F^\ast(C_G(t)/O_{p'}(C_G(t)))) = 1$.  In particular, $T$ is
normal in $O_{p'}(C_G(t))$ and, using Lemma~\ref{AutSL3Sylow}, we now see that $T$ is normal
in $C_G(t)$.

By considering $T$ as the subgroup of lower unitriangular matrices in $\SL_3(2^a)$, we see that $T$
contains exactly two elementary abelian subgroups $F_1$ and $F_2$ of order $2^{2a}$ (all the
elements of $T$ outside these two subgroups have order $4$). Thus, as $T$ is normal in $C_G(t)$,
$C_G(t)$ permutes $F_1$ and $F_2$
 and consequently $C_G(t)$ has a subgroup of index at most $2$ which normalizes $F_1$.
It follows that $N_G(F_1)\not \le H$ and that $m_p(N_G(F_1))=m_p(C_G(t)) \ge 2$. By
Lemmas~\ref{main1} and \ref{main2} (ii),  we have that $L = F^\ast(N_G(F_1)/O_{p'}(N_G(F_1)))$ is a
non-abelian simple group. Since $N_E(F_1)/O_2(N_E(F_1))$ has a section isomorphic to $\SL_2(2^a)$
with $a\ge 2$, and since $p$ divides $|\SL_2(2^a)| = (2^a-1)2^a(2^a+1)$, we get with
Proposition~\ref{Hstru1} that $F^*(N_G(F_1)/O_{p'}(N_G(F_1))) \cong \Alt(2p)$ or $\Fi_{22}$.
Furthermore, we have $N_H(F_1)/O_{p'}(N_G(F_1)) $ is isomorphic to $(\Alt(p) \times \Alt(p)):2$
with no normal components or $\Omega_8^+(2)$, neither of these has a normal subgroup isomorphic to
$\PSL_2(2^{a})$ and so we have a contradiction.
\end{proof}

The next two lemmas are needed to finally dispatch the Lie type groups as possibilities for $E$.

\begin{lemma}\label{rootinvolution} Suppose that $K$ is a simple group of Lie type defined in characteristic $2$ and $p$
is an odd prime. Let $t \in K$ be an involution in a long root subgroup of $K$. If $p$ divides
$|C_K(t)|$, then either
\begin{enumerate}
\item $p$ divides $|O^{2'}(C_K(t))|$; or
\item $K$ is isomorphic to one of $\PSL_2(2^a)$, $\PSU_3(2^a)$,
${}^2\B_2(2^{2a+1})$ or $\PSL_3(2^a)$ for some $a \ge 1$.
\end{enumerate}
\end{lemma}

\begin{proof} As $C_K(t) = O^{2'}(C_K(t))B$, where $B$ is contained in a Borel subgroup of $K$, we may assume that $p$ divides
the order of a Borel subgroup. If the Lie rank of $K$ is at least two and $K \not\cong
\PSL_3(2^{a})$, we have that $t$ is centralized by a minimal parabolic subgroup and so $p$ divides
the order of this minimal parabolic. This is the assertion.
\end{proof}

\begin{lemma}\label{classicalinvolution} Suppose that $K$ is a simple group of Lie type defined in characteristic $r$, $r$ odd, and of Lie rank at least $2$.
Let $t$ be an involution in a fundamental subgroup of $K$ and assume that $p$ is an odd prime. If
$p$ divides $|C_K(t)|$, then $p$ divides $|O^{2'}(C_K(t))|$.
\end{lemma}

\begin{proof} We have that $C_K(t) = L_1L_2S$, where $L_1$ is the fundamental subgroup, $L_2 = C_{C_K(t)}(L_1)$ and $S$
is a Sylow $2$-subgroup with $t \in Z(S)$. If $L_1L_2 = O^{2'}(L_1L_2)$, we are done. So we may
assume that $r = 3$. If $L_2 = O^{2'}(L_2)$ and $L_2 \not= 1$, then, as $3$ divides the order of
$L_2$, we also are done. Hence we have $p = 3$ and $L_2$ is a central product of groups isomorphic
to $\SL_2(3)$ or $\PSL_2(3)$. This now shows that $K $ is isomorphic to one of the following
$\PSL_3(3)$, $\PSp_4(3)$, $\mathrm P \Omega_6^\pm(3)$, $\Omega_7(3)$, $\mathrm P\Omega_8^+(3)$ or
$\mathrm G_2(3)$. But in all cases we have some $\GL_2(3)$ involved, so $3$ divides
$|O^{2'}(C_K(t))|$.
\end{proof}

For $X$ a Lie type group of Lie rank at least two in odd characteristic $r$ we list
$O^{r'}(C_X(t))$ for $t$ a classical involution in Table~\ref{centclassic}. The information is
taken from \cite[Table 4.5.1, Theorem 4.5.5]{GLS3}.

\begin{sidewaystable}

\begin{tabular}{cc||cc}
Simple group& Levi Factor of &Simple group& Levi Factor\\
\hline
$\PSL_n(q)$&$\SL_{n-2}(q)$&$\mathrm{PSU}_n(q)$&$\mathrm{SU}_{n-2}(q)$\\
$\PSp_{2n}(q)'$&$\mathrm{Sp}_{2(n-1)}(q)$&$\Omega^{\pm}_{n}(q)$&$\SL_2(q)*
\mathrm{SO}^{\pm}_{n-4}(q)$\\
$\mathrm G_2(q)'$&$\SL_2(q)$&${}^3\mathrm D_4(q)$&$\SL_2(q^3)$\\
$\F_4(q)$&$\mathrm{Sp}_6(q)$&${}^2\F_4(q)$&${}^2\mathrm{B}_2(q)$\\
${}^2\E_6(q)$&$\SU_6(q)$& $\E_6(q)$&$\SL_6(q)/\Z_{(q-1,3)}$\\$\E_7(q)$&$\mathrm{SO}_{12}^+(q)$&
$\E_8(q)$&$\E_7(q)$
\end{tabular}
\caption{Subgroups generated by root elements in the Levi complements of centralizers of
involutions in the centre of long root subgroups in Lie type groups of rank at least $2$ defined
over a field of order $q=2^a$.}\label{centlong}

\begin{tabular}{cc||cc}
Simple group $X$&$O^{r'}(C_X(t))$&Simple group $X$&$O^{r'}(C_X(x))$\\ \hline
$\PSL_n(r^{a})$&$\SL_2(r^{a}) \ast \SL_{n-2}(r^{a})$&$\mathrm{PSU}_n(r^{a})$&$\SL_2(r^{a}) \ast \mathrm{SU}_{n-2}(r^{a})$\\
$\PSp_{2n}(r^{a})$&$\SL_2(r^{a}) \ast
\mathrm{Sp}_{2(n-1)}(r^{a})$&$\Omega^{\pm}_{n}(r^{a})$&$(\SL_2(r^{a}) \ast \SL_2(r^{a}))*
\mathrm{SO}_{n-4}^\pm(r^{a})$\\
$\mathrm G_2(r^{a})$&$\SL_2(r^{a})\ast \SL_2(r^{a})$&${}^3\mathrm D_4(r^{a})$&$\SL_2(r^a)*\SL_2(r^{3a})$\\
$\F_4(r^{a})$&$\SL_2(r^{a}) \ast \mathrm{Sp}_6(r^{a})$&${}^2\E_6(r^{a})$&$\SL_2(r^{a}) \ast \SU_6(r^{a})$\\
$\E_6(r^{a})$&$\SL_2(r^{a}) \ast \SL_6(r^{a})$&$\E_7(r^{a})$&$\SL_2(r^{a}) \times \Omega_{12}^+(r^{a})$\\
$\E_8(r^{a})$&$\SL_2(r^{a}) \times \E_7(r^{a})$
\end{tabular}
\caption{The group generated by the $r$-elements in the centralizer of a classical involution $t$
in Lie type groups of rank at least $2$ and odd characteristic $r$.}\label{centclassic}

\end{sidewaystable}

\begin{lemma}\label{QS6} $E/Z(E)$ is not a simple Lie type group defined
in characteristic $2$.
\end{lemma}

\begin{proof}
Suppose that $E/Z(E)$ is a Lie type group defined in characteristic $2$. By Lemma~\ref{QS5} we have
that the Lie rank of $E$ is at least two. Further $E/Z(E) \not\cong \PSL_3(2^a)$, $a \ge 1$ by
Lemma~\ref{QS4}. Let $S \in \Syl_2(E)$ and let $t$ be an involution in the centre of the long root
group $X_\rho$ contained in $Z(S)$. Then $C_E(t)$ is a subgroup of $N_G(X_\rho)$ and $t$ is
centralized by the subgroup $K_t$ of the Levi complement of $N_G(X_\rho)$ which is generated by
root subgroups. These subgroups are given in Table~\ref{centlong} and we note that $\ov{K_t}$ is a
characteristic subgroups of $\ov{C_E(t)}$ which in turn is normal in $\ov{C_H(t)}$. By
Lemmas~\ref{rootinvolution} and \ref{QS3}, we have that $p$ divides $|O^{2'}(C_E(t))|$. Suppose
that $K_t/O_{p'}(K_t)$ is not soluble. Then, by Lemma~\ref{QS1}, we have that $F^*(\ov{C_H(t)})
\cong \Alt(p)\times \Alt(p)$ with $p \ge 5$ or $\Omega_8^+(2)$ with $p=5$. Since $\ov {K_t}$ is
normal in $\ov{C_H(t)}$, we infer that either $\ov{K_t} \cong \Alt(p) \times \Alt(p)$ and $p \ge 5$
or to $\Omega^+_8(2)$ and $p=5$. As $\Alt(p)$, $p \ge 5$, is isomorphic to a Lie type group defined
in characteristic $2$ if and only if $p=5$, we now have that $p=5$ and $K_t \cong \Alt(5) \times
\Alt(5)$ or to $\Omega_8^+(2)$. Using Table~\ref{centlong} we get that $E/Z(E) \cong
\Omega_{12}^+(2)$. Furthermore, we  have $Z(E) =O_2(E)=1$.  Therefore there is a further involution
$s \in E$ which has centralizer $C$ with $C/O_2(C)\cong \Sp_8(2)$. Now using Lemma~\ref{QS1},  we
have a contradiction. Hence we must assume that $C_E(t)$ is soluble. By Lemmas~\ref{QS2} and \ref{QS5},
$E/Z(E) \not \cong \PSp_4(2)'\cong \Alt(6)$ or to $\G_2(2)' \cong \PSU_3(3)$. Then using
Table~\ref{centlong} again, we have  $E/Z(E)$ is one of  $\PSU_4(2)$, $\PSU_5(2)$,
 ${}^2\F_4(2)'$, $\PSL_4(2)$ and $\Omega_8^+(2)$. In all the cases we check that $|C_H(t)|$ is only divisible by $p^2$ when $p=3$ and $E/Z(E) \cong \PSU_4(2)$,
$\PSU_5(2)$ or $\Omega_8^+(2)$. Suppose that  $E/Z(E) \cong \Omega_8^+(2)$. Then $E \cong
\Omega_8^+(2)$ and there is an involution $s \in E$ such that $C_E(s)/O_{2}(C_E(s))$ has a normal
subgroup isomorphic to $\Sp_4(2)'\cong \Alt(6)$. This violates Lemma~\ref{QS1}.  If $E \cong
\PSU_4(2)$,
 then there is an involution $s \in E$, which is not centralized by an elementary abelian group of order 9 contrary to Hypothesis~\ref{hypH} (ii).
 So we have $E \cong \PSU_5(2)$. In this case $\ov{C_H(t)} \cong 3^{1+2}_+.\SL_2(3)$ or $\ov{C_H(t)}\cong
 3^{1+2}_+.\GL_2(3)$.  Now Proposition~\ref{Hstru1} shows that this structure is impossible. Thus
 $E/Z(E)$ is not a Lie type group defined in characteristic $2$.
\end{proof}

\begin{lemma}\label{QS7} $E/Z(E)$ is not a Lie type group in
odd characteristic.
\end{lemma}

\begin{proof} Suppose that $E/Z(E)$ is a Lie type group defined in characteristic $r$, $r$ an odd prime.
By Lemma~\ref{QS5}, we may assume that the Lie rank of $E$ is at least two. Let $t$ be classical
involution in $t$ in $E$. By Lemmas~\ref{QS3} and \ref{classicalinvolution}, $p$ divides
$|O^{2'}(C_E(t))|$. Let $K_1$ be a subnormal subgroup of $C_E(t)$ containing $t$ with $K_1 \cong
\SL_2(r^a)$. Assume that $p$ divides $|K_1|$ and that $r^a \neq 3$. Then Lemma~\ref{QS1} implies
that $\ov K_1 \cong \Alt(p)$ or $\Omega_8^+(2)$. We conclude that $p= r^{a} = 5$ and that
$F^*(\ov{C_G(t)}) \cong \Alt(10)$. In particular, $K_1$ is not normal in $C_H(t)$. Using
Table~\ref{centclassic} now shows that $E \cong \PSL_4(5)$, $\PSp_4(5)$, $\PSU_4(5)$, or $\G_2(5)$.
The last case immediately fails, as $\Out(\G_2(5))=1$ means  $E = H$ and $K_1$ is normal in
$C_H(t)$. If $E \cong \PSp_4(5)$, then there is an involution $s\in  E$ with $O^{5'}(C_E(s)) \cong
\PSL_2(5)$, which contradicts Hypothesis~\ref{hypH} (ii). If $E \cong \PSU_4(5)$, then there is an
involution $s \in E$ with $F^*(C_E(x)/O_{5'}(C_E(x))) \cong \PSL_2(25)$ and this contradicts
Lemma~\ref{QS1}. So  we must have that $E \cong \PSL_4(5)$. Since $O_{5'}(C_G(t)) \le H$ is
normalized by $C_E(t)$, the structure of $\Aut(E)$ shows that $O_{5'}(C_G(t))\cap E = \langle
t\rangle$. It follows that $C_G(t)/O_{p'}(C_G(t))$ has a subgroup isomorphic to $C_E(t)/\langle t
\rangle$. Since the subgroup $M$ of even permutations of $\Sym(5)\wr \Sym(2)$ has $M/F^*(M)$ cyclic
of order $4$ and $C_E(t)/F^*(C_E(t))$ is a fours group, we see that $C_G(t)/O_{p'}(C_G(t)) \cong
\Sym(10)$. But then $H
> E$. In particular, there are involutions $s \in H \setminus E$ with $F^*(C_H(s)) \cong \PSL_3(5)$
in contradiction to Lemma~\ref{QS1}.

So we may assume next that $p$ does not divide the order of $\SL_2(r^{a})$ or $r^{a} = 3$. Suppose
that in the latter case we also have $p \not= 3$. Then in both cases $p $ does not divide $r^a\pm
1$. Hence $p$ divides the order of one of the factors listed in
 Table~\ref{centclassic}. Since $p$ does not divide $|K_1|$, $C_E(t)$ has a component which is a group of  Lie type in
characteristic $r$. But these components are not isomorphic to either  $\Alt(p)$ or
$\Omega_8^+(2)$, and this contradicts Lemmas~\ref{Ctnot-in} and \ref{QS1}.

So we have that $p=3=r^{a}$. Suppose that $C_E(t)$ has a component $F$, then $3$ divides $|F|$.
Therefore  Lemmas~\ref{Ctnot-in}, \ref{QS1} and Table~\ref{centclassic} give a contradiction. So we
have that $C_E(t)$ is soluble. Then reference once again to Table~\ref{centclassic} gives $E/Z(E)$
is isomorphic to one of $\PSL_3(3)$, $\PSL_4(3)$, $\PSp_4(3)$, $\PSU_4(3)$, $\Omega_7(3)$, $\mathrm
P\Omega_8^+(3)$, $\G_2(3)$. If $E \cong \Omega_7(3)$ or $\mathrm P\Omega_8^+(3)$, there is an
involution $s$ with $F^*(C_E(s)) \cong \mathrm{SO}_6^-(3)$ and, if $E \cong \PSL_4(3)$ there is an
involution $s\in E$ $F^*(C_E(s)) \cong \PSL_2(9)$. Thus in these cases we obtain a contradiction
via Lemma~\ref{QS1}. If $E \cong \PSL_3(3)$ or $\PSp_4(3)$, then
 there is an involution $s\in E$ whose centralizer is not divisible by
$9$, in the first case $|C_E(s)| = 2^5\cdot 3$, while in the second case $s = t$ and $C_E(t) \cong
\GL_2(3)$. This contradicts Hypothesis~\ref{hypH} (ii). So we have that $E/Z(E) \cong \PSU_4(3)$ or
$\G_2(3)$. Assume that $H = E$. Then, by Lemma~\ref{main1} (v),  $O_{3'}(C_G(t)) \le H$ and
Corollary~\ref{OrredH3} gives $O_2(C_E(t)) =O_{3'}(C_G(t))$ which is extraspecial $2^{1+4}_+$. It
follows that $C_G(t)/C_{C_G(t)}(O_2(C_E(t))$ is soluble. But $\ov{C_G(t)}$ is an almost simple
group and so we have that $O_2(C_E(t))$ is centralized by a Sylow $3$-subgroup of $C_E(t)$
something which is impossible. Thus $E \neq H$.
 Hence there exists an involution $s\in H\setminus E$. If $E/Z(E) \cong
\G_2(3)$, $C_E(s) \cong {}^2\G_2(3)$. If $E/Z(E) \cong \PSU_4(3)$, then $C_E(s)$ has a section
isomorphic to $\PSU_3(3)$, $\PSp_4(3)$ or $\Omega_4^-(3)$.  In all these cases we get a
contradiction to Lemma \ref{QS1} applied to $C_G(s)$.
\end{proof}

\section{Proofs of the Main Theorems}\label{lastsec}

In this final section we  assemble the proofs of our main theorems. We refer the reader to the
introduction for their statements and continue the notation of the previous sections.

\begin{proof}[Proof of Theorem~\ref{MainTheorem0}] Let $G_0 = O^{2}(G)$.
Then $H_0=H \cap G_0$ has even order and $H_0$ is strongly $p$-embedded in $G_0$ by
Lemma~\ref{main1}(ii) and (iv). As $H_0$ is normal in $G_0$ and as $H_0$ has even order by
hypothesis, we have that $G_0$ and $H_0$ together satisfy  Hypothesis~\ref{hypH}. If $G_0$ contains
a classical involution, then Theorem~\ref{clasinv} implies that $G_0$ is a $\mathcal K$-group and
then Proposition~\ref{SE-p2} implies that $F^*(G)= F^*(G_0)\cong \PSU_3(p^a)$ for some $a \ge 2$.
Thus we may suppose that $G_0$ has no classical involutions. Since $F^*(G) = O_p(G)$,
Hypothesis~\ref{EH=1} is satisfied. Therefore Theorems~\ref{L2pthm} and \ref{Ree} together imply
that $F^*(G)= F^*(G_0) \cong {}^2\G_2(3^{2n-1})$ for some $n \ge 2$.
\end{proof}

\begin{proof}[Proof of Theorem~\ref{MainTheorem1}] Suppose that $F^*(H)\neq O_p(H)$. Then as
$O_{p'}(H)=1$, we have $E=E(H) \neq 1$. Let $G_0 = O^2(G)$. Then $G_0 \ge E(H)$ and so $G_0$
satisfies hypothesis~\ref{hypH}. Therefore Lemma~\ref{Eqs} implies that $E(H)$ is a quasisimple
group. Again we may as well assume that $G$ does not contain a classical involution. In particular,
Hypothesis~\ref{EHnot1} is satisfied. By Lemma~\ref{Op=1}, we have $O_p(H)= O_{p'}(H)=1$. Now we
use the fact that $E$ is a $\mathcal K$-group and use Lemmas~\ref{QS2}, \ref{QSsporadic}, \ref{QS6}
and \ref{QS7} to deliver a contradiction. Hence we conclude that $F^*(H)= O_p(H)$.
\end{proof}

\begin{proof}[Proof of Theorem~\ref{MainTheorem}] This is merely a combination of
Theorems~\ref{MainTheorem0} and \ref{MainTheorem1}.
\end{proof}

Finally we prove Corollary~\ref{LiepCor}. For this we require the following  proposition about
centralizers of involutions in Lie type groups defined in characteristic $p$.

\begin{proposition}\label{invs5} Let $p$ be an odd prime, $X$ be a finite group and set $K= F^*(X)$.
Suppose that $K$ is a simple group of Lie type defined in characteristic $p$. Then either
$m_p(C_X(x)) \ge 2$ for all involutions $x \in X$ or one of the following holds:
\begin{enumerate}
\item $K \cong \PSL_2(p^n)$ with $n \ge 1$;
\item $K \cong \PSL_3(p)$ or $\PSU_3(p)$;
\item $K\cong \PSp_4(p)$; or
\item $K\cong {}^2\mathrm G_2(3)'\cong \PSL_2(8)$.
\end{enumerate}
\end{proposition}

\begin{proof} We may suppose that $K \not \cong \PSL_2(p^n)$. Let $x \in X$ be an involution and  set $K_x = O^{p'}(C_{K}(x))$.
Note that since the Lie type groups are generated by their (perhaps twisted) root subgroups the
only Lie type groups $L$ with $m_p(L) =1$ are $L\cong \SL_2(p)$ for arbitrary $p$ and $L \cong
{}^2\mathrm G_2(3)' \cong \PSL_2(8)$ with $p=3$. Assume that $q= p^n$ and $K$ has type
${}^d\Sigma(q)$ (using notation as in \cite[Chapter 4]{GLS3}). Then, by \cite[Definition
2.5.13]{GLS3}, $x$ induces either an inner-diagonal, graph, field or graph-field automorphism on
$K$. If $x$ induces a field automorphism, then \cite[Proposition 4.9.1]{GLS3} gives that $K_x \cong
{}^d\Sigma(q^{1/2})$. Since $m_p(C_X(x)) \le 1$, we deduce that $K_x \cong \PSL_2(p)$ and
consequently $K \cong \PSL_2(p^2)$ which is a contradiction. If $x$ induces a graph-field
automorphism of $K$, then \cite[Proposition 4.9.1]{GLS3} implies $K_x \cong {}^2\Sigma(q^{1/2})$.
Thus ${}^2\Sigma(q^{1/2})$ is $ \PSL_2(p)$ or $\SL_2(p)$ and this is impossible.  Therefore $x$ is
contained in the group generated by inner-diagonal and graph automorphisms of $K$. For each of the
possibilities for $K$ and $x$, the structure of $K_x$ is described in \cite[Table 4.5.1]{GLS3} and
we avail ourselves of this information. Thus, if $K$ has type $\A_n^\epsilon(q)$, then the only
possibility is that $n=2$ and $q=p$ and so (ii) holds. Similarly if $K$ has type $\B_m(q)$ or
$\C_m(q)$, we see again that $m=2$ and $q=p$ and thus (iii) holds. The groups of type
$\D_m^\epsilon (q)$ only need to be considered for $m \ge 4$, and in these cases $m_p(K_x) \ge 2$.
The same is true from the groups of type $\mathrm G_2(q)$ and ${}^3\mathrm D_4(q)$. For $K\cong
{}^2\G_2(3^{2n+1})'$, the only possibility is that $X \cong {}^2\G_2(3)$ and $K_x \cong \PSL_2(3)$
which gives (iv). The remaining exceptional groups are all easily seen to violate $m_p(K_x) \le 1$.
\end{proof}

\begin{proof}[Proof of Corollary~\ref{LiepCor}] If $F^*(H)$ is a Lie type group defined in characteristic $p$ and of
rank at least $3$, then Proposition~\ref{invs5} implies that $m_p(C_H(t)) \ge 2$ for all
involutions $t \in H$. Thus  Theorem~\ref{MainTheorem1} implies that $H$ cannot be strongly
$p$-embedded in any groups $G$.
\end{proof}

%\begin{lemma}\l{autf4} Let $L = {}^2F_4(q)^\prime$, $q$ even. Then $L$ has no outer automorphism of order two.
%\end{lemma}
%
%\begin{proof} If $q > 2$, this follows from \cite[(9.1)]{GLS4}. So assume that $q = 2$, then we claim that Aut$(L) = {}^2F_4(2)$.
%For the structure of $L$ see \cite{Shi}. For that let $P$ be the parabolic in $L$ with $P/O_2(P) \cong F_{20}$. If $|Out(L)|$ is divisible by 4, we may assume
%that there is a group of order 4, which centralizes $P/O_2(P)$. Hence it centralizes an element $\omega$ of order 5.
%Now we have that $O_2(P)/Z_2(O_2(P))$ is of order 16 and so we have a group $Q$ containing $Z_2(O_2(F))$, which covers this group of order 4.
%In particular we have that $C_{Q}(\omega)$ is of order 8. Further of course $[Z_2(O_2(P)), C_{Q}(\omega)] = 1$.
%As $[Q,O_2(P)] \leq Z_2(Q)$, we see that there is some $x \in C_Q(\omega)$, which centralizes $O_2(H)$ and induces an outer automorphism on $L$.
%Hence $x$ centralizes $\Omega_1(S)$, for a Sylow 2--subgroup $S$ of $P$. Let now $R$ be the other parabolic in $L$
%containing $S$. Then we see that $x$ centralizes $\Omega_1(O_2(R))$. On $Z(\Omega_1(O_2(R)))$ we have that $R$ induces a natural module
%and just trivial modules. So we see that we may assume that $x$ is centralized by an element $\nu \in R$ of order 3.
%Hence we have that $P = SC_P(x)$ and $R = SC_R(x)$. But $L$ acts flag transitively on the generalized 8-gon for $L$ and so $x$ acts trivially.
%Hence we have shown that ${}^2F_4(2)$ is the automorphism group of $L$ and now the assertion follows with \cite[Corollary 2]{Shi}. \qed
%\end{proof}

\end{document}